\documentclass[12pt]{article}  
\usepackage{graphicx,amsmath,amsthm,amssymb,dsfont}  
\numberwithin{equation}{section}
\def\ZZ{\mathbb Z}
\def\PP{\mathbb P}

\def\NN{\mathbb N}

\def\FF{\mathbb F}

\def\cD{\mathcal D} 

\def\cC{\mathcal C}

\def\fB{\mathfrak B}
\def\cL{\mathcal L}

\def\cM{\mathcal M}

\def\cF{\mathcal F}                                                              
\def\ff{\mathfrak f}

\def\fc{\mathfrak c}

\def\cC{\mathcal C}

\def\cP{\mathcal P}

\def\cD{\mathcal D}

\def\cN{\mathcal N}

\def\sk{{\sf k}}

\def\bA{{\bf A}}

\def\bT{{\bf T}}

\def\bx{{\bf x}}
\def\by{{\bf y}}
\def\bz{{\bf z}}
\def\bv{{\bf v}}
\def\bw{{\bf w}}
\def\bn{{\bf n}}

\def\bh{{\bf h}}

\def\ba{{\bf a}}

\def\bc{{\bf c}}
\def\bd{{\bf d}}

\def\bff{{\bf f}}

\def\bgamma{\boldsymbol{\gamma}}

\def\bvarphi{\boldsymbol{\varphi}}
\def\bpsi{\boldsymbol{\psi}}
\def\bff{\boldsymbol{\ff}}

\def\bs{{\bf 0}}
\def\b1{{\bf 1}}
\def\d1{\mathds{ 1}} 
\def\deg{{\rm deg}}
\def\div{{\rm div}}
\def\Div{{\rm Div}}
\def\Res{{\rm Res}}
\def\Tr{{\rm Tr}}
\def\dd{{\rm d}}
\def\mod{{\rm mod}}
\def\ad{{\rm and}}
\def\with{{\rm with}}
\def\where{{\rm where}}
\def\for{{\rm for}}
\def\fall{{\rm for\; all}}
\def\vol{{\rm vol}}
\def\Con{{\rm Con}}

\begin{document}
\hsize=14.8 true cm



\title{On the lower bound of the discrepancy  of $(t,s)$ sequences: II}
\author{Mordechay B. Levin}

\date{}

\maketitle

\begin{abstract}
Let $ (\bx(n))_{n \geq 1} $ be an $s-$dimensional Niederreiter-Xing's sequence in base $b$. Let $D((\bx(n))_{n = 1}^{N})$ be the discrepancy of the sequence $ (\bx(n))_{n = 1}^{N} $.
It is known that $N D((\bx(n))_{n = 1}^{N}) =O(\ln^s N)$ as $N \to  \infty $. 
 In this paper, we prove that this estimate  is exact. Namely, there exists a constant $K>0$, such that
 $$
    \inf_{\bw \in [0,1)^s}  \sup_{1 \leq N \leq b^m}   N D((\bx(n)\oplus \bw)_{n = 1}^{N}) \geq K m^s \quad {\rm for} \; \; m=1,2,...\;.
 $$
We also get similar results for other explicit constructions of $(t,s)$ sequences.
\end{abstract}
%
%
\section{Introduction.}
{\bf 1.1 }  Let $((\beta_{n}^{(s)})_{n \geq 1})$ be a   sequence in  unit cube $[0,1)^s$, $(\beta_{n,N}^{(s)})_{n=0}^{N-1}$ points set in $[0,1)^s$,
$J_{\by}=[0,y_1) \times \cdots \times [0,y_s) \subseteq [0,1)^s $,
\begin{equation}\label{In00}
\Delta(J_{\by}, (\beta_{n,N}^{(s)})_{k=1}^{N}  )= \#\{1\leq n
  \leq N  \;|\; \beta_{n,N}^{(s)}\in J_{\by}\}-Ny_1 \ldots y_s.
\end{equation}
We define the star {\sf discrepancy} of a 
  $(\beta_{n,N}^{(s)})_{n=0}^{N-1}$ as
\begin{equation} \label{In02}
   \emph{D}^{*}(N)=\emph{D}^{*}((\beta_{n,N}^{(s)})_{n=0}^{N-1}) = 
    \sup_{ 0<y_1, \ldots , y_s \leq 1} \; \Big| \frac{1}{  N}
  \Delta(J_{\by},(\beta_{n,N}^{(s)})_{n=1}^{N}) \Big|.
\end{equation}\\

 {\bf Definition 1.} {\it A sequence $(\beta_n^{(s)})_{n\geq 0}$ is of {\sf
low discrepancy} (abbreviated l.d.s.) if  $\emph{D}
((\beta_n^{(s)})_{n=0}^{N-1})=O(N^{-1}(\ln N)^s) $ for $ N \rightarrow \infty $.}\\

 {\bf Definition 2.} {\it A sequence  of point sets $((\beta_{n,N}^{(s)})_{n=0}^{N-1})_{N=1}^{\infty}$ is of
 low discrepancy (abbreviated
l.d.p.s.) if $ \emph{D}((\beta_{n,N}^{(s)})_{n=0}^{N-1})=O(N^{-1}(\ln
N)^{s-1}) $, for $ N \rightarrow \infty $. } \\

 For examples of such a sequence, see, e.g., [BC], [DiPi], and  [Ni].\\
In 1954, Roth proved that there exists a constant $ C_s>0 $, such
that
\begin{equation}  \nonumber
N\emph{D}^{*}((\beta_{n,N}^{(s)})_{n=0}^{N-1})>C_s(\ln N)^{\frac{s-1}{2}}, \;\;
 \quad   {\rm and} \;\; \quad \overline{\lim }{N\emph{D}^{*}((\beta_n^{(s)})_{n=0}^{N-1})(\ln
N)^{-s/2}}>0 \nonumber
\end{equation}
for all $N$-point sets $(\beta_{n,N}^{(s)})_{n=0}^{N-1}$ and all sequences
$(\beta_n^{(s)})_{n \geq 0}$.

According to the well-known conjecture (see, e.g., [BC, p.283], [DiPi, p.67], \newline[Ni, p.32]), these estimates can be improved
\begin{equation}   \label{In08}  
  {N\emph{D}^{*}((\beta_{n,N}^{(\ddot{s})})_{n=0}^{N-1}) (\ln N)^{-\ddot{s}+1}} >C_{\ddot{s}}^{'} 
     \;\;\;  {\rm and} \;\;\;
\underset{N \to \infty }{\overline{\lim }} N (\ln N)^{-\dot{s}} \emph{D}^{*}((\beta_{n}^{(\dot{s})})_{n=1}^{N})>0
\end{equation}
for all $N$-point sets $(\beta_{n,N}^{(\ddot{s})})_{n=0}^{N-1}$ and all sequences $(\beta_n^{(\dot{s})})_{n \geq 0}$ with some $C_{\ddot{s}}^{'} >0$.

In 1972, W. Schmidt proved (\ref{In08}) for $\dot{s}=1 $ and $\ddot{s}=2$. In [FaCh], 
(\ref{In08}) is proved for a class of $(t,2)-$sequences.
The review of research on this conjecture, see for example in [Bi].

A subinterval $E$ of $[0,1)^s$  of the form
$$  E = \prod_{i=1}^s [a_ib^{-d_i},(a_i+1)b^{-d_i}),   $$
 with $a_i,d_i \in Z, \; d_i \ge 0, \; 0  \le a_i < b^{d_i}$ for $1 \le i \le s$ is called an
{\it elementary interval in base $b \geq 2$}.\\

{\bf Definition 3}. {\it Let $0 \le t \le m$  be an integer. A {\sf $(t,m,s)$-{\sf net in base $b$}} is a point set 
$\bx_0,...,\bx_{b^m-1}$ in $ [0,1)^s $  such that $\# \{ n \in [0,b^m -1] | x_n \in E \}=b^t$   for every elementary interval E in base  $b$ with  
$\vol(E)=b^{t-m}$.}\\

{\bf Definition 4}. {\it 
 Let $t \geq 0$  be an integer. A sequence $\bx_0,\bx_1,...$ of points in $[0,1)^s$ is a
{\sf $(t,s)$-{\sf sequence in base} $b$} if, for all integers $k\ge 0$ and $m \geq t$, the point set
consisting of $\bx_n$ with  $ kb^m \leq  n < (k+1)b^m$      is a $(t,m,s)$-net in base $b$.}

 By [Ni, p. 56,60], $(t,m,s)$ nets and $(t,s)$ sequences  are of low discrepancy.

See reviews 
  on  $(t,m,s)$ nets and $(t,s)$ sequences  in [DiPi] and [Ni].\\

 For $x =\sum_{j \geq 1}  x_{i}p_i^{-i}$, and $y =\sum_{j \geq 1}  y_{i}p_i^{-i}$ 
where $x_{i},y_i \in  Z_b := \{0, 1,...., b-1\}$, we define the ($b$-adic) digital shifted point $v$ by 
$v = x \oplus y := \sum_{j \geq 1}  v_{i}p_i^{-i}$,
 where $v_i \equiv x_i + y_i \;\mod(b)$ and $v_i \in Z_b$.
For higher dimensions $s > 1$, let $\by = (y_1, . . . , y_s) \in [0, 1)^s$. For $\bx =
(x_1, . . . , x_s) \in [0, 1)^s$ we define the ($b$-adic) digital shifted point $\bv$ by 
$ \bv =\bx \oplus \by =(x_1 \oplus y_1, . . . ,x_s \oplus y_s)$.
 For $n_1,n_2 \in [0,b^m)$, we define 
$n_1 \oplus n_2 := (n_1 /b^m\oplus n_2)b^m)b^m$.

For $x =\sum_{j \geq 1}  x_{i}p_i^{-i}$, 
where $x_{i} \in  Z_b$,  $x_i=0 $ $(i=1,...,k)$ and $x_{k+1} \neq 0$, we define the
absolute valuation  $\left\|.  \right\|_b $ of $x$ by  $\left\|x  \right\|_b =b^{-k-1}$.
Let $\left\| n  \right\|_b =b^k$ for $n \in [b^k,b^{k+1})$.\\

	 {\bf Definition 5.} {\it A  point set   $ (\bx_{n})_{0 \leq n <b^m} $ 
	in $[0,1)^s$ is  $d-${\sf admissible} in
base $b$  if}
\begin{equation}  \label{In04}
  \min_{0 \leq k <n < b^m} \left\| \bx_n \ominus \bx_k  \right\|_b
  > b^{-m-d}  \quad {\rm where} \quad \left\| \bx  \right\|_b := \prod_{i=1}^s \left\|x^{(i)}_{j}  \right\|_b . 
\end{equation}
{\it A sequence  $ (\bx_{n})_{n \geq 0} $
	in $[0,1)^s$ is  $d-${\sf admissible} in
base $b$ if\;} $ \inf_{n >k \geq 0}
	\left\| n \ominus k 
	\right\|_b \left\| \bx_n \ominus \bx_k  \right\|_b  \\\geq b^{-d}$. \\

 Let  $ (\bx_n)_{ n \geq 0} $  be a $d-$admissible   $(t,s)$ sequence in base $b$.
 In  [Le4], we proved for all $ m \geq 9s^2(d+t)$  that
\begin{equation}  \label{In10}
    1+   \max_{1 \leq N \leq b^m} 
		N \emph{D}^{*}((\bx_{n} \oplus \bw)_{0 \leq  n < N}) \geq   
					 b^{-d}K_{d,t,s+1}^{-s} m^{s} 
\end{equation} 
with some $\bw \in [0,1)^{s}$ and $K_{d,t,s} =  4 (d+t)(s-1)^2$.

 In this paper we consider some known constructions of $(t,s)$ sequences 
 (e.g., Niederreiter's sequences, Xing-Niederreiter's sequences, Halton type $(t,s)$ 
sequences) and we prove that 
 they have $d-$admissible properties. Moreover, we prove that for these sequences 
the bound (\ref{In10}) is true for all $\bw \in [0,1)^{s}$.
This result  supports  conjecture (\ref{In08}) (see also [Be], [LaPi], [Le2] and [Le3]). 

 We describe the structure of the paper. 
	In Section 2, we fix some definitions.
  In Section 3, we state our results.
	In Section 4, we prove our outcomes.

\section{Definitions and auxiliary  results.}
{\bf 4.1 Notation and terminology for algebraic function fields.}  For
the theory of algebraic function fields, we  follow the notation and
terminology in the books  [St] and  [Sa].

 Let $b$ be an arbitrary prime power, $\sk=\FF_b$  a finite field with $b$ elements, 
$\sk(x)=\FF_b(x)$  the rational function field over $\FF_b$, and $\sk[x]=\FF_b[x]$  
  the polynomial ring over $\FF_b$.
 For $\alpha =f/g, \; f,g \in \sk[x]$, let
\begin{equation} \label{No00}
 \nu_{\infty}(\alpha) = \deg (g) - \deg (f) 
\end{equation}
be the degree valuation of $\sk(x)$. 	
	We define the field of Laurent series as
\begin{equation} \nonumber
 \sk((x)) :=  \Big\{ \sum_{i = m}^{\infty} a_i x^{i} 
       \; | \; m \in \ZZ, \; a_i \in \sk  \Big\}.
\end{equation}

 A finite extension field $F$ of $\sk(x)$ is called
 an algebraic function field  over $\sk$. Let $\sk$  is algebraically
closed in $F$. We express
this fact by simply saying that $F/\sk$  is an algebraic function field. The genus
of $F/\sk$ is denoted by $g$.

A place $\cP$ of $F$ is, by definition, the maximal ideal of some valuation
ring of $F$. We denote by $O_{\cP}$ the valuation ring corresponding to $\cP$ and we
denote by $\PP_F$ the set of places of $F$.
For a place $\cP$ of $F$, we write $\nu_{\cP}$ for the normalized discrete valuation of $F$
corresponding to $\cP$, and any element $t \in F  $ with $\nu_{\cP} (t) = 1$ is called a local parameter (prime element) at $\cP$. 

The  field $F_{\cP}:=O_{\cP}/\cP$ is called the residue field of $F$ with respect to $\cP$.
The degree of a place
$\cP$ is defined as
$\deg(\cP) = [ F_{\cP} : \sk]$.
We denote by
$\Div(F)$ the set of divisors of $F/\sk$.

 Let $y \in F\setminus \{0\}$     and denote by $Z(y)$, respectively $N(y)$, the
set of zeros, respectively poles, of $y$. Then we define the zero divisor of $y$ by
$(y)_0 = \sum_{\cP \in Z(y)}   \nu_{\cP}(y)\cP$ and the pole divisor of $y$ by
$ (y)_{\infty} = \sum_{\cP \in N(y)}   \nu_{\cP}(y)\cP$.
Furthermore, the  principal divisor of $y$ is given by
$\div(y) = (y)_0  - (y)_{\infty}$. \\

{\bf Theorem A  (Approximation Theorem).} [St, Theorem 1.3.1] {\it  Let $F/\sk$ be a function
field, $ \cP_1, . . . , \cP_n \in  \PP_F$ pairwise distinct places of $F/\sk$, $x_1, . . . , x_n \in  F$
and $r_1, . . . , r_n  \in \ZZ$. Then there is some $y \in F$ such that  }
\begin{equation} \nonumber
                  \nu_{\cP_i} (y- x_i) = r_i  \quad {\rm for} \quad  i = 1, . . . , n .  
\end{equation}

The completion of $F$ with respect to $\nu_{\cP}$ will be
denoted by $F^{(\cP)}$. Let $t$ be a local parameter of $\cP$. 
Then $F^{(\cP)}$ is isomorphic to $F_{\cP} ((t))$   (see [Sa, Theorem 2.5.20]), and
 an arbitrary element $\alpha \in F^{(P)}$ can be uniquely
expanded as  (see [Sa, p. 293])
\begin{equation} \label{No06}
 \alpha =  \sum_{i = \nu_{\cP}(\alpha)}^{\infty} S_i t^{i}   \quad  {\rm where}
 \quad  S_i = S_i(t,\alpha) \in F_{\cP} \subseteq F^{(P)}.
\end{equation}

 The derivative $\frac{\dd \alpha}{\dd t}$, or differentiation with respect to $t$, is defined by  (see [Sa, Definition 9.3.1])
\begin{equation} \label{No08}
 \frac{\dd \alpha}{\dd t} =  \sum_{i = \nu_{\cP}(\alpha)}^{\infty} iS_i t^{i-1}   .
\end{equation}

For an algebraic function field $F/\sk$, we define its set of differentials (or Hasse differentials,  H-differentials) as
\begin{equation} \nonumber
                    \Delta_F = \{y \; \dd z \; | \; y \in F, \; z \;{\rm is \;a \;separating \; element  \;for} \; F/\sk\}
\end{equation}
(see [St, Definition 4.1.7]). \\

{\bf Proposition A.} ( [St, Proposition 4.1.8] or [Sa, Theorem 9.3.13]) {\it
Let $z \in F$ be separating. Then
 every differential $\gamma \in \Delta_F$ 
can be written uniquely as  $\gamma =
y \; \dd z$ for some $y \in  F$.}

We define the order of $\alpha \;\dd \beta$ at $\cP$ by
\begin{equation} \label{No12}
         \nu_{\cP}(\alpha \; \dd \beta) : = \nu_{\cP}(\alpha \; \dd \beta/\dd t),           
\end{equation}
where $t$ is any local parameter for $\cP$  (see [Sa, Definition 9.3.8]).

 Let $\Omega_F$ be the set of all  Weil differentials of $F/\sk$. There exists a $F-$linear
isomorphism of the differential module $\Delta_F$  onto $\Omega_F$ (see [St, Theorem 4.3.2] or [Sa, Theorem 9.3.15]).

For $0 \neq  \omega \in  \Omega_F$, there exists a uniquely
determined divisor $\div(\omega) \in \Div(F)$.  Such a divisor $\div(\omega)$ is called a canonical
divisor of $F/\sk$. (see [St, Definition 1.5.11]).  For a canonical divisor $\dot{W}$, we have   (see [St, Corollary 1.5.16])
\begin{equation} \label{No14}
     \deg(\dot{W}) = 2g -2 \quad {\rm and} \quad\ell(\dot{W}) = g.          
\end{equation}
 Let $ \alpha \;\dd \beta$  be a nonzero H-differential in $F$ and let $\omega$  the corresponding
Weil differential. Then (see  [Sa, Theorem 9.3.17], [St, ref. 4.35])
\begin{equation} \label{No16}
                     \nu_{\cP} (\div(\omega) ) = \nu_{\cP}(\alpha \;\dd \beta) ,  \quad {\rm for \; all}  \quad   \cP \in \PP_F        .
\end{equation}

 Let $ \alpha \;\dd \beta$ be a  H-differential, $t$  a local parameter of $\cP$, and
\begin{equation} \nonumber
     \alpha \;\dd \beta=  \sum_{i = \nu_{\cP}(\alpha)}^{\infty} S_i t^{i} \dd t \in F^{(\cP)}.  
\end{equation}
Then the {\sf residue} of $ \alpha \;\dd \beta$ (see  [Sa, Definition 9.3.10) is defined by 
\begin{equation} \label{No18}
       \Res_{\cP} ( \alpha \; \dd \beta) := \Tr_{F_{\cP}/\sk} (S_{-1}) 
			   \in \sk .      
\end{equation}
Let
\begin{equation} \label{No19}
       \Res_{\cP,t} ( \alpha ) :=  \Res_{\cP} ( \alpha \dd t) .    
\end{equation}
 \\

{\bf Theorem B (Residue Theorem). } ([St, Corollary 4.3.3],  [Sa Theorem 9.3.14]) {\it Let $ \alpha \;\dd \beta$ be any H-differential. Then
   $\Res_{\cP} ( \alpha \; \dd \beta)=0$ for almost all places $\cP$. Furthermore,}
\begin{equation} \nonumber
           \sum_{\cP \in \PP_F}\Res_{\cP} ( \alpha \;\dd \beta)=0.        
\end{equation}

For a divisor $\cD$ of $F/\sk$, let $ \cL(\cD)$  denote the Riemann-Roch space
\begin{equation} \nonumber
           \cL(\cD) =   \cL_{F}(\cD) = \cL_{F/\sk}(\cD) = \{ y \in F\setminus {0} \; | \; \div(y) + \cD  \geq 0 \} \cup \{0\}.
\end{equation}
Then $\cL(\cD)$ is a finite-dimensional vector space over $ F$, and we denote its
dimension by $\ell(\cD)$. By [St, Corollary 1.4.12], $\ell(\cD) =\{0\}$ for $\deg(\cD) <0$.\\

{\bf Theorem C  (Riemann-Roch Theorem).}  [St, Theorem 1.5.15, and St, Theorem 1.5.17 ] {\it Let $W$ be a canonical divisor
of $F/\sk$. Then for each divisor $A \in \div(F)$, $\ell (A) = \deg(A) + 1 - g + \ell(W - A) $,  and }
\begin{equation} \nonumber
                     \ell (A) = \deg(A) + 1 - g  \quad {\rm for }  \quad  \deg(A) \geq  2 g -1.
\end{equation}  \\

Let $P \in \PP_F$, $e_P =\deg(P)$, and let $F'=FF_P$ be the {\sf compositum} field 
(see \cite[Theorem 5.4.4]{Sa}).  By \cite[Proposition 3.6.1]{St} $F_P$ is the full constant field of $F'$.

 For a place
$P \in \PP_F$, we define its {\sf conorm} (with respect to $F'/F$) as

\begin{equation} \label{No50}
             \Con_{F'/F} (P) := \sum_{P' | P} e(P'|P)P',
\end{equation}
where the sum runs over all places $P' \in \PP_{F'}$ lying over $P$
(see \cite[Definition 3.1.8.]{St}) and $e(P'|P)$ is the ramification index of $P'$ over $P$.\\

{\bf Theorem D.}  (\cite[Theorem 3.6.3]{St}) {\it 
In an algebraic constant field extension $F'=FF_P$ of $F/\sk$, 
the following hold:}
\begin{enumerate}
\item[(a)] $F'/F$ is unramified (i.e., $e(P'|P) = 1$ for all $P \in \PP_F$ and all $P' \in \PP_{F'}$
with $P'|P$).

\item[(b)]  $F'/F_P$ has the same genus as $F/\sk$.

\item[(c)]  For each divisor $A \in  \Div(F)$, we have $\deg (\Con_{F'/F} (A)) = \deg(A)$.

\item[(d)]  For each divisor $A \in  \Div(F)$, $\ell(\Con_{F'/F} (A)) = \ell(A)$.
More precisely:  Every basis of $\cL_{F/\sk}(A)$ is also a basis of 
$\cL_{F'/F_P}(\Con_{F'/F} (A))$. 
\end{enumerate} 

{\bf Theorem E. } (\cite[Proposition 3.1.9]{St}) {\it
 For $0 \neq x \in F$ let $(x)_{0}^{F}, \;(x)_{\infty}^{F}, \; \div(x)^{F} $, 
 resp. $(x)_{0}^{F'},\; (x)_{\infty}^{F'}, \;\div(x)^{F'}$
 denote
the zero, pole, principal divisor of $x$ in $\Div(F)$ resp. in $\Div(F')$. Then}
\begin{equation}\nonumber
   \Con_{F'/F} ( (x)_{0}^{F}) =(x)_{0}^{F'}, \; \Con_{F'/F} ( (x)_{\infty}^{F}) =(x)_{\infty}^{F'} \; \ad \;\Con_{F'/F} ( \div(x)^{F}) =\div(x)^{F'}.
\end{equation}

Let $\fB_1,...,\fB_{\mu}$ be all the places of $F'/F_P$ lying over $P$.
By \cite[Proposition 3.1.4.]{St}, \cite[Definition 3.1.5.]{St} and Theorem D(a), we have
\begin{equation}\label{No54}
   \nu_{\fB_i}(\alpha) = \nu_{P}(\alpha) \quad \for \quad \alpha \in 
	F,   \quad 1 \leq i \leq \mu.
\end{equation}

We will denote by $F^{(P)}$ resp. $F^{'^{(\fB_i)}}$ $(1 \leq i \leq \mu)$ the completion of $F$ 
resp. $F'$ with respect to  the valuation $\nu_P$  resp. $\nu_{\fB_i}$.
Applying  \cite[p.132, 133]{Sa}, we obtain
\begin{equation}\nonumber
  F \subseteq F^{(P)}   \subseteq F^{'^{(\fB_i)}} \quad \ad \quad      F \subseteq F'
    \subseteq F^{'^{(\fB_i)}}, \quad 1 \leq i \leq \mu.
\end{equation}

Let $t$ be a local parameter of $\cP$, and let $\alpha \in F^{(P)} $. By (\ref{No54}), we have
 $  \nu_{\fB_i}(t) =1$. Consider the local expansion (\ref{No06}). Using (\ref{No54}), we get
$\nu_{\fB_i}(\alpha) = \nu_{P}(\alpha)$.
Hence
\begin{equation}\label{No58}
   \nu_{\fB_i}(\alpha) = \nu_{P}(\alpha) \quad \for \quad \alpha \in 
	F'\cap F^{(P)}   \quad 1 \leq i \leq \mu.
\end{equation} \\

{\bf Theorem F.} (\cite[Theorem 2.24]{LiNi}) {\it Let $M$ be a finite extension of the finite field $L$, both
considered as vector spaces over $L$. Then the linear transformations from $M$
into $L$ are exactly the mappings $K_{\beta}$, $\beta \in F$ where
 $K_{\beta} = \Tr_{M/L} (\beta \alpha)$ for all
$\alpha \in F$. Furthermore, we have $K_{\beta} \neq K_{\gamma}$ whenever $\beta$ and $\gamma$ are distinct elements
of $L$.} \\

{\bf Theorem G. } (\cite[Proposition 3.3.3]{St} or \cite[Definition 2.30, and p.58]{LiNi})
 {\it  Let $L$ be a finite field and $M$ a finite extension of $L$. Consider
a basis  $\{\alpha_1,...,\alpha_m\}$ of $M/L$. Then there are uniquely determined elements
${\beta_1,...,\beta_m}$ of $M$, such that
\begin{equation}\label{No26}
   \Tr_{M/L} ( \alpha_i \beta_j) =\delta_{i,j}=  \begin{cases}
        1 \quad {\rm if} \; i=j, \\
			  0  \quad {\rm if} \; i\neq j.
			\end{cases}
\end{equation}
 The set ${\beta_1,...,\beta_m}$
 is a basis of $M/L$ as
well; it is called the dual basis of $\{ \alpha_1,...,\alpha_m  \}$ (with respect to the trace).}
\\

{\bf 4.2 Digital sequences and $(\bT, s) $ sequences} (\cite[Section 4]{DiPi}).\\

{\bf Definition 6.} (\cite[Definition 4.30]{DiPi}) { \it For a given dimension $s \geq 1$, an integer base $b \geq 2$, and a
function $\bT : \NN_0 \to \NN_0$ with $\bT(m) \leq  m$ for all $m \in \NN_0$, a sequence $(\bx_0,\bx_1, . . .)$
of points in $[0, 1)^s$ is called a $(\bT, s)$-sequence in base $b$ if for all integers $m \geq 0$
and $k \geq  0$, the point set consisting of the points $x_{kb^m}, . . . ,x_{kb^m+b^m-1}$ forms
a $(\bT(m),m, s)$-net in base $b$. }\\

{\bf Lemma A.} (\cite[Lemma 4.38]{DiPi}) {\it Let $(\bx_0,\bx_1, . . .)$ be a $(\bT, s)$-sequence in base b. Then, for
every $m$, the point set $\{\by_0, \by_1, . . . , \by_{b^m -1}\}$ with $\by_k := (\bx_k, k/b^m)$, $0 \leq k <
b^m$, is an (r(m),m, s+1)-net in base b with $r(m) := \max\{\bT(0), . . . ,\bT(m)\}$.}\\

Repeating the proof of this lemma, we obtain\\

{\bf Lemma 1}. {\it Let $(\bx_n)_{n \geq0}$ be a sequence in $[0,1)^s$, $m_n \in \NN$,
$m_i >m_j$ for $i>j$, and let $(\bx_n, n/b^{m_k})_{0 \leq n <b^{m_k}}$ be a $(t,m_k,s+1)$ net in base $b$
 for all $k \geq 1$. Then $(\bx_n)_{n \geq0}$ is a  $(t,s)$ sequence in base $b$.} \\

{\bf Lemma B}. (\cite[Lemma 3.7]{Ni}) {\it Let $(\bx_n)_{n \geq 0}$ be a sequence in $[0,1)^s$. For $N \geq 1$, let $H$ be the point set consisting of  $(\bx_n, n/N) \in  [0,1)^{s+1}$
for $n=0,...,N-1$. Then
\begin{equation}  \nonumber
 1+  \max_{1 \leq  M \leq N}  M \emph{D}^{*}((\bx_{n} )_{n=0}^{M -1})
 \geq   N \emph{D}^{*}((  \bx_n, n/N )_{n=0}^{N -1}) .
\end{equation} 
} \\

{\bf Definition 7.} (\cite[Definition 1]{DiNi}) { \it
Let $m, s \geq 1$ be integers. Let $C^{(1,m)},...,\\C^{(s,m)}$ be $m \times m$ matrices over $\FF_b$.
Now we construct $b^m$ points in $[0, 1)^s$.
 For $ n= 0, 1,...,b^m-1$, let $n =\sum^{m-1}_{j=0} a_j(n) b^{j}$ 
be the $b$-adic expansion of $n$. 
 Choose a  bijection  
$\phi : \; Z_b := \{0, 1,...., b-1\} \mapsto \FF_b$ with
$\phi(0) = \bar{0}$, the neutral element of addition in $\FF_b$. Let $|\phi(a)|:= |a|$ for 
$a \in Z_b$.
We identify $n$ with the row vector
\begin{equation} \label{Ap300}
                    \bn = (\bar{a}_0(n),...,\bar{a}_{m-1}(n)) \in  \FF^m_b   
										\quad \with \quad  \bar{a}_i(n) = \phi(a_i(n)), \; 0 \leq i \leq m -1 .
\end{equation}
We map the vectors
\begin{equation} \label{Ap301}
	y^{(i)}_{n} =(y^{(i)}_{n,1},...,y^{(i)}_{n,m})   :=  \bn C^{(i,m)\top}  \in  \FF^{m}_b
\end{equation}
to the real numbers
\begin{equation} \label{Ap302}
   x^{(i)}_n =\sum_{j=1}^m \phi^{-1} (y^{(i)}_{n,j})/b^j
\end{equation} 
to obtain the point
\begin{equation} \label{Ap303}
   \bx_n:= (x^{(1)}_n,...,x^{(s)}_n) \in [0,1)^s.
\end{equation}
 \\


The point set  $ \{\bx_0,...,\bx_{b^m-1} \}$ is called a {\sf digital net} (over $\FF_b$) (with {\sf generating matrices} $(C^{(1,m)},...,C^{(s,m)}) $).
 
For $m = \infty$, we obtain a sequence $\bx_0, \bx_1,...$ of   points in $[0, 1)^s$  which is called a {\sf digital sequence} $($over $\FF_b)$ $($with {\sf generating matrices} $(C^{(1,\infty)},...,C^{(s,\infty)}) )$.}

We abbreviate  $C^{(i,m)}$ as $C^{(i)}$ for $m \in \NN$ and for $m=\infty$.
\\

{\bf 4.3 Duality theory} ( see \cite[Section 7]{DiPi}, [DiNi], \cite{NiPi}, \cite{Skr}).\\

Let $\cN$ be
an arbitrary $\FF_b$-linear subspace of $F_b^{sm}$. Let $H$ be a matrix over $\FF_b$    consisting
of $sm$ columns such that the row-space of $H$ is equal to $\cN$. Then we define
the {\sf dual space} $\cN^{\bot} \subseteq F_b^{sm}$
 of $\cN$  to be the null space of $H$ (see [DiPi, p. 244]). In other words,  
 $\cN^{\bot} $
is the orthogonal complement of $\cN$ relative to the standard inner product
in $F_b^{sm}$,
\begin{equation} \label{Ap304}
 \cN^{\bot} = \{ A \in  F_b^{ sm}  \;  | \; B \cdot A =0 \quad {\rm for \; all} \;\; B \in \cN  \}.
\end{equation}

For  any vector $\ba = (a_1,...,a_m) \in \FF^m_b$, let 
\begin{equation} \nonumber
 v(\ba) = 0 \;\; {\rm if} \;\; \ba = \bs \quad
 \ad  \quad v_m(\ba) = \max \{j \; : \; a_j \neq 0  \}\;\; {\rm if} \;\; \ba \neq \bs.
\end{equation}
 Then we extend this definition 
 to $\FF^{ms}_b$ by writing a vector $\bA \in \FF^{ms}_b$ as the concatenation of $s$ 
 vectors of length $m$, i.e. $\bA = (\ba_1,...,\ba_s) \in \FF^{ms}_b$
 with $\ba_i \in \FF^{m}_b$
for $1 \leq i \leq s$ and putting  
\begin{equation}  \label{NiOz128a}
      V_m(\bA) = \sum_{1 \leq i \leq s} v_m(\ba_i).
\end{equation}

{\bf Definition 8}. For any nonzero $\FF^m_b$-linear subspace $\cN$ of
  $\FF^{ms}_b$, 
the {\sf minimum  distance} of  $\cN$  is defined by
\begin{equation} \nonumber 
       \delta_m(\cN)  =\min \{ V_m(\bA) \;  | \; \bA \in \cN \setminus \{ \bs\}  \}.
\end{equation} \\


Let  $C^{(1)},...,C^{(s)} \in \FF_b^{\infty \times \infty}$  be 
generating matrices of a digital sequence  $\bx_n(C)_{n \geq 0}$ over $\FF_b$.
For any $m \in \NN$, 
we denote the $m \times m$ left-upper sub-matrix of
$C^{(i)}$ by $[C^{(i)}]_m$.
The matrices $[C^{(1)}]_m,...,[C^{(s)}]_m$
 are then the generating matrices
of a digital net. We define the {\sf overall generating matrix} of
this digital net by 
\begin{equation} \label{Di08}
  [C]_m = ([C^{(1)}]_m^{\top}|[C^{(2)}]_m^{\top}|...|[C^{(s)}]_m^{\top}) \in  F_b^{m \times s m}
\end{equation}
for any $m \in \NN$.

Let $\cC_m$ denote the row space of
the matrix $[C]_m$ i.e.,
\begin{equation} \label{Ap309}
\cC_m = \Big\{  \Big(\sum_{r=0}^{m-1}   c^{(i)}_{j,r}  \bar{a}_r(n)\Big)_{0 \leq j
\leq m-1, 1 \leq i \leq s}   \; | \; 0 \leq n < b^m   \Big\}.
\end{equation}
The dual space is then given by
\begin{equation} \label{Di20}
 \cC_m^{\bot} = \{ A \in  F_b^{ sm}  \;  | \; B  \cdot A =\bs \quad {\rm for \; all} \;\; B \in \cC_m   \}.
\end{equation}

{\bf Proposition B}. (\cite[Proposition 7.22]{DiPi} {\it
 For $s \in  \NN$, $s \geq 2$, the matrices $C^{(1)}, . . . ,C^{(s)}$ generate a
digital $(\bT, s)$-sequence if and only if for all $m \in \NN$ we have
\begin{equation} \nonumber
    \bT(m) \geq  m - \delta_m(C^{\bot}_m) + 1, \quad  \fall \quad m \in \NN.
\end{equation} }


We define a weight function on $\FF^{ms}_b$ dual to the weight function $V_m$  (\ref{NiOz128a}).
For  any vector $\ba = (a_1,...,a_m) \in \FF^m_b$, let 
\begin{equation} \label{NiOz127}
 v^{\bot}(\ba) = m+1 \;\; {\rm if} \;\; \ba = \bs \quad
 \ad  \quad v^{\bot}_m(\ba) = \min \{j \; : \; a_j \neq 0  \}\;\; {\rm if} \;\; \ba \neq \bs.
\end{equation}
 Then we extend this definition 
 to $\FF^{ms}_b$ by writing a vector $\bA \in \FF^{ms}_b$ as the concatenation of $s$ 
 vectors of length $m$, i.e. $\bA = (\ba_1,...,\ba_s) \in \FF^{ms}_b$
 with $\ba_i \in \FF^{m}_b$
for $1 \leq i \leq s$ and putting  
\begin{equation} \label{NiOz128}
      V^{\bot}_m(\bA) = \sum_{1 \leq i \leq s} v^{\bot}_m(\ba_i).
\end{equation}

{\bf Definition 9}. For any nonzero $\FF^m_b$-linear subspace $\cN$ of
  $\FF^{ms}_b$, 
the {\sf maximum  distance} of  $\cN$  is defined by
\begin{equation} \label{NiOz129}
       \delta^{\bot}_m(\cN)  =\max \{ V^{\bot}_m(\bA) \;  | \; \bA \in \cN \setminus \{ \bs\}  \}.
\end{equation} \\


{\bf 4.4 Auxiliary  results}.\\

{\bf Lemma C.} (\cite[Lemma 1]{Le4}) {\it Let ${\dot{s}} \geq 2$, $d \geq 1$, $ (\bx_n)_{ 0 \leq n < b^{\tilde{m}}} $  be a $d-$admissible  $(t,\tilde{m},{\dot{s}})$ net in base $b$, 
  $d_0 =d+t$, $\hat{e} \in\NN$,  $0<\epsilon \leq (2d_0 \hat{e} (\dot{s}-1))^{-1}$, $\dot{m} = [\tilde{m} \epsilon]$,
 $\ddot{m}_i=0$,
 $\dot{m}_{i} = d_0\hat{e}\dot{m}$ $(1 \leq i \leq {\dot{s}}-1)$, 
$\ddot{m}_{\dot{s}} =\tilde{m}- ({\dot{s}}-1)\dot{m}_1 -t\geq 1 $,
 $\dot{m}_{\dot{s}} = \ddot{m}_{\dot{s}} +\dot{m}_1$, 
$B_{i} \subset \{0,...,\dot{m}-1\} $
$(1 \leq i \leq {\dot{s}}) $,  $\bw \in E_{\tilde{m}}^{{\dot{s}}}$
   and let
$\gamma^{(i)}=\gamma^{(i)}_1/b+...+\gamma^{(i)}_{\dot{m}_i}/b^{\dot{m}_i}$,
\begin{equation}\label {Ap400}
\gamma^{(i)}_{\ddot{m}_i +d_0(\hat{j}_i\hat{e} +\breve{j}_i) +\check{j}_i} =0 \;\;\for \;\;  1 \leq \check{j}_i < d_0, 
\qquad 
 \; \gamma^{(i)}_{\ddot{m}_i +d_0(\hat{j}_i\hat{e} +\breve{j}_i) +\check{j}_i}=1 \;\; \for \;\; \check{j}_i =d_0  
\end{equation}  
and $ \hat{j}_i  \in \{0,...,\dot{m}-1\} \setminus B_{i}, \; 
  0 \leq \breve{j}_i < \hat{e}, \; 1 \leq i \leq \dot{s} $, 
 $\bgamma =(\gamma^{(1)},...,\gamma^{(\dot{s})})$, 
$B =\#B_1+...+\#B_{\dot{s}}$ and  $ \tilde{m}  \geq 4 \epsilon^{-1}(\dot{s} -1)(1+\dot{s} B) +2t$.
 Let there exists $n_0 \in [0,b^{\tilde{m}})$ such that 
$ [(\bx_{n_0} \oplus \bw)^{(i)}]_{\dot{m}_i} =\gamma^{(i)}$, $1 \leq i \leq \dot{s}$. Then }\\
\begin{equation}  \label{Ap402}
\Delta((\bx_n \oplus \bw)_{0 \leq n < b^{\tilde{m}}}, J_{\bgamma})  \leq
			 		-	b^{-d} \big( \hat{e}\epsilon(2({\dot{s}}-1))^{-1} \big)^{{\dot{s}}-1} \tilde{m}^{{\dot{s}}-1} 
			+b^{t+s}d_0\hat{e}  B\tilde{m}^{\dot{s}-2}.
\end{equation} \\

{\bf Corollary 1.} {\it With notations as above.
Let $\dot{s} \geq 3$, $\tilde{r} \geq 0$, $\tilde{m}=m -\tilde{r} $,
  $ (\bx_n)_{ 0 \leq n < b^{\tilde{m}}} $  be a $d-$admissible  $(t,\tilde{m},\dot{s})$ net in base $b$, 
  $d_0 =d+t$, $\hat{e} \in\NN$, 
$\epsilon =\eta (2d_0 \hat{e} (\dot{s}-1))^{-1}$, 
	$0< \eta \leq 1$, 
	$\dot{m} = [\tilde{m} \epsilon]$,  $\ddot{m}_i=0$,
 $\dot{m}_{i} = d_0\hat{e}\dot{m}$, 
$\ddot{m}_{\dot{s}} =\tilde{m}- ({\dot{s}}-1)\dot{m}_1 -t\geq 1 $,
 $\dot{m}_{\dot{s}} = \ddot{m}_{\dot{s}} +\dot{m}_1$, 
$B_{i} \subset \{0,...,\dot{m}-1\} $,   $\bar{B}_i =\{0,...,\dot{m}-1\} \setminus B_{i}$,  $1 \leq i \leq {\dot{s}}$,  
$B =\#B_1+...+\#B_{\dot{s}}$.    Suppose that
\begin{equation} \label {Ap404}
		\{ (x^{(i)}_{n, \ddot{m}_i +d_0\hat{e}\hat{j}_i +\breve{j}_i} 
		\;| \; \hat{j}_i  \in \bar{B}_{i}, \; \breve{j}_i \in[1, d_0\hat{e}],  \;     i \in [1,\dot{s}] )  \; | \; n \in [0, b^m)           \}=Z_b^{\mu },
\end{equation} 
 with $m \geq 2t+ 8(d+t)\hat{e}(\dot{s}-1)^2\eta^{-1} + 2^{2\dot{s}}b^{d+\dot{s}+t} (d+t)^{\dot{s}} \hat{e}   ( \dot{s}-1)^{2(\dot{s}-1)} 
 \eta^{ -\dot{s} +1} B  +4(\dot{s} -1)\tilde{r}$ and $\mu =d_0 \hat{e}   (\dot{s} \dot{m} -B)$.  
 Then there exists $n_0 \in [0,b^{\tilde{m}})$ such that 
$ [(\bx_{n_0} \oplus \bw)^{(i)}]_{\dot{m}_i} =\gamma^{(i)}$, $1 \leq i \leq \dot{s}$, 
 and   for each  $\bw \in E_{\tilde{m}}^{\dot{s}}$, we have }
\begin{equation}  \nonumber
			b^{\tilde{m}} \emph{D}^{*}((\bx_{n} \oplus \bw)_{0 \leq  n < b^{\tilde{m}}}) \geq  \Big|\Delta((\bx_n \oplus \bw)_{0 \leq n < b^{\tilde{m}
		}}, J_{\bgamma}) \Big| \geq
					2^{-2} b^{-d}  K_{d,t,\dot{s}}^{-\dot{s}+1}\eta^{\dot{s}-1}  m^{\dot{s}-1}
\end{equation}
	$\;\; \with \;\;
					 K_{d,t,\dot{s}} =4(d+t)(\dot{s}-1)^2.$
\\

{\bf Proof.} Let $\bgamma(n,\bw) =\bgamma = (\gamma^{(1)},...,\gamma^{(\dot{s})})$ with 
$\gamma^{(i)}:= [(\bx_{n} \oplus \bw)^{(i)}]_{\dot{m}_i}$, $i \in [1,\dot{s}]$.
Using (\ref{Ap404}), we get that there exists $n_0 \in [0,b^{\tilde{m}})$ such that $\bgamma(n_0,\bw)$
satisfy (\ref{Ap400}). Hence (\ref{Ap402}) is true.
Taking into account (\ref{In02}) and that $\bw \in E_{\tilde{m}}^{\dot{s}}$ is arbitrary, we get the assertion in Corollary 1.
 \qed \\

Let 
$\phi : \; Z_b  \mapsto \FF_b$ be a  bijection with
$\phi(0) = \bar{0}$, and let $x^{(i)}_{n,j}=\phi^{-1}(y^{(i)}_{n,j})$ for 
$1  \leq i \leq s$, $j \geq 1$ and $n \geq 0$. We obtain from Corollary 1 :\\

{\bf Corollary 2.} {\it Let $\dot{s} \geq 3$, $\tilde{r} \geq 0$, $\tilde{m}=m -\tilde{r} $,
  $ (\bx_n)_{ 0 \leq n < b^{\tilde{m}}} $  be a $d-$admissible  $(t,\tilde{m},\dot{s})$ net in base $b$, 
  $d_0 =d+t$, $\hat{e} \in\NN$, 
$\epsilon =\eta (2d_0 \hat{e} (\dot{s}-1))^{-1}$, 
	$0< \eta \leq 1$, 
	$\dot{m} = [\tilde{m} \epsilon]$,  $\ddot{m}_i=0$,
 $\dot{m}_{i} = d_0\hat{e}\dot{m}$, 
$\ddot{m}_{\dot{s}} =\tilde{m}- ({\dot{s}}-1)\dot{m}_1 -t\geq 1 $,
 $\dot{m}_{\dot{s}} = \ddot{m}_{\dot{s}} +\dot{m}_1$, 
$B_{i} \subset \{0,...,\dot{m}-1\} $,   $\bar{B}_i =\{0,...,\dot{m}-1\} \setminus B_{i}$,  $1 \leq i \leq {\dot{s}}$,  
$B =\#B_1+...+\#B_{\dot{s}}$.    Suppose that
\begin{equation} \nonumber
		\{ (y^{(i)}_{n, \ddot{m}_i +d_0\hat{e}\hat{j}_i +\breve{j}_i} 
		\;| \; \hat{j}_i  \in \bar{B}_{i}, \; \breve{j}_i \in[1, d_0\hat{e}],  \;     i \in [1,\dot{s}] )  \; | \; n \in [0, b^m)           \}=\FF_b^{\mu },
\end{equation} 
 with $m \geq 2t+ 8(d+t)\hat{e}(\dot{s}-1)^2\eta^{-1} + 2^{2\dot{s}}b^{d+\dot{s}+t} (d+t)^{\dot{s}} \hat{e}   ( \dot{s}-1)^{2(\dot{s}-1)} 
 \eta^{ -\dot{s} +1} B  +4(\dot{s} -1)\tilde{r}$ and $\mu =d_0 \hat{e}   (\dot{s} \dot{m} -B)$.  
 Then there exists $n_0 \in [0,b^{\tilde{m}})$ such that 
$ [(\bx_{n_0} \oplus \bw)^{(i)}]_{\dot{m}_i} =\gamma^{(i)}$, $1 \leq i \leq \dot{s}$, 
 and   for each  $\bw \in E_{\tilde{m}}^{\dot{s}}$, we have }
\begin{equation}  \nonumber
			b^{\tilde{m}} \emph{D}^{*}((\bx_{n} \oplus \bw)_{0 \leq  n < b^{\tilde{m}}}) \geq  \Big|\Delta((\bx_n \oplus \bw)_{0 \leq n < b^{\tilde{m}
		}}, J_{\bgamma}) \Big| \geq
					2^{-2} b^{-d}  K_{d,t,\dot{s}}^{-\dot{s}+1}\eta^{\dot{s}-1}  m^{\dot{s}-1}.
\end{equation} \\

With notations as above, we consider the case of  $(t, s)$ sequence in base $b$: \\

{\bf Corollary 3.} {\it Let $s \geq 2$, $d \geq 1$,  $ (\bx_n)_{ n \geq 0} $  be a $d-$admissible  $(t, s)$ sequence in base $b$, 
  $d_0 =d+t$, $\hat{e} \in\NN$, $\epsilon =\eta (2d_0 \hat{e} s)^{-1}$, 
	$0< \eta \leq 1$, 
	$\dot{m} = [m \epsilon]$,  $\ddot{m}_i=0$,
 $1 \leq i \leq s$,  
$\ddot{m}_{s+1} =  t-1 +(s-1)d_0 \hat{e}\dot{m}$,
$B_{i}^{'} \subset \{0,...,\dot{m}-1\} $,   $\bar{B}_i^{'} =\{0,...,\dot{m}-1\} \setminus B_{i}^{'}$,  $1 \leq i \leq s+1$,  
$B =\#B_1^{'}+...+\#B_{s+1}^{'}$.
  Suppose that
\begin{eqnarray} \nonumber
	&	\{ (y^{(i)}_{n,\ddot{m}_i +d_0\hat{e}\hat{j}_i +\breve{j}_i} 
		\;| \; \hat{j}_i  \in \bar{B}_{i}^{'}, \; \breve{j}_i \in[1, d_0\hat{e}],  \;     i \in [1,s], \\
	&	\bar{a}_{\ddot{m}_{s+1} +d_0\hat{e}\tilde{j}_{s+1} +\check{j}_{s+1}}(n),	\;  \tilde{j}_{s+1}  \in \bar{B}_{s+1}^{'}, \; \check{j}_{s+1} \in[1, d_0\hat{e}],)  \; | \; n \in [0, b^m)           \}=\FF_b^{\mu }. \nonumber
\end{eqnarray} 
with  $\mu =d_0 \hat{e}   ((s+1) \dot{m} -B)$, and 
 $m \geq 2t+ 8(d+t)\hat{e}s^2\eta^{-1} + 2^{2s+2}b^{d+s+t+1} (d+t)^{s+1} \hat{e}s^{2s} 
 \eta^{ -s} B$.   Then }
\begin{equation}  \nonumber
 1+  \min_{0 \leq Q <b^m} \min_{ \bw \in E_m^{s}}   \max_{1 \leq N \leq b^m} 
		N \emph{D}^{*}((\bx_{n\oplus Q} \oplus \bw)_{0 \leq  n < N})  \geq   
					2^{-2} b^{-d}  K_{d,t,s+1}^{-s}\eta^s m^{s}. 
\end{equation}\\

{\bf Proof.}   Using  Lemma B, we have
\begin{eqnarray*}  \nonumber
     1+  \sup_{1 \leq N \leq b^m}  N \emph{D}^{*}((\bx_{n\oplus Q} \oplus \bw)_{0 \leq  n < N})
 \geq   b^m \emph{D}^{*}((  \bx_{n\oplus Q} \oplus \bw, n/b^m )_{0 \leq  n < b^m}) \\
		=   
		b^m \emph{D}^{*}((  \bx_n \oplus \bw, (n \ominus Q)/b^m )_{0 \leq  n < b^m}) .
\end{eqnarray*}
By  (\ref{In04}) and [DiPi, Lemma 4.38], we have that $((  \bx_n, n/b^m )_{0 \leq  n < b^m})$   is a $d-$admissible $(t,m,s+1)-$net in base~$b$.
We apply Corollary 2 with $\dot{s} =s+1$,  $\tilde{r}=0$, $B_i^{'} =B_i$, $1 \leq i < \dot{s}$,
$B_{\dot{s}}^{'} = \{\dot{m} -j-1 | j \in B_{\dot{s}}  \}$, 
$\hat{j}_{s+1} = \dot{m} -\tilde{j}_{s+1}-1$, $\breve{j}_{s+1} = 
d_0 \hat{e} -\check{j}_{s+1}+1$,
 and  $x_n^{(s+1)} =n/b^{m}$.
Taking into account that  $y^{(s+1)}_{n,m -j}   =  \bar{a}_{j}(n)$ $(0 \leq j < m)$,  we get
 $ y^{(s+1)}_{n,m -\ddot{m}_{s+1} -d_0 \hat{e} \dot{m} -1 +d_0\hat{e}\hat{j}_{s+1} +\breve{j}_{s+1}}  =  \bar{a}_{\ddot{m}_{s+1} +d_0\hat{e}\tilde{j}_{s+1} +\check{j}_{s+1}}(n)$, 
and   Corollary 3 follows.
 \qed \\

{\bf Lemma 2.} {\it Let $\dot{s} \geq 2$, 
  $d_0 \geq 1$, $\hat{e} \geq 1$, $\dot{m} \geq 1$, $\dot{m}_1 =d_0\hat{e}\dot{m}$,
	$\ddot{m}_i \in [0,m-\dot{m}_1]$ $(1 \leq i \leq \dot{s})$, $m \geq \dot{s}\dot{m}_1$,
	$\dot{m} \geq r$, and let 
\begin{equation} \label{Ap200}
    \Phi := \{ (y^{(1)}_{n,\ddot{m}_1+1},...,y^{(1)}_{n,\ddot{m}_1+\dot{m}_1}, ..., 
y^{({\dot{s}})}_{n,\ddot{m}_{\dot{s}}+1},...,
			y^{({\dot{s}})}_{n,\ddot{m}_{\dot{s}}+\dot{m}_1} )   |  n \in [0, b^m)           \} \subseteq \FF_b^{\dot{s}\dot{m}_1}.
\end{equation}
 Suppose that $\Phi$ is a $\FF_b$ linear subspace of $\FF_b^{\dot{s}\dot{m}_1}$ and $\dim_{\FF_b} (\Phi) =\dot{s} \dot{m}_1 -r$.
Then there exists $B_i \in \{0,...,\dot{m}-1\}$, $1 \leq i \leq \dot{s}$, with 
$B =\#B_1+...+\#B_{\dot{s}} \leq r$ and  

\begin{equation}\label{Ap202}
    \Psi  =\FF_b^{d_0\hat{e}(\dot{s}\dot{m}- B)}, 
\end{equation}
where
\begin{equation} \label{Ap204}
    \Psi = \{ (y^{(i)}_{n, \ddot{m}_i +d_0\hat{e}(\dot{j}_i-1) +\ddot{j}_i} 
		\;|\; \dot{j}_i  \in \bar{B}_{i}, \; \ddot{j}_i \in[1, d_0\hat{e}],  \;     i \in [1,\dot{s}] )  \; | \; n \in [0, b^m)           \}
\end{equation}
with $\bar{B}_i =\{0,...,\dot{m}-1\} \setminus B_{i}$.} \\

{\bf Proof.} Let $\hat{r} =\dot{s} \dot{m}_1 -r $, and let $\bff_1,...,\bff_{\hat{r}}$ be a basis of $\Phi$ with		
\begin{equation} \nonumber
\bff_{\mu} =(f^{(1)}_{\mu,\ddot{m}_1+1},...,f^{(1)}_{\mu,\ddot{m}_1+\dot{m}_1}, ... 
,f^{(s)}_{\mu,\ddot{m}_s+1},...,
			f^{(s)}_{\mu,\ddot{m}_s+\dot{m}_1} ), \;1 \leq \mu \leq {\hat{r}} .
\end{equation}			
Let
\begin{equation} \nonumber
  v(\bff_{\mu}) = \max \big\{ \ddot{m}_i  +(i-1)\dot{m}_1+ j \;| \;
	f^{(i)}_{\mu, \ddot{m}_i +j} \neq 0, \; j \in[1, \dot{m}_1],     i \in [1,\dot{s}]   \big\} \; \for\;	 \ \mu \in[1,{\hat{r}}].
\end{equation}			
Without loss of generality, assume now that $v(\bff_{i}) \leq v(\bff_{j})$ for $1 \leq i < j \leq {\hat{r}}$.
Let $v(\bff_{j}) = \ddot{m}_{l_1}+(l_1-1)\dot{m}_1 +l_2$, and let $\dot{\bff}_k =\bff_k -\bff_{j} f^{(l_1)}_{k,\ddot{m}_{l_1}+l_2} /  f^{(l_1)}_{j,\ddot{m}_{l_1}+l_2}$ 
for $1 \leq k \leq j-1$.\\
 We have $v(\dot{\bff}_k) < v(\bff_j)$ for all $1 \leq k \leq j-1$.

By repeating this procedure for $j=\hat{r},\hat{r}-1,...,2$, we obtain a basis $\hat{\bff}_1,...,\hat{\bff}_{\hat{r}}$  of $\Phi$ with $v(\hat{\bff}_{i}) < v(\hat{\bff}_{j})$ for $1 \leq i < j \leq {\hat{r}}$. 
Let 
\begin{equation} \nonumber
 A_i =\{ \ddot{m}_i +j \; |   \; 
v(\hat{\bff}_{\mu})  = (i-1)\dot{m}_1+ \ddot{m}_i  +j, \; 1 \leq j \leq \dot{m}_1,  1 \leq \mu \leq {\hat{r}} \}, 
 i \in [1,\dot{s}] .
\end{equation}
Taking into account that $\hat{\bff}_1,...,\hat{\bff}_{\hat{r}}$ is a basis  of $\Phi$, we get
 from (\ref{Ap200})
\begin{equation} \label{Ap206}
  \{ (y^{(i)}_{n, j} 
		\;|\; j  \in A_i,  \;     i \in [1,\dot{s}] )  \; | \; n \in [0, b^m)           \}
		=\FF_b^{\dot{s} \dot{m}_1 -r}.
\end{equation}
Now let
\begin{equation} \nonumber
 \bar{B}_i := \{ \dot{j}_i \in [0,\dot{m}_1)  
		\;|\;  \exists \ddot{j}_i  \in [1,d_0 \hat{e}],  \; \with \;  \ddot{m}_i +\dot{j}_id_0 \hat{e} +\ddot{j}_i \in A_i   )          \}, \; i \in [1,\dot{s}] .
\end{equation}
It is easy to see that
$B =\#B_1+...+\#B_{\dot{s}} \leq r$,  where $\bar{B}_i =\{0,...,\dot{m}-1\} \setminus B_{i}$.

 Bearing in mind (\ref{Ap204}), we obtain (\ref{Ap202}) from (\ref{Ap206}).
Hence Lemma 2 is proved. \qed

\section{Statements of results.}
 If $s=2$ for the case of nets, or  $s=1$ for the case of sequences,  then (\ref{In10}) follows from the W. Schmidt estimate (\ref{In08}) 
(see [Ni, p.24]).
In this paper we take  $s \geq 2$ for the case of sequences, and  
$s \geq 3$ for the case of nets. \\

{ \bf 2.1  Generalized Niederreiter sequence}.
In this subsection, we introduce a generalization of the Niederreiter sequence
due to Tezuka (see \cite[Section 6.1.2]{Te2}, \cite[Section 8.1.2]{DiPi}). 
By \cite[p.165] {Te2}, the Sobol's sequence \cite[Section 8.1.2] {DiPi}, the Faure's sequence \cite[Section 8.1.2] {DiPi}) and the original Niederreiter sequence 
\cite[Section 8.1.2] {DiPi}) are particular cases of a generalized Niederreiter sequence.

Let $b$ be a prime power and let $p_1, . . . , p_{s} \in F_b[x]$ be  pairwise coprime polynomials
 over $\FF_b$. Let $e_i = \deg(p_i) \geq 1$ for $1 \leq i \leq s$.
For each $j \geq  1$ and $1 \leq i \leq s$, the set of polynomials
$\{y_{i,j,k}(x) \; : \; 0 \leq  k < e_i\}$ needs to be linearly independent $(\mod \;p_i(x))$ over
$\FF_b$.
 For integers
$1 \leq i \leq s$, $j \geq 1$ and $0 \leq k < e_i$, consider the expansions
\begin{equation} \label{Ni00}
  \frac{ y_{i,j,k}(x)}{
p_i(x)^{j}} = \sum_{r\geq 0} a^{(i)} (j, k, r) x^{-r-1}
\end{equation}
over the field of formal Laurent series $F_b((x^{-1}))$.  Then we define the matrix
 $C^{(i)} = (c^{(i)}_{ j,r})_{j \geq 1, r \geq 0}$ by
\begin{equation}  \nonumber
c^{(i)}_{ j,r} = a^{(i)}(Q + 1, k, r) \in \FF_b \qquad \for \qquad1 \leq i \leq s,\;  j \geq 1,\; r \geq 0,
\end{equation}
where $j -1 = Qe_i + k$ with integers $Q = Q(i, j)$ and $k = k(i, j)$ satisfying
 $0 \leq k < e_i$. \\

 A digital sequence $(\bx_n)_{n \geq 0}$ over $\FF_b$  generated by the matrices $C^{(1)},...,C^{(s)}$  is called a
{\sf generalized Niederreiter sequence}  (see [DiPi, p.266]). \\

{\bf Theorem H.} (see [DiPi, p.266]) {\it The generalized Niederreiter sequence with generating matrices, defined
as above, is a digital (t, s)-sequence over $\FF_b$ with
$t =e_0-s$ and $e_0=e_1+...+e_s.$}\\

{\bf Theorem 1.} {\it  With the notations as above, 
$ (\bx_n)_{ n \geq 0} $  is  $d-$admissible with $d=e_0$.\\
 (a) For $s \geq 2$, $e=e_1 e_2 \cdots e_s$, $\eta_1 =s/(s+1)$ $ m \geq  9(d+t)es(s+1)$ and $K_{d,t,s} =  4 (d+t)(s-1)^2$, we have
 \begin{equation} \nonumber 
    1+  \min_{0 \leq Q <b^m} \min_{\bw \in E_m^{s}} \max_{1 \leq N \leq b^m} 
		N \emph{D}^{*}((\bx_{n\oplus Q} \oplus \bw)_{0 \leq  n < N}) \geq   
					2^{-2} b^{-d}K_{d,t,s+1}^{-s} \eta_1^{s}m^{s} .
\end{equation} 
(b) Let $s \geq 3$, $\eta_2 \in (0,1)$ and $ m \geq  8(d+t)e(s-1)^2 \eta_2^{-1}+
2(1+t) \eta_2^{-1}(1-\eta_2)^{-1}$. Suppose that $ \min_{m/2-t \leq j e_{i_0}   \leq m, 0 \leq k < e_{i_0}}(1-\deg( y_{i_0,j,k}(x))j^{-1}  e_{i_0}^{-1})  \geq \eta_2$ for some $i_0 \in [1,s]$. Then}
\begin{equation} \nonumber
      \min_{ \bw \in E_m^s}   b^m \emph{D}^{*}((\bx_n \oplus \bw)_{0 \leq  n < b^m}) \geq   
		2^{-2} b^{-d}K_{d,t,s}^{-s+1} \eta_2^{s-1} m^{s-1} .
\end{equation}\\

{\bf 2.2 Xing-Niederreiter sequence} (see  \cite[Section 8.4 ]{DiPi}). 
Let $F/\FF_b$ be an algebraic function field with full constant field $\FF_b$ and
genus $g = g(F/\FF_b)$. 
 Assume that $F/\FF_b$ has at least one
rational place $P_{\infty}$, and let $G$ be a positive divisor of $F/\FF_b$ with 
$\deg(G) = 2g$ and $P_{\infty} \notin  {\rm supp}(G)$. Let $P_1, . . . , P_s$ be $s$ distinct places of $F/\FF_b$ with $P_i \neq P_{\infty}$
for $1 \leq  i  \leq  s$. Put $e_i = \deg(P_i)$ for $1 \leq i\leq  s$. 

By \cite[p.279 ]{DiPi}, we have that there exists a basis   ${w_0,w_1, . . . ,w_g}$  of $\cL(G)$ over $\FF_b$ such that 
\begin{equation}  \nonumber
            \nu_{P_{\infty}}(w_u) = n_u \quad {\rm for} \quad 0 \leq u \leq g,   
\end{equation}
where $0 = n_0 <
n_1 < .... < n_g  \leq 2g$.
For each $1 \leq i  \leq s$, we consider the chain
\begin{equation}  \nonumber
          \cL(G) \subset \cL(G + P_i) \subset \cL(G + 2P_i)\subset...
\end{equation}
of vector spaces over $\FF_b$. By starting from the basis ${w_0,w_1, . . . ,w_g}$ of $\cL(G)$ 
and successively adding basis vectors at each step of the chain, we obtain
for each $n \in \NN$ a basis
\begin{equation}\label{XiNi04}
         \{w_0,w_1, . . . ,w_g, k_{i,1},k_{i,2},...,k_{i,n e_i}\}
\end{equation}
of $\cL(G + n P_i)$. We note that we then have
\begin{equation}\label{XiNi06}
    k_{i,j} \in \cL(G + ([(j-1)/e_i+1)]P_i) \quad {\rm for} \quad1 \leq i  \leq  s \quad {\rm and } \quad j \geq 1.
\end{equation}
%

 By the Riemann-Roch  theorem, there exists a local parameter $z$ at $P_{\infty}$, e.g.,  with 
\begin{equation}\label{XiNi07}
             \deg((z)_\infty)  \leq 2g+e_1 \qquad \for \qquad 
             z \in \cL(G+P_1 -P_{\infty}) \setminus \cL(G+P_1 -2P_{\infty}) . 
\end{equation}
For $r \in \NN \cup \{ 0\}$, we put
\begin{equation}\label{XiNi08}
   z_r=  \begin{cases}
        z^r \quad {\rm if} \; r \notin \{ n_0,n_1,...,n_g  \}, \\
				  w_u  \quad {\rm if} \; r= n_u \;\; {\rm for \; some} \; u \in \{0,1,...,g\}.
			\end{cases}
\end{equation}
Note that in this case $\nu_{P_{\infty}}(z_r) = r$ for all $r \in \NN \cup \{ 0\}$. For $1 \leq i \leq s$ and $j \in \NN$,
we have $k_{i,j} \in  \cL(G + n P_i)$ for some $n \in \NN$ and also 
$P_{\infty} \notin {\rm supp}(G + nP_i)$,
hence $\nu_{P_{\infty}}(k_j^{(i)}) \geq 0$.
 Thus we have the local expansions
\begin{equation}\label{XiNi10}
 k_{i,j} = \sum_{  r = 0}^{\infty} a_{j,r}^{(i)} z_r  \quad {\rm for} \;\; 
    1 \leq i \leq s \quad {\rm and} \;\; j \in \NN,
\end{equation}
where all coefficients $a_{j,r}^{(i)} \in \FF_b$. For $1 \leq i \leq s$ and $j \in \NN$, we now define the
sequences
\begin{equation}\label{XiNi12}
  \bc^{(i)}_j =(c_{j,0}^{(i)}, c_{j,1}^{(i)},...) :=(a_{j,n}^{(i)})_{n \in \NN_0 \setminus \{n_0,...,n_g \}}
\end{equation}
\begin{equation}\nonumber
   = ( \widehat{a_{j,n_0}^{(i)}},a_{j,n_0+1}^{(i)},...,\widehat{a_{j,n_1}^{(i)}},
			a_{j,n_1+1}^{(i)},....,  \widehat{a_{j,n_g}^{(i)}},a_{j,n_g+1}^{(i)},....)   \in \FF_b^{\NN},
\end{equation}
where the hat indicates that the corresponding term is deleted. We define
the matrices $C^{(1)}, . . . ,C^{(s)} \in \FF^{\NN \times \NN}_b$ by
\begin{equation}\label{XiNi14}
   C^{(i)} =(\bc^{(i)}_1,\bc^{(i)}_2,\bc^{(i)}_3,...)^{\top} \quad {\rm for } \quad 1 \leq i \leq s,
\end{equation}
i.e., the vector $\bc^{(i)}_j$
 is the $j$th row vector of $C^{(i)}$ for $1 \leq i \leq s$.\\

{ \bf Theorem I} (see \cite[Theorem 8.11]{DiPi}).  {\it With the above notations,  we have that the matrices $C^{(1)}, . . . ,C^{(s)}$
 given by (\ref{XiNi14}) are generating matrices of the Xing-Niederreiter $(t, s)$-sequence $ (\bx_n)_{ n \geq 0} $ with
$t = g +e_0 -s $ and $e_0=e_1+...+e_s$.} \\

{\bf Theorem 2.} {\it  With the above notations, 
$ (\bx_n)_{ n \geq 0} $  is  $d-$admissible, where $d=g+e_0$.\\
 (a) For  $s \geq 2$, $e=e_1...e_s$, $ m \geq  9(d+t)es^2\eta_1^{-1}$ and $K_{d,t,s} =  4 (d+t)(s-1)^{2}$, we have
 \begin{equation} \nonumber
    1+  \min_{0 \leq Q <b^m} \min_{\bw \in E_m^{s}} \max_{1 \leq N \leq b^m} 
		N \emph{D}^{*}((\bx_{n\oplus Q} \oplus \bw)_{0 \leq  n < N}) \geq   
					2^{-2} b^{-d}K_{d,t,s+1}^{-s} \eta_1^s m^{s} 
\end{equation} 
with $\eta_1 =(1+\deg((z)_\infty))^{-1} $  (see (\ref{XiNi07})).\\
(b) Let $s \geq 3$, $\eta_2 \in (0,1)$ and $m \geq   8(d+t)e(s-1)^2\eta_2^{-1} +
 2(1 +2g+ \eta_2t) \eta_2^{-1}  (1-\eta_2)^{-1}$. Suppose that $ \min_{m/2 -t \leq j   \leq m}\nu_{P_{\infty}}(k_{i_0,j})/j \geq \eta_2$, for some $i_0 \in [1,s]$. Then}
\begin{equation}\label{XiNi15a}
      \min_{ \bw \in E_m^s}   b^m \emph{D}^{*}((\bx_n \oplus \bw)_{0 \leq  n < b^m}) \geq   
		2^{-2} b^{-d}K_{d,t,s}^{-s+1} \eta_2^{s-1} m^{s-1} .
\end{equation}\\

 { \bf 2.3 Niederreiter-\"{O}zbudak nets} (see  \cite[Section 8.2 ]{DiPi}). 
Let $F/\FF_b$ be an algebraic function field with full constant field $\FF_b$ and
genus $g = g(F/\FF_b)$. 
Let $s \geq 2$, and let $P_1, . . . , P_s$ be $s$  distinct places  of $F$ with degrees $e_1, . . . , e_s$. For
$1 \leq i \leq s$, let $\nu_{P_i}$ be the normalized discrete valuation of $F$ corresponding to
$P_i$,  let $t_i$ be a local parameter at $P_i$. Further, for each $1 \leq i \leq s$, let $F_{P_i}$
be the residue class field of $P_i$, i.e., $F_{P_i} = O_{P_i}/P_i $, and let 
$\vartheta_i = (\vartheta_{i,1} ,..., \vartheta_{i, e_i}) : F_{P_i}  \to  \FF^{e_i}_b$
be an $\FF_b$-linear vector space isomorphism. Let $m > g + \sum_{i=1}^s (e_i -1)$. Choose
an arbitrary divisor $G$ of $F/\FF_b$ with $\deg(G) = ms - m + g - 1$ and define
$a_i := \nu_{P_i}(G)$ for $1 \leq i \leq s$.
For each $1 \leq i \leq s$, we define an $F_b$-linear map
$\theta_i  \; : \; \cL(G) \to  \FF^m_b$
on the Riemann-Roch space
$\cL(G) = \{ y \in F\setminus {0} : \div(y) + G  \geq 0 \} \cup \{0\}$.
We fix $i$ and repeat
the following definitions related to $\theta_i$ for each $1 \leq i \leq s$.

Note that for each $f \in \cL(G)$ we have $\nu_{P_i}(f) \geq  - a_i$, and so the local
expansion of $f$ at $P_i$ has the form
\begin{equation}\label{NiOz02}
f = \sum_{j=-a_i}^{\infty}  S_{j}(t_i,f) t^j_i , \quad {\rm with } \quad
  S_{j}(t_i,f)  \in  F_{P_i}, \; j \geq -a_i.
\end{equation}
 We denote $S_{j}(t_i,f) $ by $f_{i,j}$. Let $m_i = [m/e_i]$ and $r_i = m - e_im_i$. Note
that $0 \leq  r_i < e_i$. For $f \in \cL(G)$, the image of $f$ under $\theta^{(G)}_i$, for $1 \leq i \leq s$, is
defined as
\begin{equation}\label{NiOz04}
             \theta^{(G)}_i(f) =(\theta_{i,1}(f),...,\theta_{i,m}(f)) := (\bs_{r_i},   \vartheta_{i}(f_{i, - a_i+m_i - 1}), . . . , \vartheta_{i}(f_{i, - a_i})) \in \FF^m_b,
\end{equation}
where we add the $r_i$-dimensional zero vector $\bs_{r_i}=(0, . . . , 0) \in \FF^{r_i}_b$
 in the beginning.
Now we set
\begin{equation}\label{NiOz06}
             \theta^{(G)}(f) :=  (\theta^{(G)}_1(f), . . . , \theta^{(G)}_s(f)) \in \FF^{ms}_b,
\end{equation}
and define the $\FF_b$-linear map
\begin{equation}\nonumber
             \theta^{(G)} \; : \; \cL(G)   \to  \FF^{ms}_b , \qquad      f  \;  \mapsto \; \theta^{(G)}(f) . 
\end{equation}

 The image of $\theta^{(G)}$
 is denoted by 
\begin{equation}\label{NiOz07}
         \cN_m = \cN_m(P_1,...,P_s;G):= \{    \theta^{(G)}(f)  \in \FF^{ms}_b\;  |  \;f \in \cL(G) \}.
\end{equation}
 According to [DiPi, p.274],
\begin{equation}\nonumber
          \dim(\cN_m) = \dim(\cL(G)) \geq \deg(G) +1-g=ms-m 
  \quad {\rm for } \quad m > g-s+ e_1+...+e_s.
\end{equation}
Using the Riemann-Roch theorem, we get
\begin{equation}\label{NiOz08}
          \dim(\cN_m) = ms-m 
  \quad {\rm for } \quad m > g-s+ e_1+...+e_s,\; s \geq 3.
\end{equation}

Let $\cN_m^{\bot}=\cN_m^{\bot}(P_1,...,P_s;G)$ be the dual space of $\cN_m(P_1,...,P_s;G)$ 
(see  (\ref{Di20})). 
The space $ \cN_m^{\bot}$ can be viewed as the row space of a suitable $ m \times ms$ matrix $C$ over $\FF_b$. Finally, we consider the digital net 
$\cP_1(\cN_m^{\bot})=\{\bx_n(C)  |  n\in[0,b^m)\}$ with overall generating matrix 
$C$ (see  (\ref{Di08})).

Let $ 	\tilde{x}_i(h_i) =\sum_{j=1}^m \phi^{-1} (h_{i,j})b^{-j}$, where $h_i=(h_{i,1},...,h_{i,m}) \in  F_b^{m}$ $(i=1,...,s)$ and let $\tilde{\bx}(\bh) =(\tilde{x}_1(h_1),...,\tilde{x}_s(h_s) )$ where $\bh=(h_{1},...,h_{s})$.  From (\ref{Ap302}), (\ref{Ap303}) and (\ref{Ap309}),
 we derive
\begin{equation}\label{NiOz08b0}
  \cP_1:= \cP_1(\cN_m^{\bot} )=\{\tilde{\bx}(\bh) \; | \; 
 \bh \in \cN_m^{\bot}(P_1,...,P_s;G) \}.
\end{equation}
\\

{\bf Theorem J}  (see \cite[Corollary 8.6]{DiPi}).{\it With the above notations, we have that $\cP_1 $ is a $(t,m,s)$-net over $\FF_b$ with $t=g +e_0-s $ and $e_0=e_1+...+e_s$.}\\

To obtain a $d-$admissible net, we will consider also the following net:
\begin{equation}\label{NiOz08b1}
    \cP_2 := \{ ( \{b^{r_1}z_1\},...,\{b^{r_s}z_s\})  \; | \;   \bz=(z_1,...,z_s) \in \cP_1\}.
\end{equation} 
Without loss of generality, let
\begin{equation}\label{NiOz10b}
        e_s =\min_{1 \leq i \leq s} \; e_i.
\end{equation}
\\

{\bf Theorem 3.} {\it 
Let $s \geq 3$,  
 $ m_0= 2^{2s+3} b^{d+t+s}(d+t)^{s}(s-1)^{2s-1} (g+e_0) e \eta^{-s+1}$ and $\eta =(1+\deg((t_s)_\infty))^{-1} $. Then\\
 \begin{equation}  \nonumber
     \min_{\bw \in E_m^{s}} \max_{1 \leq N \leq b^m} 
		N \emph{D}^{*}(\cP_1 \oplus \bw) \geq   
					2^{-2} b^{-d}K_{d,t,s}^{-s+1} \eta^{-s+1} m^{s-1}, \quad \for \quad m \geq m_0,
\end{equation} 
$\cP_2 $ is a $d-$admissible $(t,m-r_0,s)$ net in base $b$ with $d=g+e_0$,  $t=g+e_0-s$,
 and
\begin{equation} \nonumber
      \min_{ \bw \in E_{m-r_0}^s}   b^m \emph{D}^{*}((\cP_2 \oplus \bw)) \geq   
		2^{-2} b^{-d}K_{d,t,s}^{-s+1} \eta^{s-1} m^{-s+1}, \quad \for \quad m \geq m_0,
\end{equation}
where $\cP_i \oplus \bw :=  \{ \bz \oplus \bw  \; | \;\bz \in\cP_i \}$.
\\
.}

{\bf 2.4 Halton-type sequence}  (see  \cite{NiYe}). 
Let $F/\FF_b$ be an algebraic function field with full constant field $\FF_b$ and
genus $g = g(F/\FF_b)$. 
We assume that $F/\FF_b$ has at least one rational place, that is, a place of degree 1. Given a dimension $s \geq 1$, we choose $s +1$ distinct places  $P_1$,...,$P_{s+1}$ of $F$ with deg$(P_{s+1}) = 1$. The
degrees of the places $P_1$,...,$P_s$ are arbitrary and we put $e_i = \deg (P_i)$ for
 $1 \leq i \leq s$. Denote
by $O_F$ the holomorphy ring given by
\begin{equation} \nonumber
               O_F = \bigcap_{ P  \neq P_{s+1}} O_P ,
\end{equation}
where the intersection is extended over all places $ P  \neq P_{s+1} $ of $F$, and $O_P$ is the valuation ring of $P$. 
We arrange the elements of $O_F$ into a sequence by using the fact that
\begin{equation} \nonumber
   O_F = \bigcup_{m=0}^{\infty} \cL(mP_{s+1}) .				
\end{equation}
The terms of this sequence are denoted by $f_0, f_1, . . .$ and they are obtained as follows. Consider the chain
\begin{equation} \nonumber
   \cL(0) \subseteq L(P_{s+1}) \subseteq L(2P_{s+1}) \subseteq \cdots				
\end{equation}
of vector spaces over $\FF_b$. At each step of this chain, the dimension either remains the same
or increases by $1$. From a certain point on, the dimension always increases by $1$ according to
the Riemann-Roch theorem. Thus we can construct a sequence $v_0, v_1, . . .$ of elements of $O_F$
such that 
\begin{equation}\label{Hal03}
        \{v_0, v_1, . . . , v_{\ell(mP_{s+1})-1} \}
\end{equation}
is a $\FF_b$-basis of $\cL(m P_{s+1})$. For $n \in  \NN$, let
\begin{equation}  \nonumber
 n = \sum_{  r = 0}^{\infty} a_r(n) b^r \quad {\rm with \; all} \;\;  a_r(n) \in Z_b
\end{equation}
be the digit expansion of $n$ in base $b$. Note that $a_r(n) = 0$ for all sufficiently large $r$. We fix
a bijection $ \phi \;: \; Z_b \to \FF_b$ with $\phi(0) = \bar{0}$. Then we define
\begin{equation}\label{Hal06}
 f_n = \sum_{  r = 0}^{\infty} \bar{a}_r(n) v_r  \in O_F\quad {\rm with} \quad
             \bar{a}_r(n)=\phi(a_r(n)) 
\quad {\rm for} \;\;  n=0,1,... \; .
\end{equation}
Note that the sum above is finite since for each $n \in \NN$ we have $a_r(n) = 0$ for all sufficiently
large $r$.
 By the Riemann-Roch theorem, we have
\begin{equation}  \label{Hal07}
 \{\tilde{f} \; |\; \tilde{f} \in \cL((m+g-1) P_{s+1}) \} =
 \{f_n \; |\; n \in [0,b^m) \} 
\quad \for \quad m \geq g.
\end{equation}
 For each $i = 1, . . . , s$, let $\wp_i$ be the maximal ideal of $O_F$ corresponding to $P_i$. Then the
residue class field $F_{P_i} := O_F /\wp_i$ has order $b^{e_i}$ (see [St, Proposition 3.2.9]). We fix a bijection
\begin{equation} \label{Ha08}
                \sigma_{P_i}  \; : \;F_{P_i} \to Z_{b^{e_i}} .
\end{equation}
 For each
$i = 1, . . . , s$, we can obtain a local parameter $t_i \in O_F$ at $\wp_i$, by applying the Riemann-Roch
theorem and choosing
\begin{equation}\label{Hal10}
               t_i \in  \cL(k P_{s+1} - P_i) \setminus  \cL(k P_{s+1} - 2P_i)
\end{equation}
for a suitably large integer $k$. We have a local expansion of 
$f_n$ at $\wp_i$ of the form
\begin{equation}\label{Hal12}
  f_n = \sum_{  j \geq 0}  f^{(i)}_{n,j} t_i^j    \quad {\rm with \; all} \;\;   f^{(i)}_{n,j}    \in F_{P_i} ,\; n=0,1,...\; .
\end{equation}
We define the map $ \xi  \; : \; O_F  \to [0, 1]^s$ by
\begin{equation}\label{Hal14}
       \xi(f_n) = \Big(  \sum_{  j = 0}^{\infty}  \sigma_{P_1} \big( f^{(1)}_{n,j}\big) 
			b^{-e_1(j+1)},...,
		  \sum_{  j = 0}^{\infty}  \sigma_{P_s} \big( f^{(s)}_{n,j}\big) (b^{-e_s(j+1)} \Big) .				
\end{equation}
Now we  define the sequence $\bx_0, \bx_1, . . .$ of points in $[0, 1]^s$ by
\begin{equation}\label{Hal16}
 \bx_n  = \xi(f_n)   \quad {\rm for} \quad  n=0,1,... \; .
\end{equation} 
 From [NiYe, Theorem 1], we get the following theorem :\\

{\bf Theorem K.}  {\it With the notation as above, we have that $(\bx_n)_{n \geq 0} $ 
is a $(t, s)$-sequence in base $b$  with $t=g+e_0 -s$ and $e_0 =e_1+...+e_s$.}\\

 By Lemma 17, $(\bx_n)_{n \geq 0} $ 
is  $d-$admissible  with $d=g+e_0$. Using \cite[Theorem 2]{Le4}, we get 
\begin{equation}  \label{Hal16c}
    1+   \max_{1 \leq N \leq b^m, } 
		N \emph{D}^{*}((\bx_{n\oplus Q} \oplus \bw)_{0 \leq  n < N}) \geq   
					2^{-2} b^{-d}K_{d,t,s+1}^{-s} m^{s} 
\end{equation}
for some $Q \in [0,b^m)$ and $\bw \in E_m^{s}$.

 In order to obtain (\ref{Hal16c}) for every $Q$ and $\bw$, we choose a specific sequence $v_0,v_1,...$
as follows. Let 
\begin{equation}\nonumber
    t_{s+1} \in \cL(([(2g+1)/e_1]+1)P_1  -P_{s+1}) \setminus \cL(([(2g+1)/e_1]+1)P_1-2P_{s+1}).
\end{equation}
It is easy to see  that
\begin{equation}\label{Hal17}
          \nu_{P_{s+1}} (t_{s+1}) =1,\;\;  \nu_{P_{i}} (t_{s+1}) \geq 0, \; i \in [2,s]  \;\; \ad \;\; \deg ( (t_{s+1})_{\infty}) \leq 2g+e_1+1.
\end{equation}
By (\ref{Hal03}) and the Riemann-Roch theorem, we have $\nu_{P_{s+1}} (v_{i}) =-i-g$ for $i \geq g$. Hence
\begin{equation}\label{Hal18}
   v_i = \sum_{  j \leq i+g}  v_{i,j} t_{s+1}^{-j}    \quad {\rm with \quad all} \quad   v_{i,j}   \in \FF_b, \quad v_{i,i+g} \neq 0, \quad i \geq g .
\end{equation}

Using the orthogonalization procedure, we can construct a sequence $v_0, v_1, . . .$ 
such that $ \{v_0, v_1, . . . , v_{\ell(mP_{s+1})-1} \}$ is a $\FF_b$-basis of $\cL(m P_{s+1})$,  
\begin{equation}\label{Hal18a}
      v_{i,i+g} =1, \quad \ad \quad  v_{i,j+g} =0 \quad \for \quad  j \in [g,i), \quad i \geq g.
\end{equation}
Subsequently, we will use just this sequence.\\

{\bf Theorem 4.}    {\it  With the above notations, 
$ (\bx_n)_{ n \geq 0} $  is  $d-$admissible, where $d=g+e_0$.\\
 (a) For $s \geq 2$,   $ m \geq 2^{2s+3} b^{d+t+s+1}(d+t)^{s+1}s^{2s} e(g+1)(e_0+s)\eta_1^{-s}$ and  \\$\eta_1 =(1+\deg((t_{s+1})_\infty))^{-1} $, we have
 \begin{equation}  \label{NiXi06}  
    1+  \min_{0 \leq Q <b^m} \min_{\bw \in E_m^{s}} \max_{1 \leq N \leq b^m} 
		N \emph{D}^{*}((\bx_{n\oplus Q} \oplus \bw)_{0 \leq  n < N}) \geq   
				2^{-2}	 b^{-d}K_{d,t,s+1}^{-s} \eta_1^s m^{s} .
\end{equation} 
(b) Let $s \geq 3$, $ m \geq 2^{2s+3} b^{d+t+s}(d+t)^{s}(s-1)^{2s-1} (g+e_0) e \eta_2^{-s+1}$, \\$e_s =\min_{1 \leq i \leq s} e_i$ and $\eta_2  =(1+\deg((t_s)_\infty))^{-1}$.  Then}
\begin{equation}\label{NiXi07}
      \min_{ \bw \in E_m^s}   b^m \emph{D}^{*}((\bx_n \oplus \bw)_{0 \leq  n < b^m}) \geq   
		2^{-2} b^{-d}K_{d,t,s}^{-s+1} \eta_2^{s-1} m^{s-1} .
\end{equation}\\


{\bf 2.5. Niederreiter-Xing sequence.} \\
Let $F/\FF_b$ be an algebraic function field with full constant field $\FF_b$ and
genus $g = g(F/\FF_b)$. 
 Assume that $F/\FF_b$ has at least $s+1$ rational places. Let $ P_1,...,P_{s+1}$ be $s+1$
distinct rational places of $F$. Let  $G_m =m(P_1+...+P_s) -(m-g+1)P_{s+1}$, and  
let $t_i$ be a local parameter at $P_i$, $1 \leq i \leq s+1$.
For any $f \in \cL(G_m)$ we have $\nu_{P_i}(f) \geq  −m$, and so the local
expansion of $f$ at $P_i$ has the form
\begin{equation}  \nonumber
f = \sum_{j=-m}^{\infty} f_{i,j} t^j_i , \quad {\rm with } \quad
  f_{i,j}  \in  \FF_b, \; j \geq -m, \; 1 \leq i \leq s.
\end{equation}
For $1 \leq i \leq s$, we define the $\FF_b$-linear map $\psi_{m,i}(f) \;:\; \cL(G_m) \to \FF^m_b$
by
\begin{equation} \nonumber
             \psi_{m,i}(f)= (   f_{i,- 1}, . . . , f_{i,- m}) \in \FF^m_b,\quad \for \quad
     f \in \cL(G_m).
\end{equation}
Let
\begin{equation}  \label{NiXi03}
        \cM_m=   \cM_m(P_1,...,P_s;G_m):= \{   (\psi_{m,1}(f), . . . , \psi_{m,s}(f))  \in \FF^{ms}_b\;  |  \;f \in \cL(G_m) \}.
\end{equation}

Let  $C^{(1)},...,C^{(s)} \in \FF_b^{\infty \times \infty}$  be the 
generating matrices of a digital sequence  $\bx_n(C)_{n \geq 0}$, and let $(\cC_m )_{m \geq 1}$ be 
the associated sequence of 
row spaces of overall generating matrices $[C]_m$, $m=1,2,...$ (see (\ref{Di08})).
  \\

{ \bf Theorem L.} (see \cite[Theorem 7.26 and Theorem 8.9]{DiPi})  {\it There exist matrices $C^{(1)},...,C^{(s)}$
such that  $\bx_n(C)_{n \geq 0}$ is a digital $(t, s)$-sequence  with
$t = g$ and $\cC_m =\cM_m^{\bot}(P_1,...,P_s; G_m)$ for $m \geq g+1$, $s \geq 2$.} \\

According to \cite[p.411]{DiNi} and \cite[p.275]{DiPi}, the construction of digital sequences of Niederreiter and
Xing \cite{NiXi}  can be achieved  by using the above approach. 
We propose the following way to get $\bx_n(C)_{n \geq 0}$. \\

We consider the $H$-differential  $dt_{s+1}$. Let $\omega$ be the corresponding
Weil differential,  $\div(\omega)$ the divisor of $\omega$, and $W:=\div(dt_{s+1})= \div(\omega)$.
By  (\ref{No14}), we have $ \deg(W)=2g-2$.
Similarly to (\ref{Hal03})-(\ref{Hal18a}),  we can construct a sequence $\dot{v}_0, \dot{v}_1, . . .$ of elements of $F$ such that $ \{\dot{v}_0, \dot{v}_1, . . . , \dot{v}_{\ell((m -g+1)P_{s+1} +W)-1} \}$
is a $\FF_b$-basis of \\$L_m:=\cL((m -g+1)P_{s+1} +W)$ and 
\begin{equation}\label{NiXi11a}
 \dot{v}_{r} \in L_{r+1} \setminus L_{r}, \; \quad  \nu_{P_{s+1}}(  \dot{v}_{r}) =-r+g-2,\; 
  r \geq g,\;\; \ad \;\; \dot{v}_{r,r+2-g} =1, \; \dot{v}_{r,j} =0
\end{equation}
for $2  \leq j <r +2- g$, where
\begin{equation} \nonumber
  \dot{v}_r   :=  \sum_{j \leq r-g+2}
 \dot{v}_{r,j}t_{s+1}^{-j}   \quad {\rm for}  \quad   
	\dot{v}_{r,j} \in \FF_b \;\; \ad \;\; r \geq g.
\end{equation}
According to Proposition A, we have that there exists $\tau_i \in F$
$(1 \leq i \leq s)$, such that $ \dd t_{s+1} = \tau_i \dd t_i   \quad \for  \quad 1 \leq i \leq s$.

Bearing in mind (\ref{No12}), (\ref{No16}) and (\ref{NiXi11a}),  we get
\begin{equation} \nonumber
   \nu_{P_i}( \dot{v}_r \tau_i) =\nu_{P_i}( \dot{v}_r \tau_i \dd t_i) =\nu_{P_i}( \dot{v}_r 
	\dd t_{s+1}) \geq
		\nu_{P_i}( \div(\dd t_{s+1}) - W ) =0, \quad  1 \leq i \leq s, \; r \geq 0.
\end{equation}
We consider the following local expansions
\begin{equation} \label{NiXi05}
  \dot{v}_r   \tau_i\; := \; \sum_{j = 0}^{\infty}
	\dot{c}^{(i)}_{j,r} t_i^{j}  , \quad {\rm where \; all }  \quad   
	\dot{c}^{(i)}_{j,r} \in \FF_b, \; 1 \leq i \leq s, \; j \geq 0.
\end{equation}
%
%
Now let $\dot{C}^{(i)} =(\dot{c}^{(i)}_{j,r}  )_{j,r \geq 0}$, $1 \leq i \leq s$,
 and let $\dot{\cC}_m$ be the 
row space of overall generating matrix $[\dot{C}]_m$ (see  (\ref{Di08})).\\

{\bf Theorem 5.} {\it  With the above notations, 
$\bx_n(\dot{C})_{n \geq 0}$ is a digital $d-$admissible $(t,s)$ sequence,
satisfying the bounds (\ref{NiXi06}) and (\ref{NiXi07}),  with $d=g+s$, 
$t=g$, and 
 $\dot{\cC}_m =\cM_m^{\bot}(P_1,...,P_s; G_m)$ for all $m \geq g+1$.} \\


{\bf 2.6 General  $d-$admissible digital  $(t,s)$ sequences}.
In \cite{KrLaPi},  discrepancy bounds for index-transformed uniformly distributed sequences was studied. In this subsection, we consider a lower bound of such a sequences.

 Let $s \geq 2$, $d \geq 1$, $t \geq 0$, $d_0=d+t$ and $m_k=s^2d_0(2^{2k+2}-1)$ for  $k=1,2,...$ .\\
Let  $C^{(s+1)} =(c^{(s+1)}_{i,j})_{ i,j \geq 1}$ be a $\NN \times \NN$ matrix over $ \FF_b$,
 and let 
$[C^{(s+1)}]_{m_k}$ be a nonsingular matrix, $k=1,2,...$ .
For $n \in[0,b^{m_k})$, let 
$\bh_k(n) =(h_{k,1}(n),...,h_{k,m_k}(n))= \bn [C^{(s+1)}]_{m_k}^{\top}$ 
and $h_k(n) =\sum_{j=1}^m \phi^{-1}(h_{k,j}(n)) b^{j-1} $ $(k \geq 1)$.
We have $h_k(l) \neq h_k(n)$ for $l \neq n$, $l,n \in[0,b^{m_k})$.
Let  $h_k^{-1}(h_k(n))=n$ for $n \in[0,b^{m_k})$.
It is easy to see that  $h_k^{-1}$   is a bijection from $[0,b^{m_k})$ to $[0,b^{m_k})$ $(k=1,2,...)$.
\\

{\bf Theorem 6.} {\it Let  $ (\bx_n)_{ n \geq 0} $  be a digital  $d-$admissible  $(t,s)$ sequence in base $b$. Then there exists a matrix $C^{(s+1)}$ and a sequence $(h^{-1} (n))_{n \geq 0}$ such that $[C^{(s+1)}]_{m_k}$ is nonsingular,
$h^{-1}(n) =h_l^{-1}(n)=h_k^{-1}(n)$ for $n \in [0,b^{m_k})$   $(l >k, \;k=1,2,...)$, 
   $ (\bx_{h^{-1} (n)})_{ n \geq 0} $  a $d-$admissible   $(t,s)$ sequence in base $b$,
and
}
\begin{equation}  \nonumber
    1+  \min_{0 \leq Q <b^{m_k},  \bw \in E_{m_k}^{s}} \max_{1 \leq N \leq b^{m_k}} 
		N \emph{D}^{*}((\bx_{h^{-1}(n)\oplus Q} \oplus \bw)_{0 \leq  n < N}) \geq   
					2^{-2} b^{-d}K_{d,t,s+1}^{-s} m_k^{s},\;\; k\geq 1\;.
\end{equation} \\

{\bf Remark 1.} 
Halton-type sequences were introduced in [Te1] for the case of rational
function fields over finite fields. 
Generalizations to the general case of algebraic function field  were obtained in [Le1] and [NiYe].
The constructions in [Le1] and [NiYe] are similar. The difference is that the construction in [NiYe] is more simple, but the construction in [Le1] a somewhat more general. \\

{\bf Remark 2.} We note that all explicit constructions of this article are expressed in terms of the residue of a differential and are similar to the Halton construction
  (see, e.g.,  (\ref{GeNi07}), (\ref{XiNi23}), (\ref{NiOz20}) and (\ref{Hal50})-(\ref{Hal68})). 
The earlier constructions of $(t,s)-$sequences using differentials, see e.g. [MaNi]. \\

{\bf Remark 3.}  A lattice  $ \Gamma \subset (\FF_b((x)))^{s+1}$ is  $d-${\sf admissible} if
\begin{equation} \nonumber 
\inf_{(\gamma_1,...,\gamma_{s+1}) \in \Gamma \setminus \{0\}} \prod_{i=1}^{s+1}b^{-\nu_{\infty}(\gamma_i)}  \geq b^{-d}.
\end{equation}
 A lattice  $ \Gamma \subset (\FF_b((x)))^{s+1}$ is said to be {\sf admissible} if   $ \Gamma $ 
is $d-${\sf admissible} with some $d>0$. 
In [Le1, Theorem 3.2], we proved the following analog of the main theorem of the duality theory (see, [DiPi, Section 7], [NiPi] and [Skr]): if a unimodular lattice  $\Gamma$ is $d-$admissible, then from the dual lattice $\Gamma^{\bot}$ we can get a $(t,s)$
 sequence with $t=d-s$. It is easy to prove (see, e.g., [Le5]) that 
 a lattice $\Gamma$ is admissible if and only if the dual lattice $\Gamma^{\bot}$ 
is admissible. In [Le4]  and in this paper we consider a more general object.
 We consider nets in $[0,1)^s$   having simultaneously  both $(t,m,s)$ properties and
 $d$-admissible properties. The $d$-admissible properties have a direct connection to the notion of the weight in the duality theory  (see  Definition~5, Definition 8 and Definition 9). Thus we can consider this paper as a part of the duality theory.



%
%
%
%
%
\section{Proof of theorems.}
\subsection{Generalized Niederreiter sequence. Proof of Theorem~1}

 Using  [Le4, Lemma 2] and \cite[Theorem 1]{Te3},  we obtain that $ (\bx_n)_{ n \geq 0}$ is $d-$admissible
  with $d=e_0$.
	
We apply Corollary 3 with
 $B_i^{'}=\emptyset$, $1 \leq i \leq s+1$, $B=0$,  $\hat{e}=e=e_1e_2\cdots e_s,$ $d_0 =d+t$,
 $\epsilon = \eta_1(2sd_0 e)^{-1}$ and $\eta_1 =s/(s+1)$.
In order to prove the first assertion in Theorem 1, it is sufficient to verify that
\begin{equation}\label{GeNi00}
    \Lambda_{1}   =\FF_b^{(s+1)d_0e[m \epsilon]},\qquad  \for \qquad   m \geq  9(d+t)es(s+1), 
\end{equation}
where
\begin{equation} \nonumber
    \Lambda_{1}  = \{ (y^{(1)}_{n,1},...,y^{(1)}_{n,d_1}, ..., y^{(s)}_{n,1},
		...,y^{(s)}_{n,d_s},
		\bar{a}_{d_{s+1,1}}(n),...,   \bar{a}_{d_{s+1,2}} (n) )  \; | \;  n \in [0, b^m)  \}
\end{equation}
with
\begin{equation}\label{GeNi01}
    d_{i}=\dot{m}_i = d_0 e [m \epsilon] \;\; (1 \leq i \leq s), \quad
  d_{s+1,1} =\ddot{m}_{s+1}+1:=t +(s-1)d_0e [m \epsilon], 
\end{equation}
$ d_{s+1,2}=\dot{m}_{s+1}: =t-1 + sd_0e [m \epsilon]$,   and  $n =\sum_{0 \leq j \leq m -1} a_j(n) b^{j}$.

Suppose that (\ref{GeNi00}) is not true. Then there exists $b_{i,j} \in \FF_b$ $(i,j \geq 1)$
 such that
\begin{equation}\label{GeNi02}
\sum_{i=1}^{s} \sum_{j=1}^{d_{i}} |b_{i,j}| +\sum_{j=d_{s+1,1}}^{d_{s+1,2}} |b_{s+1,j}| >0 
\end{equation}
and
\begin{equation}\label{GeNi04}
   \sum_{i=1}^s \sum_{j=1}^{d_{i}} b_{i,j} y_{n,j}^{(i)}  +
\sum_{j=d_{s+1,1}}^{d_{s+1,2}} b_{s+1,j} \bar{a}_j(n) =0 \quad {\rm for \; all  } \quad
  n \in [0,b^m).
\end{equation}
From (\ref{Ap301}) and (\ref{Ni00}), we have
\begin{equation}\nonumber
  y_{n,j}^{(i)}  = \sum_{r=0}^{m-1}   c_{j,r}^{(i)}  \bar{a}_{r}(n),
\end{equation}
with
\begin{equation}\label{GeNi05}
c^{(i)}_{ j,r} = a^{(i)}(Q + 1, k, r) \in \FF_b,  \qquad  j -1 = Qe_i + k, \quad
0 \leq k< e_i,
\end{equation}
$Q= Q(i,j)$, $k=k(i, j) $,    where $a^{(i)}(j, k, r)$ are defined from the expansions
\begin{equation}\nonumber
  \frac{ y_{i,j,k}(x)}{
p_i(x)^{j}} = \sum_{r\geq 0} a^{(i)} (j, k, r) x^{-r-1}.
\end{equation}
We consider the field $F=\FF_b(x)$, the valuation $\nu_{\infty}$ (see  (\ref{No00})) and the  place 
 $P_{\infty} =\div(x^{-1})$. 
By (\ref{No19}), we get
\begin{equation}\nonumber
            a^{(i)} (j, k, r) = \underset{P_{\infty},x^{-1}}\Res \big( y_{i,j,k}(x)
p_i(x)^{-j} x^{r+2}\big).
\end{equation}
Hence
\begin{equation}\label{GeNi07}
  y_{n,j}^{(i)}  = 
     \underset{P_{\infty},x^{-1}}\Res \Big(\frac{y_{i,Q(i,j)+1,k(i,j)}(x)}
{p_i(x)^{Q(i,j)+1}}   \sum_{r=0}^{m-1}  
		\bar{a}_{r}(n) x^{r+2}  \Big) = \underset{P_{\infty},x^{-1}}\Res \Big(\frac{y_{i,Q(i,j)+1,k(i,j)}(x)}
{p_i(x)^{Q(i,j)+1}}  n(x)  \Big)        
\end{equation}
with $n(x) = \sum_{j=0}^{m-1} \bar{a}_{j}(n) x^{j+2}$ $ \quad	\quad {\rm  for \; all } \;\; j \in [1,d_i], \; i \in [1,s] .$\\
We have $  \bar{a}_{j}(n) =  \underset{P_{\infty},x^{-1}}\Res (n(x) x^{-j-1} )  $.
From (\ref{GeNi04}),  we derive
\begin{equation}\label{GeNi08}
\underset{P_{\infty},x^{-1}}\Res(n(x)\alpha) =0  \; {\rm with  } \;
 \alpha =    \sum_{i=1}^s \sum_{j=1}^{d_{i}} b_{i,j}\frac{y_{i,Q(i,j)+1,k(i,j)}(x)}
{p_i(x)^{Q(i,j)+1}}  +
\sum_{j=d_{s+1,1}}^{d_{s+1,2}} b_{s+1,j} x^{-j-1}
\end{equation}
for all  $n \in [0,b^m)$.
Consider the local expansion
\begin{equation}\nonumber
          \alpha = \sum_{r= 0 }^{\infty} \varphi_r x^{-r-1} 
             \quad {\rm with   } \quad  \varphi_r \in\FF_b, \quad r \geq 0.
\end{equation}
Applying (\ref{No26}) and (\ref{GeNi08}), we derive
\begin{equation}\nonumber
 \underset{P_{\infty},x^{-1}}\Res (n(x)\alpha) =  \underset{P_{\infty},x^{-1}}\Res \Big( 
  	\sum_{\mu=0}^{m-1}  \bar{a}_{\mu}(n) x^{\mu +2} \;
\sum_{  r = 0}^{\infty} \varphi_r x^{-r-1} \Big)  = \sum_{\mu=0}^{m-1} \sum_{  r = 0}^{\infty}
    \bar{a}_{\mu}(n) 		\varphi_r
\end{equation}
\begin{equation}\nonumber
  \times   
\underset{P_{\infty},x^{-1}}\Res (x^{\mu+2-r-1}  )  =
\sum_{\mu=0}^{m-1}   \sum_{  r = 0}^{\infty}
    \bar{a}_{\mu}(n) 			\varphi_r \delta_{\mu,r}	
   =
 \sum_{\mu=0}^{m-1}   
   \bar{a}_{\mu}(n)  \varphi_{\mu} =0
\end{equation}
for all  $n \in [0,b^m)$.
Hence
\begin{equation} \label{GeNi10}
\varphi_r =0  
    \quad {\rm for  } \quad    r \in [0,m-1] \quad {\rm and} \quad \nu_{ \infty} (\alpha) \geq m.
\end{equation}
According to (\ref{GeNi05}), we obtain 
\begin{equation} \nonumber
    Q(i,j)+1 \leq Q(i,d_i)+1 \leq [(d_i-1)/e_i] +1 = d_i/e_i 
 \; {\for} \; j \in [1, d_i],  i \in [1, s].
\end{equation}
By (\ref{GeNi08}), we get
\begin{equation}\label{GeNi11a}
  \alpha \in \cL(G_1)  \quad {\rm with } \quad  G_1 = \sum_{i=1}^s d_{i}/e_i\div(p_i(x))   
					+  (d_{s+1,2}+1) \div(x)  -m P_{\infty}.
\end{equation}
From (\ref{GeNi00}) and (\ref{GeNi01}), we have for $ m \geq 2t+ 8(d+t)es(s+1)$ 
\begin{equation}\nonumber
\deg(G_1)   = \sum_{i=1}^s d_{i}
					+  d_{s+1,2} +1 -m = sd_0e [m \epsilon]  +t-1 + sd_0e [m \epsilon] +1 -m 
\end{equation}
\begin{equation}\nonumber
 \leq t-m(1 -2sd_0 e  \epsilon) = t -m(1 -\eta_1)=t -m/(s+1) <0.
\end{equation}
Hence $\alpha=0$.\\
Let g.c.d.$(x,p_{j}(x)) = 1$ for all $j \neq i$ with some $i \in [1,s]$. 
For example, let $i=1$, and let $p_{1}(x)  =x^{e_{1,1}}\dot{p}_{1}(x)$ with $e_{1,2}=\deg(\dot{p}_{1}(x))$,  $e_1= e_{1,1}   + e_{1,2}$, $e_{1,1} \geq 0$,  g.c.d.$(x,\dot{p}_{1}(x)) = 1$.
According to (\ref{GeNi08}), we get $\alpha = \alpha_1 + \alpha_2 + \alpha_3$, where
\begin{equation}\nonumber
 \alpha_1 = \sum_{i=2}^s \sum_{j=1}^{d_{i}} b_{i,j}\frac{y_{i,Q(i,j)+1,k(1,j)}(x)}
{p_i(x)^{Q(i,j)+1}}, \qquad  \alpha_2 =   \sum_{j=1}^{d_{1}} b_{1,j} \frac{\ddot{y}_{i,Q(1,j)+1,k(1,j)}(x)}
{\dot{p}_{1}(x)^{Q(1,j)+1}}  
\end{equation}
\begin{equation}\nonumber
\quad \ad \quad   \alpha_3 =   \sum_{j=1}^{d_{1}} b_{1,j} \frac{\dot{y}_{1,Q(1,j)+1,k(1,j)}(x)}
{x^{e_{1,1}(Q(1,j)+1)}} +
\sum_{j=d_{s+1,1}}^{d_{s+1,2}} \frac{b_{s+1,j}}{ x^{j+1}}
\end{equation}
with some polynomials $\dot{y}_{1,j,k}(x)$ and $\ddot{y}_{1,j,k}(x)$.

Using (\ref{GeNi01}), we obtain for $s\geq 2$ and $j \in [1,d_1]$ that
\begin{equation}\nonumber
d_{s+1,1}+1 =t+1 +(s-1)d_0e [m \epsilon]  >d_0e [m \epsilon]= d_{1} \geq 
  e_{1,1} d_{1} /e_1  \geq e_{1,1} \deg(Q(1,d_1)+1).
\end{equation}
We have that the polynomials  $p_2,...,p_s, \dot{p}_{1}$ and $x$ 
are  pairwise coprime over $\FF_b$.
By the uniqueness of the
partial fraction decomposition of a rational function, we have that $\alpha_3=0$ and $ b_{s+1,j} =0 $ for all $j \in [d_{s+1,1}, d_{s+1,2}]$.

Bearing in mind that  $p_1, . . . , p_{s}$ are  pairwise coprime polynomials
 over $\FF_b$, we obtain from  [Te3, p.242] or [Te2, p. 166,167] that 
 $ b_{i,j} =0 $ for all $j \in [1, d_{i}]$ and $i \in [1,s]$.\\
 By (\ref{GeNi02}), we have the contradiction.
Hence  assertion (\ref{GeNi00})  is  true. Thus the first assertion in Theorem 1 is proved.
  \\


%
Now consider the second assertion in Theorem 1:\\
Let, for example,  $i_0 =s$, i.e.
\begin{equation}\label{GeNi20}
  \min_{m/2 -t \leq je_s   \leq m, 0 \leq k < e_s}(1-\deg( y_{s,j,k}(x))j^{-1} e_s^{-1} ) \geq \eta_2.
\end{equation}
We apply Corollary 2 with   $\dot{s}=s \geq 3$,
 $B_i=\emptyset$, $1 \leq i \leq s$, $B=0$, $\tilde{r} =0$, $m=\tilde{m}$, $d_0 =d+t$, $\hat{e}=e=e_1e_2\cdots e_s$, 
 $\epsilon = \eta_2(2(s-1)d_0 e)^{-1}$.
In order to prove the second assertion in Theorem 1, it is sufficient to verify that
\begin{equation}\label{GeNi22}
    \Lambda_{2}  =\FF_b^{sd_0e[m \epsilon]} \qquad  \for \qquad  m \geq 
					8(d+t)e(s-1)^2 \eta_2^{-1} +2(1+t) \eta_2^{-1}(1-\eta_2)^{-1}, 
\end{equation}
where 
\begin{equation} \nonumber
    \Lambda_{2}  = \{ (y^{(1)}_{n,1},...,y^{(1)}_{n,d_1}, ..., y^{(s-1)}_{n,1},
		...,y^{(s-1)}_{n,d_{s-1}},
		y^{(s)}_{n,d_{s,1}},...,  y^{(s)}_{n,d_{s,2}} )  \; | \;  n \in [0, b^m)  \},
\end{equation}
with   
\begin{equation}\label{GeNi24}
 	d_{i} =\dot{m}_i = d_0 e [m \epsilon], \; i \in [1, s),\;\;  d_{s,1}=\ddot{m}_{s}+1: =m-t+1 -(s-1)d_0e [m \epsilon] 
\end{equation}
and  $d_{s,2} =\dot{m}_{s}: =m-t -(s-2)d_0 e [m \epsilon]$. \\
Suppose that (\ref{GeNi22}) is not true. Then there exists $b_{i,j} \in \FF_b$ $(i,j \geq 1)$
 such that
\begin{equation}\label{GeNi26}
\sum_{i=1}^{s-1} \sum_{j=1}^{d_{i}} |b_{i,j}| +\sum_{j=d_{s,1}}^{d_{s,2}} |b_{s,j}| >0 
\end{equation}
and
\begin{equation}\label{GeNi27}
   \sum_{i=1}^{s-1} \sum_{j=1}^{d_{i}} b_{i,j} y_{n,j}^{(i)}  +
\sum_{j=d_{s,1}}^{d_{s,2}} b_{s,j} y_{n,j}^{(s)} =0 \quad {\rm for \; all  } \quad
  n \in [0,b^m).
\end{equation}
Similarly to (\ref{GeNi08}), we have
\begin{equation} \nonumber
\underset{P_{\infty},x^{-1}}\Res (n(x)\alpha) =0  \quad \for \; {\rm all} \quad\; n \in [0,b^m),
 \quad  {\rm with  } \quad
 \alpha =     \alpha_1 + \alpha_2,
\end{equation}
where
\begin{equation}\label{GeNi32}
 \alpha_1 =    \sum_{i=1}^{s-1} \sum_{j=1}^{d_{i}} b_{i,j}\frac{y_{i,Q(i,j)+1,k(i,j)}(x)}
{p_i(x)^{Q(i,j)+1}}  \quad  \ad  \quad \alpha_2=
\sum_{j=d_{s,1}}^{d_{s,2}} b_{s,j} \frac{y_{s,Q(s,j)+1,k(s,j)}(x)}
   {p_s(x)^{Q(s,j)+1}}.
\end{equation}

Consider the local expansions
\begin{equation}\nonumber
          \alpha_1 = \sum_{r= 0 }^{\infty} \varphi_{1,r} x^{-r-1} \quad  \ad  \quad 
					\alpha_2 = \sum_{r= 0 }^{\infty} \varphi_{2,r} x^{-r-1} 
             \quad {\rm with   } \quad  \varphi_{i,r} \in\FF_b \quad i=1,2,\;r \geq 0.
\end{equation}
Analogously to (\ref{GeNi10}), we obtain from (\ref{GeNi27})
\begin{equation} \label{GeNi34}
\varphi_{1,r} + \varphi_{2,r} =0  
    \quad {\rm for \; all } \quad    r \in [0,m-1] .
\end{equation}
Taking into account that $j \leq (Q(s,j)+1)e_s$ and $ d_{s,1} \geq m/2 -t$, we get from  (\ref{No00}) and (\ref{GeNi20}) that
\begin{equation} \nonumber
  \nu_{\infty}\Big( \frac{y_{s,Q(s,j)+1,k(s,j)}(x)}
       {p_s(x)^{Q(s,j)+1}} \Big) = 
			(Q(s,j)+1)e_s- \deg(y_{s,Q(s,j)+1,k(s,j)}(x))=
\end{equation} 
\begin{equation} \nonumber
(Q(s,j)+1)\Big(1- \frac{\deg(y_{s,Q(s,j)+1,k(s,j)}(x))}{(Q(s,j)+1)e_s}\Big)e_s \geq	
			(Q(s,j)+1)e_s\eta_2
  \geq \eta_2 j, \quad  j \geq d_{s,1}. 
\end{equation} 
Applying (\ref{GeNi32})-(\ref{GeNi34}), we have $ \varphi_{2,r} =0  $ for $r < [\eta_2 d_{s,1}]$.
Therefore $ \varphi_{1,r} =0  \\$ for  $r < [\eta_2 d_{s,1}]$.
Hence
\begin{equation} \nonumber
\nu_{ \infty} (\alpha_1) \geq [\eta_2  d_{s,1}].
\end{equation}
Similarly to (\ref{GeNi11a}), we obtain 
\begin{equation}\nonumber
  \alpha_1 \in \cL(G_2)  \quad {\rm with } \quad  G_2 = \sum_{i=1}^{s-1} d_{i}/e_i\div(p_i(x))   
					  -[\eta_2  d_{s,1}] P_{\infty}.
\end{equation}
From (\ref{GeNi22}) and (\ref{GeNi24}), we have that $m > 2(1+t) \eta_2^{-1}(1-\eta_2)^{-1}$ and 
\begin{equation}\nonumber
\deg(G_2)   = \sum_{i=1}^{s-1} d_{i}
					   -[ d_{s,1}\eta_2] = (s-1)d_0e [m \epsilon]  -[(m-t+1 -(s-1)d_0e [m \epsilon]) \eta_2] 
\end{equation}
\begin{equation}\nonumber
 \leq (s-1)d_0e [m \epsilon] -(m-t -(s-1)d_0e [m \epsilon])\eta_2 +1
      =(1+\eta_2) (s-1)d_0e [m \epsilon]  
\end{equation}
\begin{equation}\nonumber
 - m \eta_2  +1 + t  \leq m ( (1+\eta_2)((s-1)d_0e\epsilon - \eta_2) +1 + t 
\end{equation}
\begin{equation}\nonumber
=m\eta_2((1+\eta_2)/2-1) +1+t =1+t-m\eta_2 (1-\eta_2)/2 <0.
\end{equation}
Hence $\alpha_1=0$  and $\varphi_{1,r} =0$ for $r \geq 0$.\\
Using [Te3, p.242] or [Te2, p. 166,167], we get 
 $ b_{i,j} =0 $ for all $j \in [1, d_{i}]$ and $i \in [1,s-1]$.\\
According to (\ref{GeNi34}), we have 
$ \varphi_{2,r} =0$  for  $r \in [0,m-1]$ .
Thus $\nu_{ \infty} (\alpha_2) \geq m.$  \\
From (\ref{GeNi32}), we obtain 
\begin{equation}\nonumber
  \alpha_2 \in \cL(G_3)  \quad {\rm with } \quad  G_3 =  [d_{s,2}/e_s +1]\div(p_s(x))   
					  -m  P_{\infty}.
\end{equation}
Applying (\ref{GeNi00}) and (\ref{GeNi01}), we derive for $m > 2/\epsilon$ and $s \geq 3$ 
\begin{equation}\nonumber
\deg(G_3)  \leq  m-t -(s-2)d_0 e [m \epsilon] +e_s-m <0 . 
\end{equation}
Hence $\alpha_2=0$.

By the uniqueness of the
partial fraction decomposition of a rational function, we have
from (\ref{GeNi32}) that $ b_{s+1,j} =0 $ for all $j \in [d_{s,1}, d_{s,2}]$.

 By (\ref{GeNi26}), we have a contradiction.
Thus  assertion (\ref{GeNi22})  is  true.
Therefore  Theorem 1  is proved. \qed   \\


\subsection{ Xing-Niederreiter sequence. Proof of Theorem 2}

{\bf Lemma 3.} {\it Let $P \in\PP_F$, $t$   be a local parameter of $P$ over
$F$, $k_j \in F$, $\nu_P(k_j) =j$ $(j=0,1,...)$. Then there exists $k^{\bot}_j \in F$ with $\nu_P(k^{\bot}_j) =-j$ $(j=1,2,...)$, such that}
\begin{equation}\label{Dop30}
      S_{-1} (t, k_{j_1}k^{\bot}_{j_2+1}) =\delta_{j_1,j_2}
   \qquad {\rm for } \qquad  j_1,j_2 \geq 0 .
\end{equation}

{\bf Proof.} Let $k^{\bot}_{1} =(tk_0)^{-1}$. We see 
  $\nu_P(k_j k^{\bot}_1 ) \geq 0$ for $j\geq 1$. 
Using (\ref{No06}) and (\ref{No26}), we get that (\ref{Dop30}) is true for $j_2=0$.
	Suppose that the assertion of the lemma is true for $0 \leq j_2 \leq  j_0-1$, $j_0 \geq 1$. We take
\begin{equation}\label{Dop32}
      k^{\bot}_{j_0+1} = \sum_{\mu=1}^{j_0} \rho_{\mu,j_0} k^{\bot}_{\mu}+ (tk_{j_0})^{-1},   \quad {\rm where } \quad  \rho_{\mu,j_0} = S_{-1} (t, k_{\mu-1}  (tk_{j_0})^{-1}).
\end{equation}
We see  that $ \nu_P(k^{\bot}_{j_0+1}) =-j_0-1$. 
By the condition of the lemma and the assumption of the induction, we have  $\nu_P(k_{j_1}k^{\bot}_{j_0+1}) \geq 0$  for $j_1 >j_0$ and
\begin{equation}\label{Dop34}
         S_{-1} (t,k_{j_1}k^{\bot}_{j_0+1}) =\delta_{j_1,j_0}
   \quad {\rm for } \quad  j_1 \geq j_0  .
\end{equation}

Now consider the case  $ j_1 \in[0, j_0) $.
Applying (\ref{Dop32}), we derive
\begin{equation}\nonumber
          S_{-1} (t,k_{j_1}k^{\bot}_{j_0+1}) =
					\sum_{\mu=1}^{j_0} \rho_{\mu,j_0}  S_{-1} (t, k_{j_1} k^{\bot}_{\mu}) +
					 S_{-1} (t, k_{j_1} (tk_{j_0})^{-1}) . 
\end{equation}
Using (\ref{No26}), (\ref{Dop32}) and the assumption of the induction, we get
\begin{equation}\nonumber
       S_{-1} (t,k_{j_1}k^{\bot}_{j_0+1}) =
					\sum_{\mu=1}^{j_0} \rho_{\mu,j_0} \delta_{j_1,\mu-1} +  S_{-1} (t, k_{j_1} (tk_{j_0})^{-1}) 
					= \rho_{j_1+1,j_0} -\rho_{j_1+1,j_0} =0.
\end{equation}
 Hence  (\ref{Dop34}) is true for  all $j_1 \geq 0$.
	By  induction, Lemma 3 is proved. \qed \\

{\bf Lemma 4}. {\it  $ (\bx_n)_{ n \geq 0} $ is  $d-$admissible with $d=g+e_0
$,
 where $e_0=e_1+...+e_s.$}  \\

{\bf Proof.} 
  Consider Definition 5.
	Taking into account that $ (\bx_n)_{ n \geq 0} $ is a digital sequence in base $b$, we can take $k=0$.
  Suppose that the assertion of the lemma is not true. 
By  (\ref{In04}), there exists $\tilde{n}>0$ such that	
	$\left\|\tilde{n}	\right\|_b \left\| \bx_{\tilde{n}}  \right\|_b  < b^{-d}=b^{-g-e_0}$. \\
Let	 $d_i = \dot{d}_i e_i +\ddot{d}_i $ with $0 \leq \ddot{d}_i < e_i$, 	$1 \leq i \leq s$,
 $\left\|\tilde{n}	\right\|_b =b^{m-1}$ and let
 $\left\| \bx^{(i)}_{ \tilde{n} }   \right\|_b =b^{- d_i-1} $,
 $1 \leq i \leq s$. 
Hence  $\tilde{n} \in [b^{m-1},b^m)$, $x^{(i)}_{\tilde{n},d_i+1} \neq 0$,
\begin{equation}\nonumber
          x^{(i)}_{\tilde{n},j} =0   \;\; {\rm  for \; all } \;\;  j \in [1,d_i],\;  i \in 
					[1,s] \quad {\rm  and } \quad \sum_{i=1}^s (d_i+1) -m 
					\geq d =g+e_0 .                 
\end{equation}
By (\ref{Ap301}), we have
\begin{equation}\label{XiNi16}
         y^{(i)}_{\tilde{n},j} =0    \quad {\rm  for \; all } \quad  j \in [1,\dot{d_i}e_i], \;\; i \in 
		[1,s] \quad {\rm  with } \quad \sum_{i=1}^s \dot{d_i}e_i  \geq m+g.                   
\end{equation}
 Let   
\begin{equation} \label{XiNi17}
\{ \dot{n}_0,..., \dot{n}_{g-1} \} =  \{ 0,1,..., 2g  \} \setminus \{ n_0,n_1,..., n_g  \}
\;\;\; \ad \;\;  \dot{n}_i =g+i+1 \; \for \; i \geq g.                   
\end{equation}
Let $n = \sum_{i=0}^{m-1} a_i(n) b^{i}$ with $a_i(n) \in Z_b$ $(i=0,1...)$, and let
$\bar{a}_i(n) =\phi(a_i(n))$ $(i=0,1,...)$ (see (\ref{Ap300})). 
From (\ref{Ap301}), (\ref{XiNi10}) and (\ref{XiNi12}), we get
\begin{equation}\label{XiNi18}
   y^{(i)}_{n,j}  = \sum_{\mu=0}^{m-1}  \bar{a}_{\mu}(n)  c_{j,\mu}^{(i)} \; = \;
      \sum_{\mu=0}^{m-1}  \bar{a}_{\mu}(n)  a_{j,\dot{n}_{\mu}}^{(i)} 
		\;\;	\quad {\rm  for } \;\; j \in [1,m], \; i \in 
					[1,s] .             
\end{equation}
By (\ref{XiNi08}), we have
\begin{equation}\label{XiNi18a}
 \nu_{P_{\infty}} (z_r) =r, \quad  \for \quad r \geq 0, \quad \ad \quad 
               z_{n_u}=w_u
 \quad  \with \quad u=0,1,...,g.
\end{equation}
 Using Lemma 3,   (\ref{No06}) and (\ref{No19}), we obtain that there exists a sequence $(z^{\bot}_{j})_{j \geq 1}$ such that $\nu_{P_{\infty}}(z^{\bot}_{j}) =-j$ and
\begin{equation}\label{XiNi19}
 \underset{P_{\infty},z}\Res (z_i z^{\bot}_{j+1}) =  S_{-1} (z,  z_i z^{\bot}_{j+1} ) =   \delta_{i,j}  \quad  {\rm for \;all}
   \quad  i,j \geq 0.
\end{equation}
We put
\begin{equation}\label{XiNi20}
           f_{n} = \sum_{\mu=0}^{m-1}  \bar{a}_{\mu}(n) z_{\dot{n}_{\mu} +1}^{\bot}.
\end{equation}
Hence
\begin{equation}\label{XiNi21}
    \bar{a}_{\mu}(n)=
 \underset{P_{\infty},z}\Res (f_n z_{\dot{n}_{\mu}} )  \quad  {\rm for }
   \quad  0 \leq \mu \leq m-1,\;  n \in [0,b^m) .
\end{equation}
By (\ref{No26}) and (\ref{XiNi17}), we have  $\delta_{\dot{n}_{\mu},n_u} =0$ for all 
$0 \leq u \leq g, \; \mu \geq 0$.\\
Applying (\ref{XiNi18a}) and (\ref{XiNi19}),  we derive 
\begin{equation}\label{XiNi22}
        \underset{P_{\infty},z}\Res(f_n w_u) = \underset{P_{\infty},z}\Res\Big( 
  	\sum_{\mu=0}^{m-1} \bar{a}_{\mu}(n) z_{\dot{n}_{\mu} +1}^{\bot} \; z_{n_u}\Big) 
\end{equation}
\begin{equation}\nonumber
         =\sum_{\mu=0}^{m-1} \bar{a}_{\mu}(n)  \underset{P_{\infty},z}\Res\big( 
  	z_{\dot{n}_{\mu} +1}^{\bot} \;z_{n_u}\big) =
							\sum_{\mu=0}^{m-1} \bar{a}_{\mu}(n)\delta_{\dot{n}_{\mu},n_u}			
										=0 \quad {\rm for} \quad  u=0,1,...,g, \; n\geq 0.   
\end{equation}
According to (\ref{XiNi10}) and (\ref{XiNi20}), we have
\begin{equation}\nonumber
    \underset{P_{\infty},z}\Res (f_{n} k_{i,j}) =  \underset{P_{\infty},z}\Res \Big( 
  	\sum_{\mu=0}^{m-1}  \bar{a}_{\mu}(n) z_{\dot{n}_{\mu} +1}^{\bot} \;
\sum_{  r = 0}^{\infty}  a_{j,r}^{(i)} z_{r} \Big)
\end{equation}
\begin{equation}\nonumber
   = \sum_{\mu=0}^{m-1}  \sum_{  r = 0}^{\infty}
    \bar{a}_{\mu}(n) 			a_{j,r}^{(i)}
  \underset{P_{\infty},z}\Res ( z_{\dot{n}_{\mu} +1}^{\bot} \; z_{r}  )  =
\sum_{\mu=0}^{m-1}   \sum_{  r = 0}^{\infty}
    \bar{a}_{\mu}(n) 			a_{j,r}^{(i)} \delta_{\dot{n}_{\mu},r}	
   =
 \sum_{\mu=0}^{m-1}  
   \bar{a}_{\mu}(n)  a_{j,\dot{n}_{\mu}}^{(i)} .
\end{equation}
From (\ref{XiNi18}), we get 
\begin{equation}\label{XiNi23}
  \underset{P_{\infty},z}\Res(f_{n} k_{i,j}) =y_{n,j}^{(i)} 
			\quad {\rm  for \; all } \quad  j \in [1,m],   \; i \in [1,s], \; n \in [0,b^m).          
\end{equation}
Using (\ref{XiNi16}) and (\ref{XiNi22}), we derive
\begin{equation}\nonumber
           \underset{P_{\infty},z}\Res \Big(f_{\tilde{n}} \Big(\sum_{r=0}^g b_r w_r +   
					\sum_{i=1}^s \sum_{j=1}^{\dot{d}_i e_i} b_{i,j} k_{i,j}\Big)\Big) =  0   \quad {\rm for \; all } \quad b_i, b_{i,j} \in \FF_b. 
\end{equation}
Taking into account that $(w_0,...,w_g, k_{1,1},...k_{1,\dot{d}_1e_1}, ...,k_{s,1},...,
k_{s,\dot{d}_se_s})$   is the basis of $\cL(G+  \sum_{i=1}^s\dot{d}_i P_i) $
 (see (\ref{XiNi04})), we obtain 
\begin{equation}\label{XiNi24}
          \underset{P_{\infty},z}\Res(f_{\tilde{n}} \gamma)  = 0    \quad {\rm for \; all } \quad 
				\gamma \in \cL(\dot{G}) \quad {\rm with} \quad \dot{G} =G+  \sum_{i=1}^s \dot{d}_i P_i.
\end{equation}
By (\ref{XiNi16}), we have
\begin{equation}\nonumber
          \deg(\dot{G} - (m+g+1)P_{\infty}) =2g +\sum_{i=1}^s \dot{d}_ie_i -(m+g+1) \geq
					2g +m+g -(m+g+1) =2g-1.
\end{equation}
Using the Riemann-Roch theorem, we get 
\begin{equation}\nonumber
    \ddot{G}=(\dot{G} -(m+g)P_{\infty}) \setminus (\dot{G} -(m+g+1)P_{\infty}) \neq \emptyset.
\end{equation}
We take $v \in \ddot{G}$. Hence $\nu_{P_{\infty}}(v) =m+g$.

From (\ref{XiNi08}), we derive   
  $v = \sum_{r \geq m+g} \hat{b}_{r} z_r$ with some 
$\hat{b}_{r} \in \FF_b$  $(r \geq m+g)$ and  $\hat{b}_{m+g} \neq 0$.
According to (\ref{XiNi17}), we have
 $\dot{n}_{m-1} = m+g$.
Therefore $v = \sum_{r \geq \dot{n}_{m-1}} \hat{b}_{r} z_r$.

Taking into account that $\tilde{n} \in [b^{m-1},b^m)$, we get $a_{m-1}(\tilde{n}) \neq 0$.\\ 
By (\ref{XiNi19}), (\ref{XiNi20}) and(\ref{XiNi24}), we obtain 
\begin{equation}\nonumber
     0=     \underset{P_{\infty},z}\Res (f_{\tilde{n}} v) =
				\sum_{\mu=0}^{m-1} \sum_{r \geq {\dot{n}_{m-1}}} 
				a_{\mu}(\tilde{n}) \hat{b}_r
			 \underset{P_{\infty},z}\Res (z_{\dot{n}_{\mu} +1}^{\bot} \;z_r)
				=  				\sum_{\mu=0}^{m-1} \sum_{r \geq \dot{n}_{m-1}} 
				a_{\mu}(\tilde{n}) \hat{b}_r
				\delta_{\dot{n}_{\mu},r}.
\end{equation}
Bearing in mind that $\delta_{\dot{n}_{\mu},r}=1$ for $\mu \in [0,m-1]$, $r \geq \dot{n}_{m-1}$   
if and only if $\mu=m-1$ and $r=\dot{n}_{m-1}$ (see (\ref{XiNi17})), we get $\Res_{P_{\infty},z} (f_{\tilde{n}} v)=
 a_{m-1}(\tilde{n}) \hat{b}_{\dot{n}_{m-1}} \neq 0$.
We have a contradiction. 
Hence Lemma 4  is proved. \qed  \\

{\bf Lemma 5.} {\it Let $s \geq 2$, $d_i=d_0e [m\epsilon],\; 1 \leq i \leq s$, $  d_{s+1,1} =t +(s-1)d_0e [m \epsilon] $, $ 	d_{s+1,2} =t-1 + sd_0e [m \epsilon]$, $d_0=d+t$, $t=g+e_0-s$, $e =e_1...e_s$
 and $m \geq 2/\epsilon$. Then the system $\{w_0,w_1, . . . ,w_g\}\cup\{z^{j+g+1}\}_{d_{s+1,1}\leq j \leq d_{s+1,2}}$ $ \cup
\{k_{i,j }\}_{1 \leq i \leq s, 1 \leq j \leq d_i}$ of elements of
 $F$ is linearly independent over $\FF_b$.}  \\

{\bf Proof.} 
Suppose that
\begin{equation} \nonumber
 \alpha := \sum_{j=0}^{g}  b_{0,j} w_{j}+ \sum_{i=1}^s \sum_{j=1}^{d_{i}} b_{i,j} k_{i,j} 
 + \sum_{j=d_{s+1,1}}^{d_{s+1,2}} b_{s+1,j}z^{j+g+1} =0
\end{equation}
 for some $b_{i,j} \in \FF_b$ and $ \sum_{j=0}^{g}|b_{0,j}| + \sum_{i=1}^s \sum_{j=1}^{d_{i}}|b_{i,j}| +	 \sum_{j=d_{s+1,1}}^{d_{s+1,2}} |b_{s+1,j}| >0$. 
Let 	
\begin{equation}\label{XiNi31}
\beta_1 = \sum_{j=0}^{g}  b_{0,j} w_{j}, \;   \beta_{2,i} =\sum_{j=1}^{d_{i}} b_{i,j} k_{i,j}, \;\; \beta_2 =\sum_{i=1}^s \beta_{2,i}, \; 
   \beta_3=\sum_{j=d_{s+1,1}}^{d_{s+1,2}} b_{s+1,j}z^{j+g+1}.
\end{equation}
We have
\begin{equation}\label{XiNi33}
 \alpha =  \beta_1 +\beta_2  +\beta_3 = 0.
\end{equation}
Suppose that$\sum_{i=1}^s \sum_{j=1}^{d_{i}}|b_{i,j}|=0$ and $\alpha=0$. 
By (\ref{XiNi31}) and (\ref{XiNi33}), we have $\beta_1 + \beta_3=0$ and 
$\nu_{P_{\infty}} (\beta_1 ) \geq d_{s+1,1}$.
Taking into account that $\beta_1  \in \cL(G)$ with $\deg(G) = 2g$, we obtain from the Riemann-Roch theorem that $\beta_1  =0$.
Therefore $ \sum_{j=0}^{g} |b_{0,j}|=0$ and $ \sum_{j=d_{s+1,1}}^{d_{s+1,2}} |b_{s+1,j}|=0$.
We have a contradiction.

 According to [DiPi, Lemma 8.10], we get that if  $ \sum_{j=d_{s+1,1}}^{d_{s+1,2}} |b_{s+1,j}| =0$ and $\alpha=0$, then 
	$ \sum_{j=0}^{g} |b_{0,j}|=0$ and $\sum_{i=1}^s \sum_{j=1}^{d_{i}}|b_{i,j}|=0$.  So, we will consider only the case then 
	$\sum_{i=1}^s \sum_{j=1}^{d_{i}}|b_{i,j}| >0$
  	and  $ \sum_{j=d_{s+1,1}}^{d_{s+1,2}} |b_{s+1,j}| >0$.

Let $\sum_{j=1}^{d_{h}}|b_{h,j}| >0$ for some  $h \in [1,s]$, and let $\nu_{P_h} (z) \geq 0$. \\
By the construction of  $k_{h,j}$,
 we have $\beta_{2,h} \notin \cL(G)$ and $\beta_{2,h} \neq 0$. 
Applying (\ref{XiNi06}) and (\ref{XiNi31}),
we obtain $\nu_{P} (\beta_{2,h}) \geq -\nu_{P} (G) $ for any place $P \neq P_h$ and hence we obtain that $\nu_{P_h} (\beta_{2,h}) \leq -\nu_{P_h} (G) -1$ with $\nu_{P_h} (G) \geq 0$.

On the other hand, using (\ref{XiNi06}) (\ref{XiNi31})  and (\ref{XiNi33}),  we get
\begin{equation}\nonumber
\nu_{P_h} (\beta_{2,h}) = \nu_{P_h} \Big(-\beta_{1}   - \sum_{i=1, i\neq h}^s \beta_{2,i}
  -\beta_{3} \Big) \qquad \qquad \qquad
\end{equation}
\begin{equation}\nonumber
 \qquad  \geq \min\Big(  \nu_{P_h} (\beta_{1}),\nu_{P_h} (\beta_{3}), 
	  \min_{1 \leq i \leq s, i \neq h}  \nu_{P_h} (\beta_{2,i})\Big) \geq -\nu_{P_h} (G).
\end{equation}
We have a contradiction.

Now let $\nu_{P_h} (z) \leq -1$.
 Bearing in mind that $ \sum_{j=d_{s+1,1}}^{d_{s+1,2}} |b_{s+1,j}| >0$, we obtain that $\beta_3 \neq 0$, and  $\nu_{P_h} (\beta_3) \leq -d_{s+1,1}-g-1$. 
On the other hand, using (\ref{XiNi06}) and (\ref{XiNi33}), we have
\begin{equation}\nonumber
\nu_{P_h} (\beta_{3}) = \nu_{P_h} (\beta_{1}  +\beta_{2})
  \geq -\nu_{P_h} (G) -[(d_h -1)/e_h +1]e_h \geq -2g -d_h.
\end{equation}
Taking into account that
\begin{equation}\nonumber
    d_{s+1,1} +g+1 -(2g+d_{h}) =t+g+1 +(s-2)d_0e [m \epsilon]-2g \geq t-g+1 \geq 1,
\end{equation} 
we have  a contradiction. Thus Lemma 5 is proved. \qed  \\

{\bf Lemma 6.} {\it Let  $s \geq 2$, $d_0 =d+t,$ $t=g+e_0-s$, $\epsilon = \eta_1(2sd_0 e)^{-1}$,
 $\eta_1=(1+ \deg((z)_{\infty}))^{-1}$,
\begin{equation} \nonumber
   \Lambda_{1} := \{ (y^{(1)}_{n,1},...,y^{(1)}_{n,d_{1}}, ...,
		y^{(s)}_{n,1},
		...,y^{(s)}_{n,d_{s}},
		\bar{a}_{d_{s+1,1}}(n),...,   \bar{a}_{d_{s+1,2}}(n)  )  \; |  \; n \in [0, b^m)           \},
\end{equation}
where 			
\begin{equation} \label{XiNi39}
  d_{i} =\ddot{m}_{i}:= d_0 e [m \epsilon] \;\; (1 \leq i \leq s), \quad
	d_{s+1,1} =\ddot{m}_{s+1}+1:=t +(s-1)d_0e [m \epsilon] ,  
\end{equation} 
$	d_{s+1,2} =\dot{m}_{s+1}:=t-1 + sd_0e [m \epsilon] $, $e=e_1e_2\cdots e_s$, 
 and  $n =\sum_{0 \leq j \leq m-1} a_j(n) b^{j}$.
 Then}
\begin{equation}\label{XiNi38}
    \Lambda_{1}  =\FF_b^{(s+1)d_0e[m \epsilon]}, \qquad  \with \quad  
     m \geq 9(d+t)es^2\eta_1^{-1}.
\end{equation}
\\

{\bf Proof.} Suppose that (\ref{XiNi38}) is not true. Then there exists $b_{i,j} \in \FF_b$ $(i,j \geq 1)$
 such that
\begin{equation}\label{XiNi40}
\sum_{i=1}^{s} \sum_{j=1}^{d_{i}} |b_{i,j}| +\sum_{j=d_{s+1,1}}^{d_{s+1,2}} |b_{s+1,j}| >0 
\end{equation}
and
\begin{equation}\label{XiNi41}
   \sum_{i=1}^s \sum_{j=1}^{d_{i}} b_{i,j} y_{n,j}^{(i)}  +
\sum_{j=d_{s+1,1}}^{d_{s+1,2}} b_{s+1,j} \bar{a}_j(n) =0 \quad {\rm for \; all  } \quad
  n \in [0,b^m).
\end{equation}
From (\ref{XiNi21}) and (\ref{XiNi23}),  we obtain for $n \in [0,b^m)$
\begin{equation}\nonumber
  \bar{a}_{j-1}(n) =  \underset{P_{\infty},z}\Res(f_n z_{\dot{n}_{j-1}}) \quad \ad  \quad
   y_{n,j}^{(i)} =  \underset{P_{\infty},z}\Res (f_n k_{i,j}) \quad \with \quad
	j \in [1,m], \; i \in [1,s]. 
\end{equation}
Applying (\ref{XiNi08}) and (\ref{XiNi17}), we get $\dot{n}_{j-1} =g+j$ and $z_{\dot{n}_{j-1}} =z^{g+j} $ for $j \geq d_{s+1,1}$.	
Hence	
\begin{equation}\label{XiNi44a}
   \sum_{i=1}^s \sum_{j=1}^{d_{i}} b_{i,j}  \underset{P_{\infty},z}\Res(f_n k_{i,j}) +
\sum_{j=d_{s+1,1}}^{d_{s+1,2}} b_{s+1,j}  \underset{P_{\infty},z}\Res (f_n z^{g+j+1})
 = \underset{P_{\infty},z}\Res (f_n \alpha_1) =0  
\end{equation}
with
\begin{equation}\label{XiNi44}
 \alpha_1 =  \sum_{i=1}^s \sum_{j=1}^{d_{i}} b_{i,j} k_{i,j} 
 + \sum_{j=d_{s+1,1}}^{d_{s+1,2}} b_{s+1,j}z^{g+j+1} \quad \for \quad n \in [0,b^m).
\end{equation}
Let 
\begin{equation}\nonumber
    b_{0,u} = -\sum_{i=1}^s \sum_{j=1}^{d_{i}} b_{i,j}  a^{(i)}_{j,n_u},
 \quad \beta_1 = \sum_{u=0}^g  b_{0,u}  w_u, \quad
     \beta_2 =\sum_{i=1}^s \sum_{j=1}^{d_{i}} b_{i,j} k_{i,j}, 
\end{equation}
\begin{equation} \label{XiNi46}
\beta_3 =  \sum_{j=d_{s+1,1}}^{d_{s+1,2}} b_{s+1,j}z^{g+j+1}
\quad \ad \quad  \alpha_2  = \beta_1 +\beta_2  +\beta_3 = \beta_1 +  \alpha_1.
\end{equation}
By (\ref{XiNi40}) and Lemma 5, we get
\begin{equation} \label{XiNi48}
 \alpha_2 \neq 0.
\end{equation}
Consider the local expansion
\begin{equation}\label{XiNi50}
          \alpha_2 = \sum_{r= 0 }^{\infty} \varphi_r z_r 
             \quad {\rm with   } \quad  \varphi_r \in\FF_b, \quad r \geq 0.
\end{equation}
Using (\ref{XiNi08}), (\ref{XiNi10}) and (\ref{XiNi46}), we have
\begin{equation}\label{XiNi51}
         \varphi_{n_u} =0 
             \quad {\rm for   } \quad  0 \leq u \leq g.
\end{equation}
From (\ref{XiNi22}), we derive  $ \underset{P_{\infty},z}\Res(f_n w_u) =0 $ 
$(0 \leq u \leq g)$. By (\ref{XiNi44a}) and (\ref{XiNi46}), we get
\begin{equation} \nonumber
  \underset{P_{\infty},z}\Res (f_n \beta_1) =0  
    \quad {\rm and  } \quad  \underset{P_{\infty},z}\Res (f_n \alpha_2) =0 
		  \quad {\rm for \; all  } \quad   n \in [0,b^m).
\end{equation}
Applying (\ref{XiNi19}), (\ref{XiNi20}) and (\ref{XiNi50}), we obtain
\begin{equation}\nonumber
    \underset{P_{\infty},z}\Res (f_{n} \alpha_2) =  \underset{P_{\infty},z}\Res\Big( 
  	\sum_{\mu=0}^{m-1}  \bar{a}_{\mu}(n) z_{\dot{n}_{\mu} +1}^{\bot} \;
\sum_{  r = 0}^{\infty} \varphi_r  z_{r} \Big) \qquad\qquad\qquad\qquad\qquad\qquad
\end{equation}
\begin{equation}\nonumber
   = \sum_{\mu=0}^{m-1}  \sum_{  r = 0}^{\infty}
    \bar{a}_{\mu}(n) 		\varphi_r
  \underset{P_{\infty},z}\Res ( z_{\dot{n}_{\mu} +1}^{\bot} \; z_{r}  )  =
\sum_{\mu=0}^{m-1}   \sum_{  r = 0}^{\infty}
    \bar{a}_{\mu}(n) 			\varphi_r \delta_{\dot{n}_{\mu},r}	
   =
\sum_{\mu=0}^{m-1}  
   \bar{a}_{\mu}(n)  \varphi_{\dot{n}_{\mu}} =0
\end{equation}
for all  $n \in [0,b^m)$. \\
Hence $\varphi_{\dot{n}_{\mu}} =0$ for $\mu \in [0,m-1]$.
According to (\ref{XiNi17}) and (\ref{XiNi51}), we have 
\begin{equation} \label{XiNi60}
\varphi_r =0  
    \quad {\rm for  } \quad    r \in [0,m+g].
\end{equation}
Therefore
\begin{equation}\label{XiNi62}
          \nu_{P_{\infty}} (\alpha_2) > m+g.
\end{equation}
From (\ref{XiNi06}) and (\ref{XiNi46}), we derive 
\begin{equation} \nonumber
           \beta_1 + \beta_2 \in \cL\big(G+\sum_{i=1}^s [(d_i-1)/e_i+1]P_i\big)   \quad {\rm and } \quad 
					\beta_3  \in \cL\big( ( d_{s+1,2} + g+1  )(z)_{\infty}\big).
\end{equation}
By (\ref{XiNi62}), we obtain 
\begin{equation}\nonumber
\alpha_2 \in \cL(G_1)  \;\; {\rm with } \;\;  G_1 = G+\sum_{i=1}^s [(d_i-1)/e_i+1] P_i   
					+  ( d_{s+1,2} + g+1  ) (z)_{\infty} -(m+g+1) P_{\infty}.
\end{equation}
Using (\ref{XiNi39}), we have 
\begin{equation}\nonumber
\deg(G_1)   = 2g + \sum_{i=1}^s d_{i}   
					+  ( d_{s+1,2} + g+1  ) \deg((z)_{\infty}) -(m+g+1) 
\end{equation}
\begin{equation}\nonumber
 = 2g +sd_0e [m \epsilon]  +(t+g + sd_0e [m \epsilon]) (\eta_1^{-1} -1) -(m+g+1)
\end{equation}
\begin{equation}\nonumber
   \leq 2g +(t+g) (\eta_1^{-1} -1) 
  +sd_0e m \epsilon \eta_1^{-1} -(m+g+1)
\end{equation}
\begin{equation}\nonumber
= g-1 +(t+g) (\eta_1^{-1} -1) 
	-m(1-sd_0e \epsilon \eta_1^{-1}) = g-1 +(t+g) (\eta_1^{-1} -1)   -m/2 <0 
\end{equation}
for $ m \geq 9(d+t)es^2\eta_1^{-1}> 2(g-1)+   2(t+g) (\eta_1^{-1}-1) $ and $d=g+e_0$.  
Hence $\alpha_2=0$. By (\ref{XiNi48}), we have a contradiction. 
Therefore  assertion (\ref{XiNi41})  is not true. Thus Lemma 6 is proved. \qed
\\


{\bf End of the proof of Theorem 2.}
Using Lemma 4 and Theorem I, we get that $ (\bx(n))_{ n \geq 0}$ is a $d-$admissible
 digital $(t, s) $  sequence with $d=g+e_0$ and $t=g+e_0-s$.
%
Applying Lemma 6 and Corollary 3  with
 $B_i^{'}=\emptyset$, $1 \leq i \leq s+1$, $B=0$ and  $\hat{e}=e=e_1e_2\cdots e_s$, we get the first assertion in Theorem 2.\\

%
Consider the second assertion in Theorem 2 :\\
Let, for example,  $i_0 =s$, i.e.
\begin{equation}\label{XiNi70}
 \nu_{P_\infty}\big( k_{s,j}  \big) \geq \eta_2 j \quad \for \quad j \geq  m/2-t,
 \quad \ad \quad \eta_2 \in (0,1). 
\end{equation}

From  (\ref{In04}), Lemma 4 and Theorem I, we get that $ (\bx(n))_{0 \leq n < b^m}$ is a $d-$admissible
 digital $(t, m,s) $  net with $d=g+e_0$ and $t=g+e_0-s$.\\

We apply Corollary 2 with $\dot{s} =s \geq 3$, 
 $B_i=\emptyset$, $1 \leq i \leq s$, $B=0$, $\tilde{r} =0$, $m=\tilde{m}$, $\hat{e}=e=e_1e_2\cdots e_s,$ $d_0 =d+t$,
 $t=g+e_0-s$ and $e_0=e_1+...+e_s$.
In order to prove the second assertion in Theorem 2, it is sufficient to verify that

\begin{equation}\label{XiNi72}
   \Lambda_{2}  =\FF_b^{sd_0e[m \epsilon]} \quad  \for \quad  m \geq  
	  8(d+t)e(s-1)^2\eta_2^{-1} + 2(1 +2g+ \eta_2t) \eta_2^{-1}  (1-\eta_2)^{-1} ,
\end{equation}
where 
\begin{equation} \nonumber
    \Lambda_{2} = \{ (y^{(1)}_{n,1},...,y^{(1)}_{n,d_1}, ..., y^{(s-1)}_{n,1},
		...,y^{(s-1)}_{n,d_{s-1}},
		y^{(s)}_{n,d_{s,1}},...,  y^{(s)}_{n,d_{s,2}})  \; | \;  n \in [0, b^m)  \}
\end{equation}
with   
\begin{equation}\label{XiNi72a}
   d_{i} =\dot{m}_i:= d_0 e [m \epsilon],  \; i \in [1, s), \;\;
 	d_{s,1} =\ddot{m}_{s}+1:=m-t+1 -(s-1)d_0e [m \epsilon] , 	
\end{equation}
$d_{s,2} =\dot{m}_{s}:=m-t -(s-2)d_0 e [m \epsilon] $, and $ \epsilon = \eta_2(2(s-1)d_0 e)^{-1}$.

Suppose that (\ref{XiNi72}) is not true. Then there exists $b_{i,j} \in \FF_b$ $(i,j \geq 1)$
 such that
\begin{equation}\label{XiNi76}
\sum_{i=1}^{s-1} \sum_{j=1}^{d_{i}} |b_{i,j}| +\sum_{j=d_{s,1}}^{d_{s,2}} |b_{s,j}| >0 
\end{equation}
and
\begin{equation}\nonumber
   \sum_{i=1}^{s-1} \sum_{j=1}^{d_{i}} b_{i,j} y_{n,j}^{(i)}  +
\sum_{j=d_{s,1}}^{d_{s,2}} b_{s,j} y_{n,j}^{(s)} =0 \quad {\rm for \; all  } \quad
  n \in [0,b^m).
\end{equation}
Similarly to (\ref{XiNi44a}), we get
\begin{equation} \nonumber
 \underset{P_{\infty},z}\Res (f_n \alpha_1) =0  
 \quad \for \; {\rm all} \quad\; n \in [0,b^m),
   {\rm with  } \;\;\;
 \alpha_1 =     \alpha_2 - \beta_1
\end{equation}
where $ \alpha_2 =     \beta_1 +     \beta_2 + \beta_3$,  with
\begin{equation}\label{XiNi78}
    \beta_1 = \sum_{u=0}^g  b_{0,u}  w_u,\quad \beta_2 =    \sum_{i=1}^{s-1} \sum_{j=1}^{d_{i}} b_{i,j} k_{i,j}  \quad 
	\ad  \quad  \beta_3=
    \sum_{j=d_{s,1}}^{d_{s,2}} b_{s,j} k_{s,j}
\end{equation}
and $ b_{0,u} = -\sum_{i=1}^{s-1} \sum_{j=1}^{d_{i}} b_{i,j}  a^{(i)}_{j,n_u}
 - \sum_{j=d_{s_1}}^{d_{s_2}} b_{s,j}  a^{(s)}_{j,n_u} $.
Consider the local expansions
\begin{equation}\nonumber
          \beta_1 +     \beta_2 = \sum_{r= 0 }^{\infty} \dot{\varphi}_r z_r \quad  \ad  \quad 
					 \beta_3 = \sum_{r= 0 }^{\infty} \ddot{\varphi}_r z_r 
             \quad {\rm with   } \quad  \varphi_{i,r} \in\FF_b \quad i=1,2,\;r \geq 0.
\end{equation}
Analogously to (\ref{XiNi60}), we obtain
\begin{equation} \label{XiNi80}
\dot{\varphi}_r + \ddot{\varphi}_r =0  
    \quad {\rm for  } \quad    r \in [0,m+g] .
\end{equation}
Using  (\ref{XiNi70}), (\ref{XiNi72a}) and (\ref{XiNi78}), we get
\begin{equation} \nonumber
  \nu_{P_\infty}\big( k_{s,j}  \big) \geq \eta_2 j \quad \for \quad 
	j \geq d_{s,1} \geq m/2-t,
	\quad  \ad \quad \ddot{\varphi}_r =0  \; \for \;r \leq [\eta_2 d_{s,1}]-1.
\end{equation}    
Therefore $ \dot{\varphi}_r =0  $ for  $r \leq [\eta_2 d_{s,1}]-1$.
Hence
\begin{equation} \nonumber
\nu_{P_{\infty}} (\beta_1 +\beta_2) \geq [\eta_2  d_{s,1}].
\end{equation}
By (\ref{XiNi78}), we obtain 
\begin{equation}\nonumber
 \beta_1 +     \beta_2 \in \cL(G_2)  \quad {\rm with } \quad  
	G_2 = G+ \sum_{i=1}^{s-1} [(d_{i} -1)/e_i +1]P_i   
					  -[\eta_2  d_{s,1}] P_{\infty}.
\end{equation}
According to (\ref{XiNi72}) and (\ref{XiNi72a}), we have
\begin{equation}\nonumber
\deg(G_2) = 2g+ \sum_{i=1}^{s-1} d_{i}
					   -[\eta_2 d_{s,1}] = 2g+ (s-1)d_0e [m \epsilon]  -[\eta_2(m-t+1-(s-1)d_0e [m \epsilon])] 
\end{equation}
\begin{equation}\nonumber
 \leq 2g+ (s-1)d_0e [m \epsilon] -\eta_2(m-t+1 -(s-1)d_0e [m \epsilon]) +1
      =(1+\eta_2) (s-1)d_0e [m \epsilon] 
\end{equation}
\begin{equation}\nonumber
- m \eta_2 + 2g   +1 +\eta_2( t-1)  \leq m \eta_2 ( (1+\eta_2)/2 -1) +1 +2g+\eta_2 t <0
\end{equation}
 for $m > 2(1 +2g+ \eta_2t) \eta_2^{-1}  (1-\eta_2)^{-1}$.
Hence $ \beta_1 +     \beta_2=0$.\\

By [DiPi, Lemma 8.10] (or Lemma 5), we get that 
 $ b_{i,j} =0 $ for all $j \in [1, d_{i}]$, $i \in [1,s-1]$  and  $ b_{0,j} =0 $ for  $j \in [0, g]$.\\
From (\ref{XiNi80}) we have 
$ \ddot{\varphi}_r =0$  for  $r \in [0,m+g]$ .
Thus $\nu_{ P_{\infty}} (\beta_3) \geq m+g+1.$  \\
Applying (\ref{XiNi78}), we derive 
\begin{equation}\nonumber
  \beta_3 \in \cL(G_3)  \quad {\rm with } \quad  G_3 = G+ [(d_{s,2}-1)/e_s +1] P_s  
					  -(m +g+1 )  P_{\infty}.
\end{equation}
By (\ref{XiNi72a}), we obtain 
\begin{equation}\nonumber
\deg(G_3)   = 2g+  m-t -(s-2)d_0 e [m \epsilon] +e_s-m -g-1 \leq g-t-1+e_s -(s-2)d_0 e [m \epsilon]
<0  
\end{equation}
for $m \geq \epsilon^{-1}$ and $s \geq 3$.
Hence $\beta_3=0$.
Using   (\ref{XiNi04}) and  (\ref{XiNi78}), we get  that $ b_{s,j} =0 $ for all $j \in [d_{s,1}, d_{s,2}]$.

 By (\ref{XiNi76}), we have a contradiction.
Thus  assertions (\ref{XiNi72}) and (\ref{XiNi15a})  are  true.
Therefore Theorem 2  is proved. \qed   \\

\subsection{Niederreiter-\"{O}zbudak nets. Proof of Theorem 3}  
Let
\begin{equation}\label{NiOz10a}
          m=m_i e_i +r_i, \;\;\; {\rm with }  \;\;\; 0 \leq r_i <e_i,\; 1\leq i \leq s \; \ad  \;\tilde{r}_0=\sum_{i=1}^{s-1} r_i, \;\; r_0=\sum_{i=1}^{s} r_i.
\end{equation}

{\bf Lemma 7.} {\it There exists a divisor $\tilde{G}$  of $F/\FF_b$ with $\deg(\tilde{G}) =  g -1 +\tilde{r}_0$,
such that $\nu_{P_i}(\tilde{G}) = 0 $  for $1 \leq i \leq s$, and }
\begin{equation}\nonumber
           \cN_m(P_1,...,P_s;G)   =  \cN_m(P_1,...,P_s;\hat{G}), \quad {\rm where }  \quad
					\hat{G}=m_1P_1+...+m_{s-1}P_{s-1} +\tilde{G} .
\end{equation}

{\bf Proof.} We have $\nu_{P_i}(G) =a_i$ and $\nu_{P_i}(t_i) =1$ for $1 \leq i \leq s$.
Using the Approximation Theorem, we obtain that there exists $y \in F$, such that
\begin{equation}\label{NiOz10}
           \nu_{P_i}(y - t_i^{a_i -m_i}) = a_i  +1, \quad     {\rm for}  \quad 1 \leq i \leq s-1,    \quad  \nu_{P_s}(y - t_s^{a_s}) = a_s +m_s +1.
\end{equation}
Let  $\dot{f}=fy$ and $\hat{G} =G-\div(y)$.
We note 
\begin{equation}\label{NiOz12}
  f \in \cL(G) \Leftrightarrow \div(f) +G \geq 0     \Leftrightarrow \div(fy) +G -\div(y) \geq 0 \Leftrightarrow \dot{f}=fy \in \cL(\hat{G}).
\end{equation}
  It is easy to see that   $\nu_{P_i}(\hat{G}) = m_i $  $(1 \leq i \leq {s-1})$, 
	 $\nu_{P_s}(\hat{G}) = 0 $ and 
 $\deg(\hat{G}) = \deg(G)=m(s-1)+ g-1$. Let $\tilde{G} =\hat{G} -m_1P_1-...-m_{s-1}P_{s-1}$. We get
 $\nu_{P_i}(\tilde{G}) = 0 $  for $1 \leq i \leq s$. 
Hence
\begin{equation} \nonumber
\deg(\tilde{G}) =m(s-1)+ g-1 -e_1m_1 -... -e_{s-1} m_{s-1} = g-1 +\tilde{r}_0 .
\end{equation}
Let $\dot{f}_{i,j} =S_{j}(t_i,\dot{f}) $  (see (\ref{NiOz02})). By  (\ref{NiOz10}), we have
\begin{equation}\nonumber
\dot{f}_{i,  -j} = f_{i, -a_i +m_i-j }
  \;\; 1 \leq i \leq s-1, \;\;\; {\rm and } \;\;\;  \dot{f}_{s, m_s -j} =
    f_{s,-a_s +m_s -j  }  \;\;   \with \;\; 1 \leq j \leq m_s. 
\end{equation}
Using notations (\ref{NiOz04}), we get
\begin{equation}\nonumber
  \theta_i^{(\hat{G})}(\dot{f}) = (\bs_{r_i}, \vartheta_i(\dot{f}_{i, -1}), ...,
	\vartheta_i(\dot{f}_{i, -m_i})) = (\bs_{r_i},   \vartheta_{i}(f_{i, - a_i+m_i - 1}), . . . , \vartheta_{i}(f_{i, - a_i}))   =\theta_i^{(G)}(f) 
\end{equation}
 for $1 \leq i \leq s-1$, and
\begin{equation}\nonumber
  \theta_s^{(\hat{G})}(\dot{f}) = (\bs_{r_s}, \vartheta_s(\dot{f}_{s, m_s -1}), ...,
	\vartheta_s(\dot{f}_{s, 0})) = (\bs_{r_s},   \vartheta_{s}(f_{s, - a_s+m_s - 1}), . . . , \vartheta_{s}(f_{s, - a_s}))= 
\end{equation}	
$  \theta_s^{(G)}(f)  $. 
By (\ref{NiOz06}),  we have 
\begin{equation} \nonumber
  \theta^{(\hat{G})}(\dot{f}) := (  \theta_1^{(\hat{G})}(\dot{f}),..., \theta_s^{(\hat{G})}(\dot{f})) =
		(  \theta_1^{(G)}(f),..., \theta_s^{(G)}(f)) =\theta^{(G)}(f)  
\end{equation}
for  all $f \in \cL(G)$.   From (\ref{NiOz07}) and (\ref{NiOz12}) , we obtain the assertion of  Lemma 7. \qed \\

By Lemma 7, we can take  $\hat{G}$ instead of  $G$. Hence
\begin{equation}\label{NiOz13a}
					G=m_1P_1+...+m_{s-1}P_{s-1} +\tilde{G},  \quad {\rm and}  \quad
					 a_i=m_i, \;\; 
1 \leq i \leq s-1,\;\;\;  a_s=0.
\end{equation} 
Let $\vartheta_i = (\vartheta_{i,1} ,..., \vartheta_{i, e_i})$.
From (\ref{NiOz04}), we get   for $ 0 \leq \check{j}_i \leq m_i-1, \; 1 \leq \hat{j}_i \leq e_i,$
that
\begin{equation}\nonumber
             \theta^{(G)}_i(f)=(\theta_{i,1}(f),...,\theta_{i,m}(f))  =
						(\bs_{r_i}, \vartheta_{i}(f_{i, - 1 }), . . . 
						,  \vartheta_{i}(f_{i, - m_i})) ,  \;1 \leq i \leq s-1,
\end{equation} 
 with $\theta_{i,r_i+\check{j}_i e_i+\hat{j}_i}(f) = \vartheta_{i,\hat{j}_i}(f_{i,-\check{j}_i -1}) $, 
and
\begin{equation} \label{NiOz16}
             \theta^{(G)}_s(f) =(\theta_{s,1}(f),...,\theta_{s,m}(f)) = 
						(\bs_{r_s}, \vartheta_{s}(f_{s,m_s - 1}), . . . , \vartheta_{s}(f_{s,0})) ,  
\end{equation}
 with $\theta_{s,r_s+\check{j}_s e_s+\hat{j}_i}(f) = \vartheta_{s,\hat{j}_s}(f_{s,m_s -
\check{j}_s-1})$.\\

{\bf Lemma 8.} {\it Let $\vartheta_i =(\vartheta_{i,1} ,..., \vartheta_{i, e_i}) \; : \; F_{P_i}  \to  \FF^{e_i}_b$
be an $\FF_b$-linear vector space isomorphism. Then there exists  an $\FF_b$-linear vector space isomorphism  $ \vartheta^{\bot}_i=(\vartheta^{\bot}_{i,1} ,..., \vartheta^{\bot}_{i, e_i})    : F_{P_i}  \to  \FF^{e_i}_b$ such that}
\begin{equation}\nonumber
        \Tr_{F_{P_i}/\FF_b} (\dot{x}\ddot{x}) = \sum_{j=1}^{e_i} \vartheta_{i,j}(\dot{x})  \vartheta^{\bot}_{i,j}(\ddot{x})
   \quad { for \; all} \quad \dot{x},\ddot{x}\in  F_{P_i},\quad 1 \leq i \leq s .
\end{equation}

{\bf Proof.}
Using Theorem F, we get that there exists  $\beta_{i,j} \in  F_{P_i}$   such that
\begin{equation}\label{NiOz23ac}
       \vartheta_{i,j}(y)  =   \Tr_{F_{P_i}/\FF_b} (y \beta_{i,j})   \quad {\rm for }  \quad  1 \leq j \leq e_i,  
\end{equation}
and $(\beta_{i,1},..., \beta_{i,e_i} ) $  is the basis of $F_{P_i}$ over $\FF_b$ ($ 1 \leq i \leq s $).
Applying Theorem G, we obtain that there exists a  basis $(\beta^{\bot}_{i,1},..., \beta^{\bot}_{i,e_i}) $  of  $F_{P_i}$ over $\FF_b$  such that
\begin{equation}\nonumber
        \Tr_{F_{P_i}/\FF_b} ( \beta_{i,j_1} \beta^{\bot}_{i,j_2}  )  = \delta_{j_1,j_2} \quad {\rm with }  \quad  1 \leq j_1, j_2 \leq e_i.
\end{equation}
Let   $\dot{x} =\sum_{j=1}^{e_i} \dot{\gamma}_j \beta^{\bot}_{i,j} $, 
  $\ddot{x} =\sum_{j=1}^{e_i} \ddot{\gamma}_j \beta_{i,j} $ 
	and let
\begin{equation}\label{NiOz23ab}
      \vartheta^{\bot}_{i,j}(\ddot{x}):  =  \ddot{\gamma}_{j}= \Tr_{F_{P_i}/\FF_b} (\ddot{x}\beta^{\bot}_{i,j}) . 
\end{equation}	
By (\ref{NiOz23ac}), we have $\dot{\gamma}_{j} = \vartheta_{i,j}(\dot{x}) $.  Now, we get
\begin{equation}\nonumber
       \Tr_{F_{P_i}/\FF_b} (\dot{x}\ddot{x}) = \sum_{j_1,j_2=1}^{e_i} \dot{\gamma}_{j_1} \ddot{\gamma}_{j_2}   \Tr_{F_{P_i}/\FF_b} (\beta^{\bot}_{i,j_1}  \beta_{i,j_2} )
			=\sum_{j=1}^{e_i} \dot{\gamma}_{j} \ddot{\gamma}_{j} = \sum_{j=1}^{e_i} \vartheta_{i,j}(\dot{x})  \vartheta^{\bot}_{i,j}(\ddot{x}).
\end{equation}
Hence Lemma 8 is proved. 
\qed  \\

We consider the $H$-differential  $dt_s$. Let $\omega$ be the corresponding
Weil differential,  $\div(\omega)$ the divisor of $\omega$, and $W:=\div(dt_s) = \div(\omega)$.
By  (\ref{No12}) and (\ref{No16}), we have
\begin{equation}\label{NiOz23a}
    \deg(W)=2g-2 \quad \ad \quad \nu_{P_s}(W)=\nu_{P_s}(d t_s)=\nu_{P_s}(d t_s/d t_s)=0. 
\end{equation}
Using notations of Lemma 7, we define
\begin{equation}\label{NiOz23b}
G^{\bot} =m_sP_s - \tilde{G} +W, \;\;\; \where \;\;\; \deg(\tilde{G})=g-1+\tilde{r}_0   \quad \ad \quad \nu_{P_i} (\tilde{G}) =0
\end{equation}
for  $1 \leq i \leq s$.
Let $a_i^{\bot} :=\nu_{P_i} (G^{\bot} -W)$ for  $1 \leq i \leq s$.
We obtain from (\ref{NiOz23b}) that 
$a_i^{\bot}=0$  for  $1 \leq i \leq s-1$ and $a_s^{\bot}=m_s$.
Let $f^{\bot} \in \cL(G^{\bot})$, then $ \div(f^{\bot}) + W +G^{\bot} -W \geq~0   $
 and $\nu_{P_i}( \div(f^{\bot} ) +W )\geq  - \nu_{P_i}( G^{\bot} -W  )$.
Applying (\ref{No16}),   we get 
\begin{equation}\label{NiOz23c}
\nu_{P_i}( f^{\bot}    \dd t_s) = \nu_{P_i}( f^{\bot}) +    \nu_{P_i}(W) \geq  -\nu_{P_i} (G^{\bot} -W)=- a_i^{\bot}, \;\; \with \;\; a_i^{\bot}=0,
\end{equation}
$1 \leq i \leq s-1$, and  $a_s^{\bot}=m_s$ for $f^{\bot} \in \cL(G^{\bot})$.
According to Proposition A, we have that there exists $\tau_i \in F$, such that
\begin{equation}\label{NiOz23d}
   \dd t_s = \tau_i \dd t_i, \qquad 1 \leq i \leq s.
\end{equation}
From (\ref{No12}) and (\ref{NiOz23c}),   we get 
\begin{equation}  \nonumber
\nu_{P_i}( f^{\bot}  \tau_i) = \nu_{P_i}( f^{\bot}  \tau_i \dd t_i) = 
\nu_{P_i}( f^{\bot}  \dd t_s)
\geq - a_i^{\bot}, \qquad 1 \leq i \leq s.
\end{equation}
 By (\ref{No06}), we have the local expansions
\begin{equation}\label{NiOz18}
   f^{\bot}    \tau_i\; := \; \sum_{j = -a_i^{\bot} }^{\infty}
	S_{j}(t_i,f^{\bot} \tau_i) t_i^{j}  , \quad {\rm where \; all }  \quad   S_{j}(t_i,f^{\bot} \tau_i) \in F_{P_i}
\end{equation}
for $1 \leq i \leq s$ and $f^{\bot} \in \cL(G^{\bot})$.  We denote $S_{j}(t_i,f^{\bot} \tau_i) $ by $f^{\bot}_{i,j}$.

Using  (\ref{No18}), (\ref{No19}) and (\ref{NiOz23ab}),   we denote
\begin{equation}\label{NiOz20}
     \vartheta^{\bot}_{i,\hat{j}_i} (  f^{\bot}_{i,\check{j}_i} ) := 
		\Tr_{F_{P_i}/\FF_b} ( \beta^{\bot}_{i,\hat{j}_i } f^{\bot}_{i,\check{j}_i})
= \underset{P_i,t_i}\Res( \beta^{\bot}_{i,\hat{j}_i}   t_i^{-\check{j}_i -1}  f^{\bot} \tau_i)  
\end{equation}
 and $ \vartheta_{i}^{\bot} =(\vartheta_{i,1}^{\bot},...,\vartheta_{i,e_i}^{\bot}) $
with $1 \leq \hat{j}_i \leq e_i$, $-a_i^{\bot} \leq \check{j}_i \leq -a_i^{\bot} + m_i-1$, $1 \leq i \leq s$.\\
For $f^{\bot}  \in \cL(G^{\bot})$, the image of $f^{\bot}$ under $\dot{\theta}^{\bot}_i$, for $1 \leq i \leq s$, is
defined as
\begin{equation}\nonumber
             \dot{\theta}^{\bot}_i(f^{\bot}) =
						( \dot{\theta}^{\bot}_{i,1}(f^{\bot}) ,...,
		\dot{\theta}^{\bot}_{i,m}(f^{\bot})) := 
		( \vartheta_{i}^{\bot}(f^{\bot}_{i,-a^{\bot}_i}), . . . , \vartheta^{\bot}_{i}(f^{\bot}_{i, -a^{\bot}_i+m_i-1}) ,\bs_{r_i}) \in \FF^m_b,
\end{equation}
It is easy to verify that
\begin{equation}\label{NiOz20a}
\dot{\theta}^{\bot}_{i, \check{j}_ie_i +\hat{j}_i}(f^{\bot}) = 
	\vartheta_{i,\hat{j}_i}^{\bot}(f^{\bot}_{i,\check{j}_i}), \quad {\rm for }  \quad
	1 \leq \hat{j}_i \leq e_i,
\; 0 \leq \check{j}_i \leq m_i-1,
\end{equation}
\begin{equation}\label{NiOz21}
          1 \leq i \leq  s-1            \quad \ad \quad 
\dot{\theta}^{\bot}_{s, \check{j}_s e_s +\hat{j}_s}(f^{\bot}) = 
	\vartheta_{s,\hat{j}_s}^{\bot}(f^{\bot}_{s, - m_s+\check{j}_s}), \;  0 \leq \check{j}_s 
	\leq m_s-1.
\end{equation}
Let
\begin{equation} \label{NiOz20ab}
						\dot{\theta}^{(G,{\bot})}(f^{\bot}) \;: = \; \big(\dot{\theta}^{\bot}_1(f^{\bot}), . . . , 	
				\dot{\theta}^{\bot}_s(f^{\bot})\big) \in \FF^{ms}_b.	
\end{equation}
 Let $\bvarphi_i =(\varphi_{i,1},...,\varphi_{i,r_i})$ with $\varphi_{i,j} \in \FF_b$   $(1 \leq j \leq r_i,\;1 \leq i \leq s)$, and let
\begin{equation}\label{NiOz19}
   \Phi =\{ \bvarphi =(\bvarphi_1,...,\bvarphi_s) \; |  \; \bvarphi_i \in \FF_b^{r_i}, \;  i=1,...,s  \} \;\; \with \;\; \dim(\Phi)= r_0=\sum_{i=1}^{s} r_i.
\end{equation}
Now, we set
\begin{equation} \label{NiOz22}
						\theta^{(G,{\bot})}(f^{\bot},\bvarphi)  \; : = \; \big(\theta^{\bot}_1(f^{\bot},\bvarphi), . . . , 	
				\theta^{\bot}_s(f^{\bot},\bvarphi)\big) \in \FF^{ms}_b	,
\end{equation}
where
\begin{equation}\nonumber
             \theta^{\bot}_i(f^{\bot},\bvarphi) =( \theta^{\bot}_{i,1}(f^{\bot},\bvarphi) ,...,
		\theta^{\bot}_{i,m}(f^{\bot},\bvarphi)) := 
		(\bvarphi_{i}, \dot{\theta}^{\bot}_{i,1}(f^{\bot}), ...,\dot{\theta}^{\bot}_{i,m-r_i}(f^{\bot}) ) \in \FF^m_b.
\end{equation}
We define the $\FF_b$-linear maps
\begin{equation}\label{NiOz23}
             \theta^{(G,{\bot})} \; : \; \big(\cL(G^{\bot}), \Phi\big)  \to  \FF^{ms}_b , 
						\quad      (f^{\bot}, \bvarphi)   \mapsto  \theta^{(G,{\bot})}(f^{\bot}, \bvarphi)
\end{equation}
\begin{equation}\nonumber 
\ad \qquad
\dot{\theta}^{(G,{\bot})}  \; :  \; \cL(G^{\bot})  \to  \FF^{ms}_b , \quad  
    f^{\bot}   \mapsto \dot{\theta}^{(G,{\bot})}(f^{\bot}) 						. 
\end{equation}
The images of $\theta^{(G,{\bot})}$ and $\dot{\theta}^{(G,{\bot})}$ are denoted by
\begin{equation}\label{NiOz24}
         \Xi_m := \{    \theta^{(G,{\bot})}(f^{\bot}, \bvarphi)  \;  |  \;f^{\bot}\in 
				\cL(G^{\bot}), \;\bvarphi \in \Phi \}
\end{equation}
\begin{equation} \nonumber
        \ad \qquad
         \dot{\Xi}_m := \{    \dot{\theta}^{(G,{\bot})}(f^{\bot})  \; | \; f^{\bot}
				\in \cL(G^{\bot}) \}.				
\end{equation}\\

{\bf Lemma 9} {\it With notation as above, we have  ${\rm ker}(\theta^{(G,{\bot})}) =\bs$ and}

\begin{equation} \nonumber
        \delta^{\bot}_m(\dot{\Xi}_m)  \leq m+ g-1  +e_0  -r_0.
\end{equation}

{\bf Proof.} Consider  (\ref{NiOz23a})-(\ref{NiOz23d}). Let $f^{\bot} \in \cL(G^{\bot}) \setminus \{ 0\}$, and let
\begin{equation}\label{NiOz132}
 \nu_{P_i}(f^{\bot} \tau_i) = d_i  \quad {\rm for  }  \quad 1 \leq i \leq s-1, \qquad
 \nu_{P_s}(f^{\bot} ) = d_s-m_s.
\end{equation}
We see that
\begin{equation}\label{NiOz134}
\div(f^{\bot}) + G^{\bot} \geq 0,  \quad {\rm with  }  \quad G^{\bot}= m_sP_s -\tilde{G} +W
             \quad {\rm and  }  \quad W = (\dd t_s).
\end{equation}
Hence
\begin{equation}\label{NiOz136}
\nu_{P}\big( \div(f^{\bot}) +m_s P_s -\tilde{G} +W\big) \geq 0,
  \quad {\rm for \; all }  \quad P \in \PP_F.
\end{equation}
By (\ref{No12}) and (\ref{No16}),   we obtain $ \nu_{P_i}(W)= \nu_{P_i}(d t_s) =
 \nu_{P_i}( \tau_i)$, $1 \leq i \leq s$.\\
Bearing in mind (\ref{NiOz132}) and that $\nu_{P_i}(\tilde{G}) =0$ for $i \in[1, s]$,  
 we get
\begin{equation}   \nonumber
    \nu_{P_i}( \div(f^{\bot}) +m_s P_s -\tilde{G} +W) =d_i\geq 0,  \quad  \quad 1 \leq i \leq s .
\end{equation}
Therefore
\begin{equation} \nonumber
\nu_{P_i}(\div(f^{\bot}) +\dot{G}) \geq 0
 \;\; \for  \;\; f^{\bot} \in \cL(G^{\bot}) \setminus \{ 0\},  
 \; {\rm where  } \; \dot{G} = G^{\bot} -\sum_{i=1}^s d_iP_i
\end{equation}
and $G^{\bot}= m_sP_s -\tilde{G} +W$.
Taking into account that $f^{\bot} \in \cL(G^{\bot}) \setminus \{ 0\} $, we obtain
\begin{equation} \nonumber
   0 \leq \deg(\dot{G}) = \deg\Big(G^{\bot} -\sum_{i=1}^s d_iP_i \Big) =\deg(G^{\bot}) -\sum_{i=1}^s d_ie_i.
\end{equation}
By (\ref{NiOz23a}), (\ref{NiOz23b}) and (\ref{NiOz10a}), we get
\begin{equation}  \nonumber
          \sum_{i=1}^s d_ie_i  \leq \deg(m_sP_s - \tilde{G} +W) =m_s e_s - (g-1 + \tilde{r}_0)
					+2g-2 =
					m - r_0+g-1. 
\end{equation}
According to (\ref{NiOz18}), (\ref{NiOz20}) and (\ref{NiOz132}),  we obtain
\begin{equation} \nonumber
       f^{\bot}_{i,a^{\bot}_i+j}= 0 \quad {\rm for  }  \quad 0 \leq j <d_i 
			\quad \ad \quad f^{\bot}_{i,a^{\bot}_i+d_i} 
			\neq 0, \quad  1 \leq i \leq s.
\end{equation}
From    (\ref{NiOz127}), (\ref{NiOz21}) and Lemma 8,  we have
\begin{equation} \nonumber
        v_m^{\bot} (\dot{\theta}^{\bot}_i(f^{\bot})) \leq  (d_i+1)e_i   \quad {\rm for  }  \quad 1 \leq i \leq s.
\end{equation}
Applying (\ref{NiOz20ab}) and (\ref{NiOz128}), we derive
\begin{equation}  \nonumber
 V^{\bot}_m\big(\dot{\theta}^{(G,{\bot})}(f^{\bot})\big)  \leq \sum_{i=1}^s 
(d_i+1)e_i \leq m +g-1 +e_0  -r_0.
\end{equation}
By (\ref{NiOz129}), $\delta^{\bot}_m(\dot{\Xi}_m)  \leq m +g-1 +e_0-r_0$. 
 Taking into account (\ref{NiOz127}) and that $s \geq 3$, we get
${\rm ker}(\theta^{(G,{\bot})}) =\bs$.

Therefore  Lemma 9 is proved. \qed  \\

{\bf Lemma 10.} {\it With notation as above, we have that $\dim(\Xi_m) =m. $} \\

{\bf Proof.}   
By (\ref{NiOz23a}) and (\ref{NiOz23b}), we have 
\begin{equation}\nonumber
 \deg(G^{\bot})=\deg(m_sP_s -\tilde{G} +W) =m_se_s -\deg(\tilde{G}) +2g-2=
 m -r_s +2g-2 -\tilde{r}_0  - g+1.
\end{equation}
Using (\ref{NiOz10a}) and the Riemann-Roch theorem, we obtain for $m \geq g+e_0-1 \geq g+r_0$
that
\begin{equation}\nonumber 
  \dim(\cL(G^{\bot})) = \deg(m_sP_s -\tilde{G} +W) -g+1 =m -r_0+2g-2  -2g+2 =m-r_0.
\end{equation}
From (\ref{NiOz19}), we have $\dim(\Phi) =r_0$.
Hence   
\begin{equation}\nonumber
 \dim\big((\cL(G^{\bot}),  \Phi  )\big) =
 \dim(\cL(G^{\bot})) +\dim(\Phi) =m-r_0 +r_0 =m .
\end{equation}
By Lemma 9, we get ${\rm ker}(\theta^{(G,{\bot})}) =\bs $.
Bearing in mind that $\theta^{(G,{\bot})}\big((\cL(G^{\bot}),  \Phi  )\big) =\Xi_m$, we obtain the assertion of  Lemma 10. \qed  \\

{\bf Lemma 11.} {\it Let $f \in \cL(G)$, and $f^{\bot} \in \cL(G^{\bot})$. Then
}
\begin{equation}\label{NiOz150}
  \sum_{i=1}^s \underset{P_i}\Res(f f^{\bot} \dd t_s)=0,
\end{equation}
\begin{equation}\label{NiOz152}
 \underset{P_i}\Res(f f^{\bot} \dd t_s) =\sum_{j=0 }^{m_i -1}      \Tr_{F_{P_i}/\FF_b}\big(    f_{i,-j-1} \;
 f^{\bot}_{i,j}   \big), \qquad  1 \leq i \leq s-1
\end{equation}  
\begin{equation}\label{NiOz154}
\quad \ad \quad
  \underset{P_s}\Res(f f^{\bot} \dd t_s) = \sum_{j=0}^{m_s-1}   \Tr_{F_{P_s}/\FF_b}
 \big( f_{s,m_s-j-1} \;
 f^{\bot}_{s,-m_s+j}   \big).
\end{equation}

{\bf Proof.} By (\ref{NiOz13a}) and (\ref{NiOz23b}), we have 
\begin{equation}  \nonumber
					G=m_1P_1+...+m_{s-1}P_{s-1} +\tilde{G},  \quad {\rm and}  \quad
		G^{\bot} =m_sP_s - \tilde{G} +W		.
\end{equation} 
Bearing in mind that $ \div(f) +  G \geq 0 $, $ \div(f^{\bot}) +  G^{\bot} \geq 0 $ and that $W= \div(\dd t_s)$, we obtain
\begin{equation}\nonumber
 \div(f) +  \sum_{i=1}^{s} m_i P_i + \tilde{G}  + \div(f^{\bot})  -\tilde{G} +W =
 \div(f) +  \div(f^{\bot}) +\sum_{i=1}^s m_i P_i +\div(\dd t_s) \geq 0.
\end{equation}
From (\ref{No16}), we derive
\begin{equation}\nonumber
 \nu_P(f f^{\bot} \dd t_s) = \nu_P(f f^{\bot})+\nu_P(\div( \dd t_s))  \geq 0 \quad {\rm and }  \quad  \underset{P}\Res(f f^{\bot} \dd t_s)=0 
\end{equation}
for all $P \in \PP_f \setminus \{ P_1,...,P_s\}$.\\
Applying the Residue Theorem, we get    assertion (\ref{NiOz150}).\\
By (\ref{NiOz02})  and (\ref{NiOz18}), we derive
\begin{equation}\nonumber
  \underset{P_s}\Res(f f^{\bot} \dd t_s) =   \underset{P_s}\Res\Big(  \sum_{j_1=0}^{\infty}  S_{j_1}(t_s,f)  t_s^{j_1}  \;  \sum_{j_2=-m_s}^{\infty}    
 S_{j_2}(t_s, f^{\bot} )  t_s^{j_2} \dd t_s\Big) 
\end{equation}
\begin{equation}\nonumber
=  \sum_{j_1=0}^{\infty}   \sum_{j_2=-m_s}^{\infty}     \underset{P_s}\Res\big(    S_{j_1}(t_s,f) \;
 S_{j_2}(t_s, f^{\bot} )  t_s^{j_1+j_2} \dd t_s \big) 
\end{equation}
\begin{equation}\nonumber
 =\sum_{0 \leq j_1 \leq m_s -1, \;    j_1+j_2=-1}      \Tr_{F_{P_s}/\FF_b}\big(    S_{j_1}(t_s,f) \;
 S_{j_2}(t_s, f^{\bot} )  \big) 
\end{equation}
\begin{equation}\nonumber
 =\sum_{j=0 }^{m_s -1}      \Tr_{F_{P_s}/\FF_b}\big(    S_{m_s-j-1}(t_s,f) \;
 S_{-m_s+j}(t_s, f^{\bot} )  \big) = \sum_{j=0 }^{m_s -1}      \Tr_{F_{P_s}/\FF_b}\big(f_{s,m_s-j-1} \;
 f^{\bot}_{s,-m_s+j}   \big).
\end{equation}
Hence assertion (\ref{NiOz154}) is proved. \\
Analogously, using (\ref{NiOz23d}),  we have
\begin{equation}\nonumber
\underset{P_i}\Res(f f^{\bot} \dd t_s) =  \underset{P_i}\Res(f f^{\bot} \tau_i\dd t_i) =
\underset{P_i}\Res\Big(  \sum_{j_1=-m_i}^{\infty}  S_{j_1}(t_i,f)  t_i^{j_1}   \; \sum_{j_2=0}^{\infty}    
 S_{j_2}(t_i, f^{\bot} \tau_i)  t_i^{j_2} \dd t_i\Big) 
\end{equation}
\begin{equation}\nonumber
=\sum_{0 \leq j_2 \leq m_i -1, \;    j_1+j_2=-1}     \Tr_{F_{P_i}/\FF_b}\big(    S_{j_1}(t_i,f) \;
 S_{j_2}(t_i, f^{\bot} \tau_i)  \big),  
\end{equation}
\begin{equation}\nonumber
          =\sum_{j=0 }^{m_i -1}      \Tr_{F_{P_i}/\FF_b}\big(    f_{i,-j-1} \;
 f^{\bot}_{i,j}   \big), \quad {\rm for }  \quad  1 \leq i \leq s-1.
\end{equation}
Thus  Lemma 11 is proved. \qed    \\

{\bf Lemma 12.} {\it With notation as above, we have $\Xi_m =\cN^{\bot}(P_1,...,P_s,G)$.
}\\

{\bf Proof.} Using (\ref{NiOz08}) and Lemma 10, we have 
\begin{equation}  \nonumber
\dim_{\FF_b}(\cN_m) = ms-m  \quad \ad \quad \dim_{\FF_b}(\Xi_m) =m.
\end{equation}
From (\ref{NiOz07}), (\ref{NiOz23}) and (\ref{NiOz24}), we get that $\cN_m, \Xi_m \subset \FF_b^{ms}$.

By (\ref{Ap304}), in order to obtain the assertion of the lemma, it is sufficient to prove that $A \cdot B =0$ for all $A \in \cN_m$ and $B	\in \Xi_m$.

According to (\ref{NiOz04}), (\ref{NiOz07}), (\ref{NiOz16})  and   (\ref{NiOz21}) -  (\ref{NiOz24}), it is enough to verify that
\begin{equation}\label{NiOz160}
 A \cdot B 
 = \sum_{i=1}^s \eth_i =0 \quad {\rm with  }  \quad \eth_i =\sum_{j=1}^m
  \theta_{i,j}(f)	  \theta^{\bot}_{i,j}((f^{\bot},\bvarphi))
		\quad {\rm for \; all  }  \quad f \in \cL(G),
\end{equation}
and $(f^{\bot},\bvarphi)		\in (\cL(G^{\bot}),\Phi)$.
From (\ref{NiOz16}) and (\ref{NiOz20}) -  (\ref{NiOz21}), we derive
\begin{equation}\label{NiOz162}
\eth_i =\sum_{\check{j}_i=0}^{m_i -1}  \varkappa_{i,j_1} \quad {\rm with  }  \quad   \varkappa_{i,\check{j}_i}=\sum_{\hat{j}_i=1}^{e_i} 
  \theta_{i,r_i + \check{j}_i e_i+\hat{j}_i}(f)	 \; \; \theta^{\bot}_{i,r_i + \check{j}_i e_i+\hat{j}_i}((f^{\bot},\bvarphi)).
\end{equation}
Using (\ref{NiOz16}) and (\ref{NiOz21})-(\ref{NiOz22}), we have
  for $ \check{j}_i \in [0,m_i -1], \;  \hat{j}_i \in [1, e_i]$
\begin{equation}\nonumber
 \theta_{s,r_s+\check{j}_s e_s+\hat{j}_s}(f) = \vartheta_{s,\hat{j}_s}(f_{s,m_s -\check{j}_s -1})
	\quad {\rm and  }  \quad
       \theta^{\bot}_{s,r_s+\check{j}_s e_s+\hat{j}_s}((f^{\bot},\bvarphi)) 	=
		\vartheta^{\bot}_{s,\hat{j}_s}(f^{\bot}_{s,-m_s +\check{j}_s}) ,
\end{equation}
\begin{equation}\nonumber
        \theta_{i,r_i+\check{j}_i e_i+\hat{j}_i}(f) = \vartheta_{i,\hat{j}_i}(f_{i,-\check{j}_i -1}) 
	 \;\; {\rm and}  \;\;
\theta^{\bot}_{i,r_1+\check{j}_ie_i+\hat{j}_i}((f^{\bot},\bvarphi)) 	=
		\vartheta^{\bot}_{i,\hat{j}_i}(f^{\bot}_{i, \check{j}_i}), \; 1 \leq i \leq s-1.
\end{equation} 
By Lemma 8 and  (\ref{NiOz162}), we obtain
\begin{equation}\nonumber
  \varkappa_{s,\check{j}_s} = \sum_{\hat{j}_i=s}^{e_s} 
  \vartheta_{s,\hat{j}_s} (f_{s, m_s- \check{j}_s - 1})	 \;  \vartheta^{\bot}_{s,\hat{j}_s} (f^{\bot}_{s, -m_s+ \check{j}_s})=  
     \Tr_{F_{P_s}/\FF_b}
 \big( f_{s,m_s-\check{j}_s-1} \;
 f^{\bot}_{s,-m_s+\check{j}_s}   \big)
\end{equation}
and
\begin{equation}\nonumber
  \varkappa_{i,\check{j}_i} = \sum_{\hat{j}_i=1}^{e_i} 
  \vartheta_{i,\hat{j}_i} (f_{i,- \check{j}_i - 1})	 \;  \vartheta^{\bot}_{i,\hat{j}_i} (f^{\bot}_{i,  \check{j}_i})=  
     \Tr_{F_{P_i}/\FF_b}
 \big( f_{i,-\check{j}_i-1} 
 f^{\bot}_{i,\check{j}_i}   \big) \quad {\rm for }  \quad  1 \leq i \leq s-1.
\end{equation}
From (\ref{NiOz152}), (\ref{NiOz154}) and (\ref{NiOz162}),  we get
\begin{equation}\nonumber
  \eth_i =\underset{P_i}\Res(f f^{\bot} \dd t_s) \quad {\rm for }  \quad  1 \leq i \leq s.
\end{equation}
Applying Lemma 11, we get  assertion (\ref{NiOz160}). Hence Lemma 12 is proved. \qed \\
\\
%
%
%
%
Let  
\begin{equation}\label{NiOz200}
 G_i = \tilde{G} + q_iP_i -q_sP_s \; \with \; 
q_s =[\frac{g+\tilde{r}_0 }{e_s}] +1 \;\;\; \ad \;\;\; q_i =[\frac{g-\tilde{r}_0 +q_s e_s}{e_i}] +1
\end{equation}
for $i \in [1,  s-1]$. 
By (\ref{NiOz23b}), we have $\deg(\tilde{G}) =g-1 + \tilde{r}_0$
    and $\nu_{P_i}(\tilde{G}) =0$, $i \in [1,  s]$.
It is easy to see that $\deg(G_i)  \geq 2g-1$, $i \in [1,  s-1]$. Let $z_i=\dim(\cL(G_i))$,
 and let ${u^{(i)}_{1},...,u^{(i)}_{z_i}}$ be a basis of $\cL(G_i)$  over $\FF_b$, $i \in [1,  s-1]$.

For each $i \in [1,  s-1]$, we consider the chain
\begin{equation}  \nonumber
          \cL(G_i) \subset \cL(G_i + P_i) \subset \cL(G_i + 2P_i)\subset...
\end{equation}
of vector spaces over $\FF_b$. By starting from the basis ${u^{(i)}_{1},...,u^{(i)}_{z_i}}$ 
of $\cL(G_i)$ 
and successively adding basis vectors at each step of the chain, we obtain
for each $n \geq q_i$ a basis
\begin{equation}\label{NiOz204}
         \{u^{(i)}_{1},...,u^{(i)}_{z_i}, k^{(i)}_{q_i,1},...,k^{(i)}_{q_i,e_i},...,
	 k^{(i)}_{n,1},...,k^{(i)}_{n,e_i}			\}
\end{equation}
of $\cL(G_i + (n-q_i+1) P_i)$. We note that we then have
\begin{equation}\label{NiOz206}
    k^{(i)}_{j_1,j_2} \in \cL(G_i + (j_1-q_i+1)P_i) \; \ad \;
		\nu_{P_i} (k^{(i)}_{j_1,j_2}) =-j_1-1, \;\nu_{P_s} (k^{(i)}_{j_1,j_2}) \geq q_s \
\end{equation} 
for $j_1 \geq  q_i$,
$1 \leq j_2 \leq e_i, \;1 \leq i \leq s-1$.\\
Let $\check{G} = \tilde{G}+gP_s$. We see that $\deg(\check{G}) =g-1 + \tilde{r}_0 +g e_s \geq 2g-1$. 
Let  ${u^{(0)}_{1},...,u^{(0)}_{z_0}}$ be a basis of $\cL(\check{G})$  over $\FF_b$.
In a similar way, we construct a basis \\
$   \{u^{(0)}_{1},...,u^{(0)}_{z_0}, k^{(i)}_{0,1},...,k^{(i)}_{0,e_i},...,
	 k^{(i)}_{(q_i-1),1},...,k^{(i)}_{(q_i-1),e_i} \}$ of $\cL(\check{G} + q_iP_i)$
 with
\begin{equation}\label{NiOz207}
    k^{(i)}_{j_1,j_2} \in \cL(\check{G} + (j_1+1)P_i) \; \ad \; \nu_{P_i} (k^{(i)}_{j_1,j_2}) =-j_1 -1\;\;{\rm for} \;\;  j_1 \in [0,q_i),
\end{equation} 
 $1 \leq j_2 \leq e_i, \;1 \leq i \leq s-1$.\\

Now, consider the chain
\begin{equation}\nonumber
    \cL(q_s P_s -\tilde{G}+W) \subset  \cL((q_s+1) P_s -\tilde{G}+W) \subset ... \subset
			\cL(G^{\bot} - P_s) \subset \cL(G^{\bot}), 
\end{equation}
where $ G^{\bot} = m_sP_s -\tilde{G}+W$ and $q_s =[(g+\tilde{r}_0 )/e_s] +1$. By (\ref{NiOz23a}) and (\ref{NiOz23b}), we have $\deg(\tilde{G}) =g-1 + \tilde{r}_0$,
   $\deg(W) =2g-2$ and $\nu_{P_s}(\tilde{G}) =\nu_{P_s}(W)=0$. Hence $\deg(q_s P_s -\tilde{G}+W) \geq 2g-1$.
Let  ${u^{(s)}_{1},...,u^{(s)}_{z_s}}$ be a basis of $\cL(q_s P_s -\tilde{G}+W)$  over $\FF_b$.
In a similar way, we construct a basis 
$\{u^{(s)}_{1},...,u^{(s)}_{z_s}, k^{(s)}_{q_s,1},...,k^{(s)}_{q_s,e_s},...,
	 k^{(s)}_{n,1},...,k^{(i)}_{n,e_s} \}$ of $\cL((n+1) P_s -\check{G} + W)$
 with
\begin{equation}\label{NiOz207a}
    k^{(s)}_{j_1,j_2} \in \cL((j_1+1)P_s-\check{G} + W) \;\; \ad \;\; \nu_{P_s} (k^{(s)}_{j_1,j_2}) =-j_1-1 \;\; {\rm for} \;\; j_1 \geq q_s
\end{equation} 
and $j_2 \in [1, e_s]$. 
By (\ref{NiOz204})-(\ref{NiOz207}), we have the following local expansions
\begin{equation}\label{NiOz218}
  k^{(i)}_{j_1,j_2} \; := \; \sum_{r = -j_1 }^{\infty} \varkappa^{(i,j_2)}_{j_1,r}
	t_i^{r-1}   \quad \for  \quad
	\varkappa^{(i,j_2)}_{j_1,r} \in F_{P_i},
	\;\; \; \; i \in [1,s].
\end{equation}
\\

{\bf Lemma 13}. {\it Let $j_i \geq 0$ for $i \in [1,s-1]$ and let  $j_s \geq q_s$. Then
$\{\varkappa^{(i,1)}_{j_i,-j_i},...,\varkappa^{(i,e_i)}_{j_i,-j_i}\}$ is a basis of $F_{P_i}$
over $\FF_b$ for $i \in [1,s]$.}\\

{\bf Proof.}   Let $i \in [1,s-1]$ and let $j_i \geq q_i$.  
Suppose that there exist $a_1,...,a_{e_i} \in \FF_b$, such that 
$\sum_{1 \leq j \leq e_i} a_i \varkappa^{(i,j)}_{j_i,-j_i} =0$ and 
$(a_1,...,a_{e_i}) \neq (\bar{0},...,\bar{0})$.
By (\ref{NiOz218}), we get $\nu_{P_i} (\alpha) \geq -j_i$, where $\alpha:=
\sum_{1 \leq j_2 \leq e_i} a_i k^{(i)}_{j_i,j_2} $. Hence $\alpha \in \cL(G_i + (j_i -q_i)P_i)$. We have a contradiction with the construction of the basis vectors (\ref{NiOz204}).

Similarly, we can consider the cases $i \in [1,s-1]$, $j_i \in [0, q_i-1]$ and $i=s$.
Therefore Lemma 13 is proved. \qed \\

{\bf Lemma 14}. {\it Let $d_i \geq 1$ be an integer $(i=1,...,s-1)$ and   $f^{\bot} \in G^{\bot}$. 
Suppose that 
$\Res_{P_s,t_s}(f^{\bot}  k^{(i)}_{j_1,j_2} ) =0 $ for
  $j_1 \in [0, d_i-1], j_2 \in [1, e_i]$  and    $i  \in [1, s-1]$. 
Then}  
\begin{equation}\label{NiOz208}
  \vartheta^{\bot}_{i,j_2} (  f^{\bot}_{i,j_1} ) =0 \quad
\for  \quad  j_1 \in [0, d_i-1], \; j_2 \in [1, e_i] \;    \ad  \;  i  \in [1, s-1].
\end{equation}\\

{\bf Proof.}  
By (\ref{NiOz134}), (\ref{NiOz136}),   (\ref{NiOz200}), (\ref{NiOz206}) and (\ref{NiOz207}), we have 
$\nu_{P}\big( \div(f^{\bot}) +m_s P_s -\tilde{G} +W\big) \geq 0$, for  all $P \in \PP_F$
 and 
$k^{(i)}_{j_1,j_2}\in \cL(\tilde{G} +a_{j_1}P_s + (j_1+1)P_i) $ with some integer $a_{j_1}$. \\
From (\ref{No12}), (\ref{No16}) and  (\ref{No18}), we derive
\begin{equation}\nonumber
 \nu_P(f^{\bot}  k^{(i)}_{j_1,j_2} \dd t_s)  \geq 0 \quad {\rm and }  
\quad  \underset{P}\Res(f^{\bot}  k^{(i)}_{j_1,j_2} \dd t_s)=0 \quad {\rm for\; all }  \quad 
   P \in \PP_F \setminus \{ P_i,P_s\}.
\end{equation}
Applying (\ref{NiOz23d}) and the Residue Theorem, we get 
\begin{equation} \nonumber
   \underset{P_i, t_i}\Res(f^{\bot} \tau_i k^{(i)}_{j_1,j_2} ) =\underset{P_i}\Res(f^{\bot}  k^{(i)}_{j_1,j_2} \dd t_s)=- 
\underset{P_s}\Res(f^{\bot}  k^{(i)}_{j_1,j_2} \dd t_s) =- 
\underset{P_s,t_s}\Res(f^{\bot}  k^{(i)}_{j_1,j_2}) 
\end{equation}
for  all $0 \leq j_1, \;1 \leq j_2 \leq e_i,\;  1 \leq i  \leq  s-1 $.

By  (\ref{NiOz18}), (\ref{NiOz218}) and the conditions of the lemma, we obtain
\begin{equation} \nonumber
  -\underset{P_s,t_s}\Res(f^{\bot}  k^{(i)}_{j_1,j_2} )=  \underset{P_i,t_i}\Res(f^{\bot} \tau_i k^{(i)}_{j_1,j_2} ) = 
	 \underset{P_i,t_i}\Res\Big(\sum_{j = 0}^{\infty} f^{\bot}_{i,j}
	t_i^{j}  \sum_{r = -j_1 }^{\infty} \varkappa^{(i,j_2)}_{j_1,r}
	t_i^{r-1}  \Big)
\end{equation}
\begin{equation}\label{NiOz220}
 =\sum_{j = 0}^{\infty}  \sum_{r = -j_1}^{\infty} 
 \Tr_{F_{P_i}/\FF_b} ( f^{\bot}_{i,j} \varkappa^{(i,j_2)}_{j_1,r}  ) \delta_{j,-r}  =  \sum_{j = 0}^{j_1}
\Tr_{F_{P_i}/\FF_b} ( f^{\bot}_{i,j} \varkappa^{(i,j_2)}_{j_1,-j}  )=0 
\end{equation}
 for $0 \leq j_1 \leq d_i-1, \; 1 \leq j_2 \leq e_i, \;    \ad \;\; 1 \leq i  \leq  s-1$.\\
Consider  (\ref{NiOz220}) for $j_1=0$. We have 
$\Tr_{F_{P_i}/\FF_b} ( f^{\bot}_{i,0} \varkappa^{(i,j_2)}_{0,0}  )=0$ for all $j_2 \in [1.e_i]$.
By Lemma 13, we obtain that $f^{\bot}_{i,0}=0$. Suppose that $f^{\bot}_{i,j}=0$ for 
$0 \leq j <j_0$. Consider (\ref{NiOz220}) for $j_1=j_0$.
We get 
$\Tr_{F_{P_i}/\FF_b} ( f^{\bot}_{i,j_0} \varkappa^{(i,j_2)}_{j_0,-j_0}  )=0$ for all $j_2 \in [1.e_i]$. Applying Lemma 13, we have that $f^{\bot}_{i,j_0}=0$.
By induction, we obtain that $f^{\bot}_{i,j}=0$ for all $j \in [0,d_i-1]$ and $i \in [1,s-1]$.
Now, using (\ref{NiOz20}), we get that assertion (\ref{NiOz208}) is true. 
Hence Lemma 14 is proved.~\qed\\

{\bf Lemma 15.} {\it Let  $s \geq 3$ , $\{\beta^{\bot}_{s, 1}, ..., \beta^{\bot}_{s, e_{s}}\}$ be a basis of $F_{P_s}/\FF_b$,
\begin{equation}  \nonumber
   		\Lambda_1 = \Big\{ \big( \underset{P_s,t_s}\Res(f^{\bot}  k^{(i)}_{j_1,j_2} ) \big)_{d_{i,1} \leq j_1 \leq  d_{i,2}, 1 \leq j_2 \leq e_i, 1 \leq i \leq s-1},
			\qquad\qquad\qquad\qquad\qquad\quad
\end{equation}
\begin{equation}						
\quad\qquad\qquad\qquad	\big( \underset{P_s,t_s}\Res(\beta^{\bot}_{s,j_2}f^{\bot}
 t_s^{m_s -j_1 -1} ) \big)_{d_{s,1} \leq j_1 \leq  d_{s,2}, 1 \leq j_2 \leq e_s}		
			\; | \; f^{\bot} \in \cL(G^{\bot}) \Big\} \nonumber
\end{equation}
with $ d_{s,1} =m_s +1-[t/e_s]-(s-1)d_0 \dot{m} e/e_s $, 	$\dot{m} =[\tilde{m} \epsilon]$,
 $\tilde{m} =m-r_0$,		
\begin{equation}\label{NiOz253}
  	d_{s,2} =m_s -2-[t/e_s]- (s-2)d_0 \dot{m} e/e_s, \;\;\; d_{i,1}=q_i, \;
		d_{i,2} = d_0\dot{m}]e/e_i-1,
\end{equation}
 $i \in [ 1,s-1]$, $d_0 =d+t,$  $e=e_1e_2\cdots e_s$, $\epsilon = \eta(2(s-1)d_0 e)^{-1}$,
$\eta =(1+\deg((t_s)_\infty))^{-1} $.
 Then}
\begin{equation}\label{NiOz254}
   		\Lambda_1 = \FF_b^{\chi}, \; \with \; \chi =\sum_{i=1}^s (d_{i,2} - d_{i,1} +1)e_i
			\; for \; 
  m > 2(g-1 +e_0)e_s +2t(\eta^{-1} -1).
\end{equation}
\\

{\bf Proof.}
Suppose that (\ref{NiOz254}) is not true. 
Then there exists $b^{(i)}_{j_1,j_2} \in \FF_b$ $(i,j_1,j_2 \geq 1)$
 such that
\begin{equation}\label{NiOz258}
\sum_{i=1}^{s} \sum_{j_1=d_{i,1}}^{d_{i,2}} \sum_{j_2=1}^{e_{i}}|b^{(i)}_{j_1,j_2}|  >0 
\end{equation}
and
\begin{equation}\label{NiOz260}
   \sum_{i=1}^{s-1} \sum_{j_1=d_{i,1}}^{d_{i,2}} \sum_{j_2=1}^{e_{i}} b^{(i)}_{j_1,j_2}
		\underset{P_s,t_s}\Res(f^{\bot}  k^{(i)}_{j_1,j_2} ) +
\sum_{j_1=d_{s,1}}^{d_{s,2}} \sum_{j_2=1}^{e_{s}} b^{(s)}_{j_1,j_2} 
\underset{P_s,t_s}\Res(\beta^{\bot}_{s,j_2}f^{\bot} t_s^{m_s -j_1 -1} ) =0 
\end{equation}
for  all  $f^{\bot} \in \cL(G^{\bot})$.
Let $\alpha = \alpha_1 + \alpha_2$ with
\begin{equation}\label{NiOz262}
  \alpha_1=  \sum_{i=1}^{s-1} 
	\sum_{j_1=d_{i,1}}^{d_{i,2}} \sum_{j_2=1}^{e_{i}} b^{(i)}_{j_1,j_2}  k^{(i)}_{j_1,j_2}  \quad \ad \quad \alpha_2=
\sum_{j_1=d_{s,1}}^{d_{s,2}} \sum_{j_2=1}^{e_{s}} b^{(s)}_{j_1,j_2} \beta^{\bot}_{s,j_2}  t_s^{m_s -j_1 -1} .
\end{equation}
By  (\ref{NiOz260}), we have 
\begin{equation}\label{NiOz264}
   \underset{P_s,t_s}\Res(f^{\bot} \alpha )  =0 \quad {\rm for \;
	all}  \quad f^{\bot} \in \cL(G^{\bot}).
\end{equation}
From  (\ref{NiOz206}), we get $\nu_{P_s}(\alpha) \geq q_s$.
Consider the local expansion 
\begin{equation}  \nonumber
           \alpha = \sum_{r= q_s }^{\infty} \varphi_{r} t_s^r 
             \quad {\rm with   } \quad  \varphi_{r} \in F_{P_s} 
						\quad \for \quad r \geq q_s. 
\end{equation}
Suppose that $m_s>j_0:=\nu_{P_s}(\alpha)$. Therefore $\varphi_{j_0} \neq 0$.
From (\ref{NiOz207a}), we obtain that $ k^{(s)}_{j_0,j_2} \in \cL(G^{\bot}) $
 for all $j_2 \in [1,e_s]$.
Applying (\ref{NiOz218}) and (\ref{NiOz264}), we derive
\begin{equation}\nonumber
   \underset{P_s,t_s}\Res ( k^{(s)}_{j_0,j_2} \alpha) = \underset{P_s,t_s}\Res \Big( 
  \sum_{j = -j_0 }^{\infty} \varkappa^{(s,j_2)}_{j_0,j}
	t_s^{j-1} 
\sum_{r= j_0}^{\infty} \varphi_{r} t_s^r \Big)=\Tr_{F_{P_s}/\FF_b} (\varkappa^{(s,j_2)}_{j_0,-j_0}\varphi_{j_0}) =0 
\end{equation}
 for all $j_2 \in [1,e_s]$.
By Lemma 13, $\{\varkappa^{(s,1)}_{j_0,-j_0},...,\varkappa^{(s,e_s)}_{j_0,-j_0}\}$ is a basis of $F_{P_s}$. Hence $\varphi_{j_0} =0 $.
We have a contradiction. Thus  $\nu_{P_s}(\alpha) \geq m_s$.

We consider the compositum field   $F'=FF_{P_s}$.
Let $\fB_1,...,\fB_{\mu}$ be all the places of $F'/F_{P_s}$ lying over $P_s$. 
From   (\ref{No58}), we get
\begin{equation}\label{NiOz267}
   \nu_{\fB_i}(\alpha) \geq m_s \quad \for \quad  i=1,...,\mu.
\end{equation}
According to   (\ref{NiOz200}) and (\ref{NiOz206}), we obtain 
\begin{equation}\nonumber
       \alpha_1 \in \cL_F(A_1) =\cL(A_1), \quad \with \quad A_1:=\tilde{G}-q_s P_s+\sum_{i=1}^{s-1}(d_{i,2}+1)P_i.
\end{equation}
Applying  Theorem D(d), we have
\begin{equation}  \nonumber
       \alpha_1 \in \cL_{F'}(\Con_{F'/F} (A_1)).
\end{equation}
By (\ref{NiOz262}), we derive
\begin{equation}  \nonumber
       \alpha_2 \in \cL_{F'}(A_2), \quad \with \quad 
	A_2=((t_s)^{F'}_{\infty})^{m_s-d_{s,1}-1}	.
\end{equation}
Using (\ref{NiOz267}), we get
\begin{equation} \nonumber
       \alpha \in \cL_{F'}(A_1+A_2 - m_s \sum_{i=1}^{\mu} \fB_i).
\end{equation}
From  (\ref{No50}), Theorem D(a) and Theorem E, we derive $\Con_{F'/F} (P_s)=\sum_{i=1}^{\mu} \fB_i$, $\Con_{F'/F} ((t_s)^{F}_{\infty})=(t_s)^{F'}_{\infty}$ and 
\begin{equation}\nonumber
       \alpha  \in \cL_{F'}(A_3), \quad \with \quad 
			A_3 = \Con_{F'/F} \big(A_1+
	(m_s-d_{s,1}-1)	 (t_s)^{F}_{\infty}  -
		m_sP_s \big) .
\end{equation}
Applying  Theorem D(c) and (\ref{NiOz200}), we have
\begin{align*}
  & \deg(A_3) = \deg\Big(\tilde{G}+\sum_{i=1}^{s-1}(d_{i,2}+1)P_i + 	(m_s-d_{s,1}-1) (t_s)^{F}_{\infty} -m_s	P_s \Big)
 \\
	&	\leq  g-1+\tilde{r}_0+(s-1)d_0 e \dot{m}+  (m_s-d_{s,1}-1)\deg((t_s)_{\infty}) -m_s e_s \\
	&  \leq  g-1+e_0-e_s +(s-1)d_0e \dot{m} 
	+([t/e_s] + (s-1)d_0 \dot{m} e/e_s -2)(\eta^{-1} -1)  \\
	&-m_s e_s \leq g-1+e_0+ (t/e_s-2)(\eta^{-1}-1)   + (s-1)d_0 e \dot{m} (1+(\eta^{-1}-1)/e_s) -m \\
&\leq  g-1+e_0+ t(\eta^{-1}-1)/e_s   -m\big((2e_s)^{-1} + (1-\eta/2)(1-1/e_s)\big)\leq \beta -m/(2e_s) <0 
\end{align*}
for $m>2e_s \beta, 	$
with $\beta = g-1+e_0+ t(\eta^{-1}-1)/e_s$ and $\epsilon=\eta (2(s-1)d_0e)^{-1}$.\\
Hence $\alpha=0$. \\
\\
%
Suppose that
 $\sum_{i=1}^{s-1}\sum_{j_1=d_{i,1}}^{d_{i,2}}\sum_{j_2=1}^{e_{i}}|b^{(i)}_{j_1,j_2}|=0$.
 Then $\alpha_2 =0$ and $\sum_{j_2=1}^{e_{s}} b^{(s)}_{j_1,j_2} \beta^{\bot}_{s,j_2} =~0$ for all 
$ j_1 \in [ d_{s,1},d_{s,2}]$. Bearing in mind that  $(\beta^{\bot}_{s,j_2})_{1 \leq j_2 \leq e_2}$ is a basis of
 $F_{P_s}/\FF_b$, we get $\sum_{j_1=d_{s,1}}^{d_{s,2}} \sum_{j_2=1}^{e_{s}} |b^{(s)}_{j_1,j_2}| =0$. By (\ref{NiOz258}), we have a contradiction.

Therefore there exists $h \in [1,s-1]$ with 
 \begin{equation}\label{NiOz270a}
 \sum_{j_1=d_{h,1}}^{d_{h,2}}\sum_{j_2=1}^{e_{h}}|b^{(h)}_{j_1,j_2}|  > 0.
\end{equation}
Let $\fB_{h,1},...,\fB_{h,\mu_h}$ be all the places of $F'/F_{P_s}$ lying over $P_h$.
Let  
\begin{equation}\nonumber
\alpha_{1,i}=  
	\sum_{j_1=d_{i,1}}^{d_{i,2}} \sum_{j_2=1}^{e_{i}} b^{(i)}_{j_1,j_2}  k^{(i)}_{j_1,j_2}, 
	  \quad i=1,...,s-1.
\end{equation}
Let $\nu_{P_h}(t_s) \geq 0$ or $\alpha_2 =0$. Therefore $\nu_{\fB_{h,j}}(\alpha_2) \geq 0$
 for $1 \leq j \leq \mu_h$. Taking into account that $\alpha_1 = -\alpha_2$, we get
 $\nu_{\fB_{h,j}}(\alpha_1) \geq 0$
 for $1 \leq j \leq \mu_h$, and  $\nu_{P_h}(\alpha_1) \geq 0$.\\ 
 Using (\ref{NiOz23b}), (\ref{NiOz200}), (\ref{NiOz206}) and (\ref{NiOz253}), we obtain   $\nu_{P_h}(\alpha_{1,i}) \geq 0$  for $1 \leq i \leq s-1, i \neq h$. Bearing in mind (\ref{NiOz270a}) and that
$ \{u^{(h)}_{1},...,u^{(h)}_{z_h}, k^{(h)}_{q_h,1},...,k^{(h)}_{q_h,e_h},...,
	 k^{(h)}_{n,1},...,k^{(h)}_{n,e_h}			\}$ is a basis of  $\cL(G_h + (n-q_h+1) P_h)$, we get 
\begin{equation}\nonumber
 \alpha_{1,h} \in \cL(G_h + (j-q_h+1) P_h) \setminus
	\cL(G_h + (j-q_h) P_h)\;\;{\rm with\; some}\; j \geq q_h.
\end{equation}	
 By  (\ref{NiOz200}) and (\ref{NiOz206}), we get  $\nu_{P_h}(\alpha_{1,h}) \leq -1$. 
We have a contradiction.\\

Now let $\nu_{P_h} (t_s) \leq -1$ and  $\alpha_2 \neq 0$. 
We have $\nu_{P_h}(\alpha_{1,h}) \geq -d_{h,2}-1$,  $\nu_{P_h}(\alpha_{1}) \geq -d_{h,2}-1$
 and $\nu_{\fB_{h,j}} (\alpha_1) \geq - d_{h,2}-1$,
 $j=1,...,\mu_h$. 
On the other hand, using   (\ref{NiOz262}) and (\ref{No58}), we have  $\nu_{\fB_{h,j}} (\alpha_2) \leq - (m_s-d_{s,2}-1)$,
 $j=1,...,\mu_h$.
According to (\ref{NiOz10b}) and (\ref{NiOz253}), we obtain $s\geq 3$, $ e_h \geq e_s $ and
\begin{equation}\nonumber
m_s-d_{s,2}-1 -d_{h,2}-1 =  [t/e_s]+1 +(s-2)d_0e \dot{m} e/e_s -d_0 \dot{m}e/e_h \geq 1.
\end{equation}
We have  a contradiction.  
Thus  assertion (\ref{NiOz260})  is not true.
 Hence (\ref{NiOz254}) is true and Lemma 15 follows. \qed \\
\\
%
{\bf End of the proof of Theorem 3.}

Using (\ref{Ap302}), (\ref{NiOz08b0}), (\ref{NiOz22})-(\ref{NiOz24}) and Lemma 12, we have
\begin{equation} \label{NiOz300}
  \cP_1 =\{\tilde{\bx}(f^{\bot}, \bvarphi)=(\tilde{x}_1(f^{\bot}, \bvarphi),...,\tilde{x}_s(f^{\bot}, \bvarphi))\; \; | \;\; f^{\bot} \in \cL(G^{\bot}), \bvarphi \in \Phi \} 
\end{equation}
with
\begin{equation} \nonumber
 \tilde{x}_i(f^{\bot}, \bvarphi) = \sum_{j=1}^m \phi^{-1} (\theta^{\bot}_{i,j}(f^{\bot},\bvarphi))b^{-j} = 
\sum_{j=1}^{r_i} \phi^{-1} (\varphi_{i,j}) b^{-j} +
b^{-r_i}\sum_{j=1}^{m-r_i} \phi^{-1} (\dot{\theta}^{\bot}_{i,j}(f^{\bot}))b^{-j}.
\end{equation}
By (\ref{NiOz08b1}), we have
\begin{equation} \label{NiOz302b}
  \cP_2 =\{\dot{\bx}(f^{\bot})=(\dot{x}_1(f^{\bot}),...,\dot{x}_s(f^{\bot})) \; | \; f^{\bot} \in \cL(G^{\bot}) \} 
\end{equation}
with
\begin{equation} \label{NiOz302c}
    \dot{x}_i(f^{\bot}) = 
\sum_{j=1}^{m-r_i} \phi^{-1} (\dot{\theta}^{\bot}_{i,j}(f^{\bot}))b^{-j}, \quad 1 \leq i \leq s.
\end{equation}\\

{\bf Lemma 16.} {\it  With notation as above,
 $\cP_2 $ is a $d-$admissible $(t,m-r_0,s)$ net in base $b$ with $d=g+e_0$, and
 $t=g+e_0-s$. }\\

{\bf Proof.} Let $J =\prod_{i=1}^{s} [ A_i/b^{d_i}, (A_i +1)/b^{d_i})$ with $d_i \geq 0$, 
and $0 \leq A_i <b^{d_i}$, $1 \leq i \leq s$, and let 
$J_{\bpsi} =\prod_{i=1}^{s} [ \psi_i/b^{r_i}+A_i/b^{r_i+d_i}, 
\psi_i/b^{r_i}+ (A_i +1)/b^{r_i+d_i})$ with $\psi_i/b^{r_i} =\psi_{i,1}/b +...
+\psi_{i,r_i}/b^{r_i}$, $\psi_{i,j} \in Z_b$, $1 \leq i \leq s$, $d_1+...+d_s =m-r_0 -t$ .\\
It is easy to see, that
\begin{equation} \nonumber
     \dot{\bx}(f^{\bot}) \in J \Longleftrightarrow  \tilde{\bx}(f^{\bot}, \bvarphi) \in J_{\bpsi}
		 \quad \with \quad \psi_{i,j} = \phi^{-1} (\varphi_{i,j}), \; 1\leq j \leq r_i, \; 1 \leq i \leq s.
\end{equation}
Bearing in mind that $\cP_1$ is a $(t,m,s)$ net with $t=g+e_0-s$, we have
\begin{equation} \nonumber
    \sum_{f^{\bot} \in \cL(G^{\bot})} \d1(J,\dot{\bx}(f^{\bot})) =
	 \sum_{f^{\bot} \in \cL(G^{\bot}), \bvarphi \in \Phi } \d1(J_{\bpsi}, \bx(f^{\bot}, \bvarphi))	=b^t.
\end{equation}
Therefore  $\cP_2 $ is a $(t,m-r_0,s)$ net in base $b$ with $t=g+e_0-s$.

Using (\ref{NiOz24}), Definition 5 and Definition 9, we can get $d$ from the following equation  $- \delta_m^{\bot}(\dot{\Xi}_m) = -(m-r_0)-d+1$. Applying Lemma 9, we obtain $-(m+g-1 +e_0 -r_0) \leq -(m-r_0)-d+1$.
Hence $d \leq g+e_0$.
Thus Lemma 16 is proved. \qed  \\ 

Let $V_i \subseteq \FF_b^{\mu_i}$ be a vector space over $\FF_b$, $\mu_i \geq 1 $, $i=1,2$.
Consider a linear map $h:\; V_1 \to V_2$. By the first isomorphism theorem, we have
\begin{equation} \label{NiOz250}
    \dim_{\FF_b}(V_1) = \dim_{\FF_b} ({\rm ker}(h)) + \dim_{\FF_b}({\rm im}(h)).
\end{equation} \\
%
Let
\begin{equation}\nonumber	
   		\Lambda_1^{'} = \Big\{ \big( \underset{P_s,t_s}\Res(f^{\bot}  k^{(i)}_{j_1,j_2} ) 
		\big)_{0 \leq j_1 \leq  d_{i,2}, 1 \leq j_2 \leq e_i, 1 \leq i \leq s-1},
			\qquad\qquad\qquad\qquad\qquad\qquad\qquad\qquad\quad\\	
\end{equation}
\begin{equation}						
\quad\qquad\qquad\qquad	\big( \underset{P_s,t_s}\Res(\beta^{\bot}_{s,j_2}f^{\bot}
 t_s^{m_s -j_1 -1} ) \big)_{d_{s,1} \leq j_1 \leq  d_{s,2}, 1 \leq j_2 \leq e_s}		
			\; | \; f^{\bot} \in \cL(G^{\bot}) \Big\} \nonumber
\end{equation}
and
\begin{equation}\nonumber	
  		\Lambda_2 = \Big\{ 			
	\big( \underset{P_s,t_s}\Res(\beta^{\bot}_{s,j_2}f^{\bot} t_s^{m_s -j_1 -1} ) \big)_{d_{s,1} \leq j_1 \leq  d_{s,2}, 1 \leq j_2 \leq e_s}		
 		\; |  \;  \underset{P_s,t_s}\Res(f^{\bot}  k^{(i)}_{j_1,j_2} )=0 	\qquad\qquad\qquad\qquad\qquad\qquad\qquad\qquad
\end{equation}
\begin{equation} \nonumber		
		\qquad\qquad\qquad\qquad \for \quad  0 \leq j_1 \leq  d_{i,2}, 1 \leq j_2 \leq e_i, 1 \leq i \leq s-1, \;			
			f^{\bot} \in \cL(G^{\bot}) \Big\} 
\end{equation}
with $ d_{s,1} =m_s +1-[t/e_s]-(s-1)d_0 \dot{m} e/e_s $, 			
\begin{equation}\label{NiOz253a}
  	d_{s,2} =m_s -2-[t/e_s]- (s-2)d_0 \dot{m} e/e_s, \;\;\; d_{i,1}=q_i, \; d_{i,2} = d_0
		\dot{m} e/e_i-1,
\end{equation}
 $i \in [ 1,s-1]$, $d_0 =d+t,$  $e=e_1e_2\cdots e_s$, $\epsilon = \eta(2(s-1)d_0 e)^{-1}$, 
$\eta =(1+\deg((t_s)_\infty))^{-1} $, $\dot{m}=[\tilde{m} \epsilon]$, $\tilde{m}=m-r_0$, 
  $ m > 2(g-1 +e_0)e_s +2t(\eta^{-1} -1)$,  $d=g+e_0$ and $t=g+e_0-s$. 
	
By (\ref{NiOz250}), (\ref{NiOz253a}) and Lemma 15, we have  $\dim_{\FF_b}(\Lambda_1^{'}) \geq \dim_{\FF_b}(\Lambda_1)$ and 
\begin{equation}\nonumber
			\dim_{\FF_b}(\Lambda_2)= \dim_{\FF_b}(\Lambda_1^{'}) -
			\dim_{\FF_b}\Big( \Big\{ \big( \underset{P_{s},t_{s}}\Res(f^{\bot} k^{(i)}_{j_1,j_2} ) \big)_{\substack{0 \leq j_1 \leq  d_{i,2}, 1 \leq j_2 \leq e_i\\ 1 \leq i \leq s-1}}   |  f^{\bot} \in \cL(G^{\bot}           \Big\} \Big)
\end{equation}
\begin{equation}  \nonumber
		\geq \dim_{\FF_b}(\Lambda_1) -\sum_{i=1}^{s-1} (d_{i,2}  +1)e_i  \geq (d_{s,2} -d_{s,1}+1)e_s - \sum_{i=1}^{s-1}q_ie_i=  d_0e\dot{m} -2e_s  -\sum_{i=1}^{s-1}q_ie_i.
\end{equation}
Let
\begin{equation}\nonumber	
 		\Lambda_3 = \Big\{ 			
	\big( \underset{P_s,t_s}\Res(\beta^{\bot}_{s,j_2}f^{\bot} t_s^{m_s -j_1 -1} ) \big)_{d_{s,1} \leq j_1 \leq  d_{s,2}, 1 \leq j_2 \leq e_s}		
 		\; | \; \vartheta^{\bot}_{i,j_2} (  f^{\bot}_{i,j_1} ) =0 \qquad\qquad\qquad \qquad\qquad
		\qquad\qquad\qquad\qquad
\end{equation}		
\begin{equation} \nonumber		
\qquad\qquad	\qquad\qquad \for \; 0 \leq j_1 \leq  d_{i,2}, 1 \leq j_2 \leq e_i, 1 \leq i \leq s-1 \; | \;			
	f^{\bot} \in \cL(G^{\bot}) \Big\}. 
\end{equation}
 Using Lemma 14, we get $\Lambda_3 \supseteq \Lambda_2$ and 
$\dim_{\FF_b}(\Lambda_3) \geq \dim_{\FF_b}(\Lambda_2)$.
Let
\begin{equation}\nonumber
  		\Lambda_4 = \Big\{ 			
		\big( 	\vartheta^{\bot}_{i,j_2} (  f^{\bot}_{i,j_1} )  \big)_{ 0 \leq j_1 \leq  d_{i,2}, 1 \leq j_2 \leq e_i, 1 \leq i \leq s-1}	 \; | \; 		
			f^{\bot} \in \cL(G^{\bot}) \Big\}. 	
\end{equation}		
Taking into account that 
 $\cP_2 $ is a  $(t,m-r_0,s)$ net in base $b$, we  get from (\ref{NiOz302b})
 that $\dim_{\FF_b}(\Lambda_4) = (s-1)d_0e \dot{m}$.
Let
\begin{equation}\nonumber
   		\Lambda_5 = \Big\{ \Big( \vartheta^{\bot}_{i,j_2} (  f^{\bot}_{i,j_1} ) 
  \Big)_{\substack{0 \leq j_1 \leq  d_{i,2}, 1 \leq j_2 \leq e_i\\ 1 \leq i \leq s-1}},	
\Big( \underset{P_s,t_s}\Res(\beta^{\bot}_{s,j_2}f^{\bot} t_s^{m_s -j_1 -1} )  \Big)_{\substack{d_{s,1}\leq j_1 \leq  d_{s,2}\\ 1 \leq j_2 \leq e_s}}
 \Big|  f^{\bot} \in \cL(G^{\bot}) \Big\}. 
\end{equation}
By (\ref{NiOz200})and  (\ref{NiOz250}), we have  
\begin{equation}\nonumber
			\dim_{\FF_b}(\Lambda_5) = \dim_{\FF_b}(\Lambda_3) +\dim_{\FF_b}(\Lambda_4)  \geq 
			s d_0e \dot{m} -2e_s 
			-2(s-1)(g+e_0).
\end{equation}
Let   $\dot{m}_{1} = d_0 e \dot{m}$,  $\dot{m} = [\tilde{m} \epsilon]$,
  $\ddot{m}_i=0$,	$ i \in [1,  s-1] $  and 
$\ddot{m}_{s} =m-t -(s-1)\dot{m}_1  $.
Bearing in mind that 
$ \dot{\theta}^{\bot}_{i, \check{j}_ie_i +\hat{j}_i}(f^{\bot}) = 
	\vartheta_{i,\hat{j}_i}^{\bot}(f^{\bot}_{i,\check{j}_i})$  for 
	$	1 \leq \hat{j}_i \leq e_i,
\; 0 \leq \check{j}_i \leq m_i-1$, 	
	$i \in [1,s-1]$ (see  (\ref{NiOz20a})), we obtain  
\begin{equation}\label{NiOz253a2}
   					\Big( 	\dot{\theta}^{\bot}_{i,\ddot{m}_i+j}(f^{\bot})  \Big)_{
			1 \leq j \leq  \dot{m}_1, 1 \leq i \leq s-1}  \supseteq 
			\Big( \vartheta^{\bot}_{i,j_2} (  f^{\bot}_{i,j_1} ) 
  \Big)_{0 \leq j_1 \leq  d_{i,2}, 1 \leq j_2 \leq e_i,  1 \leq i \leq s-1}.
\end{equation}
From (\ref{NiOz253a}), we have  $\ddot{m}_s <  d_{s,1} e_s$ and
$  (d_{s,2}+1) e_s < \ddot{m}_s +\dot{m}_1$.
Taking into account that
\begin{equation} \nonumber
  \dot{\theta}^{\bot}_{s, j_1 e_s +j_2}(f^{\bot}_{s, -m_s +j_1}) 
 =  \vartheta^{\bot}_{s, j_2}(f^{\bot}) 
			=\Res_{P_s,t_s}(\beta^{\bot}_{s,j_2}f^{\bot} t_s^{m_s -j_1 -1} )
\end{equation}
(see  (\ref{NiOz20}) and (\ref{NiOz21})),  we get
\begin{equation} \label{NiOz253a1}
  \Big( 	\dot{\theta}^{\bot}_{s,\ddot{m}_s+j}(f^{\bot})  \Big)_{
	1 \leq j \leq  \dot{m}_1}
 \supseteq
    \Big( \underset{P_s,t_s}\Res(\beta^{\bot}_{s,j_2}f^{\bot} t_s^{m_s -j_1 -1} )  
		  \Big)_{d_{s,1} \leq j_1 \leq  d_{s,2}, 1 \leq j_2 \leq e_s} .  
\end{equation}  
Let
\begin{equation}\nonumber
   		\Lambda_6 = \Big\{ \Big(
			\Big( 	\dot{\theta}^{\bot}_{i,\ddot{m}_i+j}(f^{\bot})  \Big)_{
			1 \leq j \leq  \dot{m}_1, 1 \leq i \leq s}\Big)
\; \Big| \; f^{\bot} \in \cL(G^{\bot}) \Big\}.
\end{equation}
By   (\ref{NiOz253a2}) and (\ref{NiOz253a1}), we derive
\begin{equation}\nonumber
		\dim_{\FF_b}(\Lambda_6) \geq \dim_{\FF_b}(\Lambda_5) \geq
		sd_0e \dot{m} -2e_s -2(s-1)(g+e_0) .
\end{equation}
Applying  (\ref{Ap302}), (\ref{NiOz08b1}), (\ref{NiOz302b}) and Lemma 2, we get that there exists  \\$B_i \in \{0,...,\dot{m}-1\}$, $1 \leq i \leq s$ 
such that 
\begin{equation}\label{NiOz257a}
    \Lambda_7 =\FF_b^{sd_0e\dot{m}-d_0e B} \qquad \qquad \for \quad \dot{m} \geq 1 ,
\end{equation}
where  $B =\#B_1+...+\#B_{s} \leq 4(s-1)(g+e_0)$ and  
\begin{equation} \nonumber
  		  \Lambda_7  = \Big\{   \Big( 	
\dot{\theta}^{\bot}_{i,\ddot{m}_i+\dot{j}_i d_0 e +\ddot{j}_i}(f^{\bot})  \; | \; 
	 \dot{j}_i  \in \bar{B}_{i}, \; \ddot{j}_i \in [1, d_0 e],  \;     i \in [1,s] \Big)  \; | 
 f^{\bot} \in \cL(G^{\bot})          \Big\}
\end{equation}
with $\bar{B}_i =\{0,...,\dot{m}-1\} \setminus B_{i}$.\\
From (\ref{NiOz302c}), we have
\begin{equation} \nonumber
  		 \Big\{   \Big( 	
\dot{x}_{i,\ddot{m}_i+\dot{j}_i d_0 e +\ddot{j}_i}(f^{\bot})   |  
	 \dot{j}_i  \in \bar{B}_{i},  \ddot{j}_i \in [1, d_0 e],       i \in [1,s] \Big)  \; | 
 f^{\bot} \in \cL(G^{\bot})          \Big\} =Z_b^{sd_0e\dot{m}-d_0e B}.
\end{equation}
We apply Corollary 2   with $\dot{s} =s$, $\tilde{r} =r_0$, $\tilde{m} =m-r_0$,
$\epsilon = \eta (2(s-1)d_0 e)^{-1}$ and $\hat{e}=e=e_1e_2\cdots e_s$.

Let $\dot{\bgamma}(f^{\bot},\dot{\bw}) =\dot{\bgamma} = (\dot{\gamma}^{(1)},...,
 \dot{\gamma}^{(\dot{s})})$ with 
$\dot{\gamma}^{(i)}:= [(\dot{\bx}(f^{\bot})  \oplus \dot{\bw})^{(i)}]_{\dot{m}_i}$, $i \in [1,s]$.
Using (\ref{NiOz302c}) and (\ref{NiOz257a}), we get that there exists $f^{\bot} \in G^{\bot}$ such that $ \dot{\bgamma}(f^{\bot},\dot{\bw})$
satisfy (\ref{Ap400}). Bearing in mind Lemma 16, we get from  Corollary 2 that
\begin{equation}  \label{NiOz270b}
\Big| \Delta((\dot{\bx}(f^{\bot}) \oplus \dot{\bw})_{f^{\bot} \in G^{\bot}}, J_{\dot{\bgamma}})
  \Big| \geq
			 		2^{-2}	b^{-d} K_{d,t,s}^{-s+1} \eta^{s-1}m^{s-1} .
\end{equation}
for  $ m \geq 2^{2s+3} b^{d+t+s}(d+t)^{s}(s-1)^{2s-1} (g+e_0) e \eta^{-s+1}$.

Taking into account (\ref{In02}), and that $\dot{\bw} \in E_{m-r_0}^{s}$ is arbitrary, we get the  second assertion in Theorem 3.\\

Consider the first assertion in Theorem 3.\\
Let $\tilde{\bgamma} = (\tilde{\gamma}^{(1)},...,\tilde{\gamma}^{(s)})$ with 
$\tilde{\gamma}^{(i)} = b^{-r_i}\dot{\gamma}^{(i)}$, $i \in [1,s]$, and let 
 $\tilde{\bw} = (\tilde{w}^{(1)},...,\tilde{w}^{(s)}) \in E_m^s  $ with 
$ \tilde{w}^{(i)}_{j+r_i} = \dot{w}^{(i)}_{j}$ for $j \in [1,m-r_0]$, $i \in [1,s]$.
By (\ref{NiOz300}) and (\ref{NiOz302b}), we have  
\begin{equation} \nonumber 
 \tilde{x}_i(f^{\bot}, \bvarphi) \oplus \tilde{w}^{(i)}\in
		[0,\tilde{\gamma}_i)
      \Longleftrightarrow \dot{x}_i(f^{\bot})  \oplus \dot{w}^{(i)} \in
			[0,\dot{\gamma}_i) 
		 \;\; \ad \;\;  \phi^{-1} (\varphi_{i,j})\oplus \tilde{w}_{i,j}=0 
\end{equation}
for $j \in [1, r_i],  \;  i \in[1, s]$.
Hence
\begin{equation} \nonumber
\sum_{ \bvarphi \in \Phi } (\d1([\bs,\tilde{\bgamma}), \tilde{\bx}(f^{\bot}, \bvarphi) \oplus 
\tilde{\bw})	-
	\tilde{\gamma}_0) =
   \d1([\bs,\dot{\bgamma}),\dot{\bx}(f^{\bot}) \oplus \dot{\bw} ) -\dot{\gamma}_0
	 , 
\end{equation}
where  $[\bs,\dot{\bgamma}) =\prod_{i=1}^s [0,\dot{\gamma}^{(i)})$, $ [\bs, \tilde{ \bgamma}) =\prod_{i=1}^s [0,\tilde{\gamma}^{(i)})$,  $ \tilde{\gamma}_0 =\tilde{\gamma}^{(1)}... \tilde{\gamma}^{(s)}$ and \\
 $ \dot{\gamma}_0 =\dot{\gamma}^{(1)}... \dot{\gamma}^{(s)}$.
Therefore
\begin{equation} \nonumber
 \sum_{f^{\bot} \in \cL(G^{\bot}),\bvarphi \in \Phi} \big(\d1([\bs,\tilde{\bgamma}), \tilde{\bx}(f^{\bot}, \bvarphi)  \oplus \tilde{\bw})	-
	 \tilde{\gamma}_0 \big)
 =  \sum_{f^{\bot} \in \cL(G^{\bot})} \big(  \d1([\bs,\dot{\bgamma}),\dot{\bx}(f^{\bot}) \oplus \dot{\bw}) - 
  \dot{\gamma}_0   \big)  .
\end{equation}
Using (\ref{In00}), (\ref{In02}) and (\ref{NiOz270b}), we get the  first assertion in Theorem 3.\\
Thus Theorem 3  is proved. \qed   \\

\subsection{Halton-type sequences. Proof of Theorem 4.}
%
%
%
%
Using  (\ref{Hal14}) and (\ref{Hal16}), we  define the sequence $(\bx_{n,j}^{(i)})_{j \geq 1}$ by
\begin{equation}\label{Hal25}
\sum_{j_2=1}^{e_i} x^{(i)}_{n,j_1 e_i +j_2} b^{-j_2+e_i} :=\sigma_{P_i}(f^{(i)}_{n,j_1}), \quad
x^{(i)}_n :=  \sum_{j=0}^{\infty} \frac{x^{(i)}_{n,j}}{ b^j }=\sum_{j_1=0}^{\infty} \sum_{j_2=1}^{e_i} \frac{x^{(i)}_{n,j_1 e_i +j_2}}{ b^{j_1e_i+j_2} },     
\end{equation}
$ 1 \leq i \leq s$, with $ (x^{(1)}_n,...,x^{(s)}_n )= \bx_n  = \xi(f_n)$,  and  $n=0,1,...$ .\\

{\bf Lemma 17}. {\it  $ (\bx_n)_{ n \geq 0} $ is  $d-$admissible with $d=g+e_0$,
 where $e_0=e_1+...+e_s.$}  \\

{\bf Proof.} 
  Suppose that the assertion of the lemma is not true. 
By  (\ref{In04}), there exists $\dot{n}> \dot{k}$ such that	
	$\left\|\dot{n}\ominus \dot{k}
	\right\|_b \left\| \bx_{\dot{n}} \ominus \bx_{\dot{k}}  \right\|_b  < b^{-d}$.\\
Let	 $d_i+1 = \dot{d}_i e_i +\ddot{d}_i $ with $1 \leq \ddot{d}_i \leq e_i$, 	$1 \leq i \leq s$,
$n= \dot{n}\ominus \dot{k} $, $\left\|n	\right\|_b =b^{m-1}$ and let
 $\left\| \bx^{(i)}_{\dot{n}} \ominus \bx^{(i)}_{\dot{k}}  \right\|_b =b^{-d_i-1} $,
 $1 \leq i \leq s$. 
Hence $m-1 -  \sum_{i=1}^s (d_i+1)  \leq -d-1$, and
\begin{equation}\label{Hal26}
         m+g-1 -  \sum_{i=1}^s \dot{d}_i e_i   \leq m+g-1 -  \sum_{i=1}^s (d_i+1)  +e_0   \leq -d-1  +g +e_0 <0.          
\end{equation}
We have
\begin{equation}\label{Hal27}
   a_{m-1}(n) \neq 0, \; a_{r}(n) = 0, \; \for \; r \geq m, \quad 
	x^{(i)}_{\dot{n},d_i+1} \neq x^{(i)}_{\dot{k},d_i+1}, \;
	x^{(i)}_{\dot{n},r} = x^{(i)}_{\dot{k},r}               
\end{equation}
for $ r \leq d_i, \;  1 \leq i \leq s$.
From  (\ref{Hal25}), we get
\begin{equation} \nonumber
     f^{(i)}_{\dot{n}, j_1} = f^{(i)}_{\dot{k}, j_1} \quad \ad \quad 
		 f^{(i)}_{n, j_1} =   0 \quad \for \quad 0 \leq j_1 < \dot{d}_i,\; 1 \leq i \leq s.
\end{equation} 
Suppose that $f^{(i)}_{n, \dot{d}_i} =   0$, then $f^{(i)}_{\dot{n}, \dot{d}_i} = f^{(i)}_{\dot{k}, \dot{d}_i}$
 and $x^{(i)}_{\dot{n},j} = x^{(i)}_{\dot{k},j}   $ for $1 \leq j \leq (\dot{d}_i+1)e_i$.
Taking into account that $d_i+1 \leq (\dot{d}_i+1)e_i$, we have a contradiction.
Therefore $f^{(i)}_{n, d_i} \neq   0$, for all $1 \leq i \leq s$.
Applying (\ref{Hal12}), we derive 
$\nu_{P_i}(f_{n}) =\dot{d}_i$, $1 \leq i \leq s$. \\
Using (\ref{Hal03})-(\ref{Hal07}) and (\ref{Hal27}), we obtain $f_n \in \cL((m+g-1)P_{s+1} -\sum_{i=1}^s
\dot{d}_iP_i) \setminus \{ 0\}$. \\
By  (\ref{Hal26}), we get
\begin{equation}\nonumber
 \deg\big((m+g-1)P_{s+1} -\sum_{i=1}^s
\dot{d}_i P_i  \big) = m+g-1 -\sum_{i=1}^s  \dot{d}_i e_i <0.
\end{equation} 	
Hence $f_n =0$. 
We have a contradiction. 
Thus Lemma 17  is proved. \qed  \\


Consider the $H-$differential $\dd t_{s+1}$.  By Proposition A,
we have that there exists $\tau_i$ with $\dd t_{s+1} = \tau_i \dd t_i$, $1 \leq i \leq s$.
Let $W = \div(\dd t_{s+1})$,  and let
\begin{equation}\label{Hal40}
 G_i =  W + q_i P_i-gP_{s+1}, \qquad \with \quad q_i =[(g+1)/e_i+1], 
	 \;\; 1 \leq i \leq s.
\end{equation}
It is easy to see that $\deg(G_i) \geq 2g-2+ g+1 -g = 2g-1$, $1 \leq i \leq s$. Let $z_i=\dim(\cL(G_i))$,
 and let ${u^{(i)}_{1},...,u^{(i)}_{z_i}}$ be a basis of $\cL(G_i)$  over $\FF_b$, $1 \leq i \leq s$.\\
For each $1 \leq i  \leq s-1$, we consider the chain
\begin{equation}  \nonumber
          \cL(G_i) \subset \cL(G_i + P_i) \subset \cL(G_i + 2P_i)\subset...
\end{equation}
of vector spaces over $\FF_b$. By starting from the basis ${u^{(i)}_{1},...,u^{(i)}_{z_i}}$ 
of $\cL(G_i)$ 
and successively adding basis vectors at each step of the chain, we obtain
for each $n \geq q_i$ a basis
\begin{equation}  \nonumber
         \{u^{(i)}_{1},...,u^{(i)}_{z_i}, k^{(i)}_{q_i,1},...,k^{(i)}_{q_i,e_i},...,
	 k^{(i)}_{n,1},...,k^{(i)}_{n,e_i}			\}
\end{equation}
of $\cL(G_i + (n-q_i+1) P_i)$. 
We note that we then have
\begin{equation}\label{Hal46}
    k^{(i)}_{j_1,j_2} \in \cL(G_i + (j_1-q_i+1)P_i) \setminus \cL(G_i + (j_1-q_i)P_i) 
\end{equation}  
for $q_i \leq j_1,
		 \; 1 \leq j_2 \leq e_i$ and $1 \leq i  \leq  s.$ 
Hence
\begin{equation}\nonumber
 \div(k^{(i)}_{j_1,j_2}) +W -gP_{s+1} + (j_1+1)P_i \geq 0 \;\;   and \;\; \nu_{P_{s+1}} (k^{(i)}_{j_1,j_2}) + \nu_{P_{s+1}} (W) \geq g .
\end{equation}
From  (\ref{No12}) and (\ref{No16}), we obtain
\begin{equation}\nonumber
    \nu_{P_{s+1}} (k^{(i)}_{j_1,j_2}) = \nu_{P_{s+1}} (k^{(i)}_{j_1,j_2}\dd t_{s+1}) = 
				\nu_{P_{s+1}} (k^{(i)}_{j_1,j_2}) + \nu_{P_{s+1}} (W) .
\end{equation}
Therefore 
\begin{equation}\label{Hal47}
    \nu_{P_{s+1}} (W) =0  \quad \ad \quad   \nu_{P_{s+1}} (k^{(i)}_{j_1,j_2})  \geq g.
\end{equation}

 Now, let $\check{G}_i = W +(e_i+1)P_{s+1} -P_i$. We see that $\deg(\check{G}_i) = 2g-1$. 
Let  ${\dot{u}^{(i)}_{1},...,\dot{u}^{(i)}_{\dot{z}_i}}$ be a basis of $\cL(\check{G}_i)$  over $\FF_b$.
In a similar way, we construct a basis $   \{\dot{u}^{(i)}_{1},...,\dot{u}^{(i)}_{\dot{z}_i}, k^{(i)}_{0,1},...,k^{(i)}_{0,e_i},...,
	 k^{(i)}_{q_i-1,1},...,k^{(i)}_{q_i-1,e_i} \}$ of $\cL(\check{G} + q_iP_i)$
 with
\begin{equation}\label{Hal48}
    k^{(i)}_{j_1,j_2} \in \cL(\check{G} + (j_1+1)P_i)\setminus\cL(\check{G} + j_1P_i)
\; {\rm for} \; j_1 \in [0,q_i),
		 j_2 \in [1,e_i],   i  \in [1, s] .
\end{equation} \\


{\bf Lemma 18}. {\it  Let $\{\beta^{(i)}_{1}, ..., \beta^{(i)}_{e_i}\}$ be a basis of $F_{P_i}/\FF_b$,
 $s \geq 2$,	 $d_i \geq 1$ be integer $(i=1,...,s)$ and   $n \in [0,b^m)$. 
Suppose that 
$\Res_{P_{s+1},t_{s+1}}(f_n  k^{(i)}_{j_1,j_2} ) =0 $ for\\
  $j_1 \in [0, d_i-1],\; j_2 \in [1, e_i]$ and $i  \in [1, s]$. 
Then}  
\begin{equation} \nonumber
\Tr_{F_{P_i}/\FF_b}(\beta^{(i)}_{j_2}f_{n, j_1}^{(i)} ) =0 \quad
for  \quad  j_1 \in [0, d_i-1], \; j_2 \in [1, e_i] \;    \ad  \;  i  \in [1, s].
\end{equation}\\

{\bf Proof.}  
Using  (\ref{Hal46}) and (\ref{Hal48}), we get
\begin{equation}\nonumber
 \nu_{P_i} (k^{(i)}_{j_1,j_2}) =-j_1-1 - \nu_{P_i} (W) \quad \for \quad j_1 \geq 0,  \; j_2 \in [1, e_i] \;    and  \;  i  \in [1, s]. 
\end{equation} 
From  (\ref{No12}) and (\ref{No16}), we obtain
\begin{equation}\label{Hal49a}
        \nu_{P_i} (\tau_i) = \nu_{P_i} (\tau_i \dd t_i)=\nu_{P_i} ( \dd t_{s+1})
				=\nu_{P_i} (\div( \dd t_{s+1}))= \nu_{P_i} ( W).
\end{equation}  
Hence 
\begin{equation}\label{Hal49b}
    \nu_{P_i} (k^{(i)}_{j_1,j_2}\tau_i) = -j_1-1 \quad \for \quad j_1 \geq 0, 
		\; j_2 \in [1, e_i] \;    and  \;  i  \in [1, s].
\end{equation}
By (\ref{Hal46}) and (\ref{Hal48}), we have
\begin{equation}\label{Hal49c}
  \div(  k^{(i)}_{j_1,j_2})  +\div(\dd t_{s+1})
	+(j_1 +1)P_i +a_{j_1}P_{s+1}  \geq 0 
\end{equation}
for $ j_1 \geq 0, \; j_2 \in [1, e_i]$,   $i  \in [1, s]$ and some $ a_{j_1} \in \ZZ$.
According to  (\ref{Hal03}) and (\ref{Hal07}), we get $f_n \in \cL((m+g-1)P_{s+1})$.
Therefore
\begin{equation}\nonumber
 \nu_P(f_n  k^{(i)}_{j_1,j_2}  \dd t_{s+1})  \geq 0 \quad {\rm and }  
\quad  \underset{P}\Res(f_n  k^{(i)}_{j_1,j_2}  \dd t_{s+1})=0 \quad {\rm for\; all }  \quad 
   P \in \PP_f \setminus \{ P_i,P_{s+1}\}.
\end{equation}
Applying the Residue Theorem, we derive 
\begin{equation}  \label{Hal50}
   \underset{P_{i}}\Res(f_n  k^{(i)}_{j_1,j_2}  \dd t_{s+1})=- 
 \underset{P_{s+1}}\Res(f_n  k^{(i)}_{j_1,j_2}  \dd t_{s+1}) 
\end{equation}
for  $j_1 \geq 0,  \; j_2 \in [1, e_i]$    and    $i  \in [1, s]$.
Using (\ref{Hal49b}), we get the following local expansion
\begin{equation} \nonumber
\tau_i  k^{(i)}_{j_1,j_2}  := \; \sum_{r = -j_1 }^{\infty} \varkappa^{(i,j_2)}_{j_1,r}
	t_{i}^{r-1}  , \quad {\rm where \; all }  \quad  \varkappa^{(i,j_2)}_{j_1,r} \in \FF_b
	\; \ad \; \varkappa^{(i,j_2)}_{j_1,j_1}\neq 0
\end{equation}
for $j_1 \geq 0,  \; j_2 \in [1, e_i]$    and    $i  \in [1, s] $. 
By  (\ref{Hal12}) and (\ref{Hal50}), we obtain
\begin{equation} \nonumber
  -\underset{P_{s+1},t_{s+1}}\Res(f_n  k^{(i)}_{j_1,j_2} )=  
	\underset{P_{i},t_{i}}\Res(f_n \tau_i k^{(i)}_{j_1,j_2} ) = 
	 \underset{P_{i},t_{i}}\Res\Big(\sum_{j = 0}^{\infty} f^{(i)}_{n,j}
	t_i^{j}  \sum_{r = -j_1 }^{\infty} \varkappa^{(i,j_2)}_{j_1,r}
	t_i^{r-1}  \Big)
\end{equation}
\begin{equation}\label{Hal59}
 =\sum_{j = 0}^{\infty}  \sum_{r = -j_1}^{0} 
 \Tr_{F_{P_i}/\FF_b} ( f^{(i)}_{n,j} \varkappa^{(i,j_2)}_{j_1,r}  ) \delta_{j,-r}  =  \sum_{j = 0}^{j_1}
\Tr_{F_{P_i}/\FF_b} ( f^{(i)}_{n,j} \varkappa^{(i,j_2)}_{j_1,-j}  )=0 
\end{equation}
 for $0 \leq j_1 \leq d_i-1, \; 1 \leq j_2 \leq e_i\;    \ad \; 1 \leq i  \leq  s$.
Similarly to the proof of Lemma 14, we get from 
  (\ref{Hal59}) the assertion of Lemma 18.   \qed\\

{\bf Lemma 19.} {\it Let  $s \geq 2$, $d_0 =d+t,$ $\epsilon = \eta_1(2sd_0 e)^{-1}$,
 $\eta_1=(1+ \deg((t_{s+1})_{\infty}))^{-1}$,
\begin{equation} \nonumber
 \Lambda_1 = \Big\{ \Big(\Big( \underset{P_{s+1},t_{s+1}}\Res(f_n k^{(i)}_{j_1,j_2} ) 
\Big)_{\substack{d_{i,1} \leq j_1 \leq  d_{i,2} \\ 1 \leq j_2 \leq e_i}, 1 \leq i \leq s},
			\bar{a}_{d_{s+1,1}}(n),...,   \bar{a}_{d_{s+1,2}}(n) \Big)   |  n \in [0, b^m)           \Big\}
\end{equation}
with $e=e_1e_2\cdots e_s$, $e_{s+1}=1$, $ d_{s+1,1}  =t+(s-1)d_0 [m \epsilon]e $, 			
\begin{equation}\label{Hal60a}
  	d_{s+1,2} =t-1+sd_0 [m \epsilon]e, \;\;\; d_{i,1}=q_i, \; d_{i,2} = d_0  [m \epsilon]e/e_i-
		g-1 \;\for\;i \in [ 1,s],
\end{equation}  
and $m \geq  |2g -2 +2(t+g-2) (\eta_1^{-1} -1)|+2t+ 	2 / \epsilon $. Then}
\begin{equation}\label{Hal60b}
   		\Lambda_1 = \FF_b^{\chi} \qquad \with \qquad  \chi =\sum_{i=1}^{s+1} (d_{i,2} - d_{i,1} +1)e_i.
\end{equation}
\\

{\bf Proof.} Suppose that (\ref{Hal60b}) is not true. We get  
that  there exists $b^{(i)}_{j_1,j_2} \in \FF_b$ $(i,j_1,j_2 \geq 1)$
 such that
\begin{equation}\label{Hal64}
\sum_{i=1}^{s} \sum_{j_1=d_{i,1}}^{d_{i,2}} \sum_{j_2=1}^{e_{i}}|b^{(i)}_{j_1,j_2}| 
+\sum_{j_1=d_{s+1,1}}^{d_{s+1,2}} |b^{(s+1)}_{j_1}| >0 
\end{equation}
and
\begin{equation}\label{Hal66}
\sum_{i=1}^{s} \sum_{j_1=d_{i,1}}^{d_{i,2}} \sum_{j_2=1}^{e_{i}} b^{(i)}_{j_1,j_2} 
\underset{P_{s+1},t_{s+1}}\Res(f_n k^{(i)}_{j_1,j_2})+
\sum_{j_1=d_{s+1,1}}^{d_{s+1,2}}  b^{(s+1)}_{j_1} 
  \bar{a}_{j_1}(n) =0  
\end{equation}
$\fall \;  n \in [0, b^m)$. 
From (\ref{Hal03})-(\ref{Hal07}), we obtain the  following local expansion
\begin{equation}\label{Hal66a}
   f_n =\dot{f}_n +\ddot{f}_n =\sum_{r \leq   m+g-1 } f^{(s+1)}_{n,r}
	t_{s+1}^{-r},   \quad \with  \quad
 \ddot{f}_n=  \sum_{i=g}^{m-1} \bar{a}_{i}(n)v_i,
\end{equation}
and $\dot{f}_n=  \sum_{i=0}^{g-1} \bar{a}_{i}(n)v_i$, $\where \;  n \in [0, b^m)$. Let $r \geq g$.\\
Using (\ref{Hal03})-(\ref{Hal07}) and   (\ref{Hal18}), we derive  that  $\nu_{P_{s+1}}(\dot{f}_n) \geq -2g+1$, 
$\nu_{P_{s+1}}(\dot{f}_n t_{s+1}^{r+g-1})$ $ \geq 0$ and
\begin{equation}\nonumber
 f^{(s+1)}_{n,r+g}  = \underset{P_{s+1},t_{s+1}}\Res(f_n t_{s+1}^{r+g-1}) = 	 
\underset{P_{s+1},t_{s+1}}\Res(\ddot{f}_n t_{s+1}^{r+g-1}) = 
	\underset{P_{s+1},t_{s+1}}\Res \Big( \sum_{i=g}^{m-1} \bar{a}_{i}(n) 
\end{equation}
\begin{equation}\nonumber
  \times \sum_{j \leq i+g}v_{i,j}
	t_{s+1}^{-j+ r+g-1}  \Big) = \sum_{i=g}^{m-1} \bar{a}_{i}(n) \sum_{j \leq i+g}v_{i,j}
	\delta_{j,r+g} = \sum_{m-1 \geq i  \geq r} \bar{a}_{i}(n) v_{i,r+g} 
	\; \for \; r \geq g	.
\end{equation}
Taking into account that $v_{i,i+g} =1$ and $v_{i,r+g} =0$ for $i>r \geq g $ (see (\ref{Hal18a})),
we get 
\begin{equation}  \label{Hal67}
   f^{(s+1)}_{n,r+g} = \bar{a}_{r}(n) \quad
	\for \quad r \geq g \quad  \ad \quad n \in [0,b^m).
\end{equation}
By (\ref{Hal66}),  we have
\begin{equation}\nonumber
\sum_{i=1}^{s} \sum_{j_1=d_{i,1}}^{d_{i,2}} \sum_{j_2=1}^{e_{i}} b^{(i)}_{j_1,j_2} 
\underset{P_{s+1},t_{s+1}}\Res(f_n k^{(i)}_{j_1,j_2})+
\sum_{j_1=d_{s+1,1}}^{d_{s+1,2}}  b^{(s+1)}_{j_1} 
\underset{P_{s+1},t_{s+1}}\Res(f_n t_{s+1}^{j_1+g-1} ) =0 
\end{equation}
 for  all  $n \in [0, b^m)$. 
Hence
\begin{equation}\label{Hal68}
\underset{P_{s+1},t_{s+1}}\Res(f_n  \alpha )  =0 \quad {\rm for \;
	all}  \quad n \in [0, b^m), \quad \where \quad \alpha = \alpha_1 + \alpha_2,
\end{equation}
\begin{equation}\nonumber
 \alpha_{1}= 	\sum_{i=1}^{s} \alpha_{1,i},\;\; \alpha_{1,i}= \sum_{j_1=d_{i,1}}^{d_{i,2}} \sum_{j_2=1}^{e_{i}} b^{(i)}_{j_1,j_2}  k^{(i)}_{j_1,j_2}, \;\; \ad \;\;  \alpha_2=
\sum_{j_1=d_{s+1,1}}^{d_{s+1,2}}  b^{(s+1)}_{j_1}
t_{s+1}^{j_1+g-1}.
\end{equation}
According to (\ref{Hal47}), we get the following local expansion
\begin{equation} \nonumber
  k^{(i)}_{j_1,j_2}  := \; \sum_{r = g+1}^{\infty} \varkappa^{(i,j_2)}_{j_1,r}
	t_{s+1}^{r-1}  , \quad {\rm where \; all }  \quad  \varkappa^{(i,j_2)}_{j_1,r} \in \FF_b,
\end{equation}
and 
\begin{equation}  \label{Hal72}
          \alpha = \sum_{r= g+1 }^{\infty} \varphi_r t_{s+1}^{r-1}
             \quad {\rm with   } \quad  \varphi_r \in\FF_b, \quad r \geq g+1.
\end{equation}
Using (\ref{No26}) and (\ref{Hal66a})-(\ref{Hal68}), we have
\begin{equation}\nonumber
  \underset{P_{s+1},t_{s+1}}\Res(f_n \alpha) =
	\underset{P_{s+1},t_{s+1}}\Res\Big(\sum_{j \leq m+g-1}f_{n,j}^{(s+1)}  t_{s+1}^{-j}  \sum_{r=g+1}^{\infty} \varphi_r  t_{s+1}^{r-1}\Big)
\end{equation}
\begin{equation}\nonumber
= \sum_{j \leq m+g-1}f_{n,j}^{(s+1)}   \sum_{r=g+1}^{\infty} \varphi_r  \delta_{j,r}
           = \sum_{j=g+1}^{m+g-1} f_{n,j}^{(s+1)} \varphi_{j}   =
	 \sum_{r=g+1}^{m+g-1}  \bar{a}_{r}(n)  \varphi_{r} =0 .
\end{equation}
for $n \in [0, b^m))$.
Hence
\begin{equation} \nonumber
         \varphi_{r} =0 
             \quad {\rm for   } \quad  g+1 \leq r \leq m+g-1.
\end{equation}
By (\ref{Hal72}), we obtain 
\begin{equation}  \nonumber
          \nu_{P_{s+1}} (\alpha) \geq m+g-1.
\end{equation}
Applying (\ref{Hal40}), (\ref{Hal46}) and (\ref{Hal68}), we derive 
\begin{equation} \nonumber
           \alpha \in \cL(G_1), \; \with \; G_1=W +\sum_{i=1}^s d_{i,2}P_i  + (d_{s+1,2}+g-1) ( t_{s+1})_{\infty}
					- (m+g-1) P_{s+1} .
\end{equation}
From (\ref{Hal60a}), we have
\begin{equation}\nonumber
\deg(G_1)   = 2g-2 + \sum_{i=1}^s d_{i,2} e_i   
					+  (d_{s+1,2} +g-1) \deg((t_{s+1})_{\infty}) -(m+g-1) 
\end{equation}
\begin{equation}\nonumber
 \leq 2g-2 +sd_0e [m \epsilon]  +(t-1 + sd_0e [m \epsilon] +g-1) (\eta_1^{-1} -1) -(m+g-1)
\end{equation}
\begin{equation}\nonumber
\leq g -1 +(t+g-2) (\eta_1^{-1} -1)   +sd_0e m \epsilon \eta_1^{-1}  -m= g -1 +(t+g-2) (\eta_1^{-1} -1)  -m/2  <0 
\end{equation}
for $m > 2g -2 +2(t+g-2) (\eta_1^{-1} -1)$.
Hence $\alpha=0$.\\

%
Suppose that $\sum_{i=1}^{s}\sum_{j_1=d_{i,1}}^{d_{i,2}}\sum_{j_2=1}^{e_{i}}|b^{(i)}_{j_1,j_2}|=0$.
 Then $\alpha_2 =0$. From (\ref{Hal68}), we derive $ b^{(s+1)}_{j_1}  =0$ for all 
 $ j_1 \in [d_{s+1,1},d_{s+1,2}]$.  
According to (\ref{Hal64}), we have a contradiction.
Hence there exists $h \in [1,s]$ with 
\begin{equation}\label{Hal81}
          \sum_{j_1=d_{h,1}}^{d_{h,2}}\sum_{j_2=1}^{e_{h}}|b^{(h)}_{j_1,j_2}|  > 0.
\end{equation}
Let $h>1$. By  (\ref{Hal17}) and (\ref{Hal68}), we get
 $\nu_{P_h}(t_{s+1}) \geq 0$ and  $\nu_{P_h}(\alpha_2) \geq 0$.
Applying  (\ref{No08}) and (\ref{No12}), we derive $\nu_{P_h}(W) =\nu_{P_h}(\dd t_{s+1}) =
\nu_{P_h}(\dd t_{s+1}/\dd t_h) \geq 0$.

By (\ref{Hal49c}), we have 
 $\nu_{P_h}(\alpha_{1,j}) \geq - \nu_{P_h}(W)$  for $1 \leq j \leq s, j \neq h$.
 Taking into account that $\alpha_{1,h} = -\sum_{1 \leq j \leq s, j \neq h} \alpha_{1,j} -\alpha_2$, we get
 $\nu_{P_h}(\alpha_{1,h}) \geq - \nu_{P_h}(W)$ .\\

 Using (\ref{Hal49a}) and (\ref{Hal49b}), we obtain 
$\nu_{P_h}(k^{(h)}_{j_1,j_2}) =-j_1 -1- \nu_{P_h}(W)$.
 Bearing in mind (\ref{Hal81}) and that
$ \{u^{(i)}_{1},...,u^{(i)}_{z_i}, k^{(i)}_{q_i,1},...,k^{(i)}_{q_i,e_1},...,
	 k^{(i)}_{n,1},...,k^{(i)}_{n,e_1}			\}$ is a basis of $\cL(G_i + (n-q_i+1) P_i)$, we get
\begin{equation}\nonumber
     \alpha_{1,h} \in \cL(G_i + (d_{i,2}-q_i+1) P_i) \setminus
	\cL(G_i + (d_{i,1}-q_i) P_i). 
\end{equation}	
From  (\ref{Hal60a}) and (\ref{Hal68}), we derive  $\nu_{P_h}(\alpha_{1,h}) \leq - \nu_{P_h}(W)  -1$. 
We have a contradiction.\\

Now let $h=1$ and (\ref{Hal81})  is not true for $h \in [2,s]$.
Hence $\alpha_{1,1} =-\alpha_2  $ and 
 $\nu_{P_{s+1}} (\alpha_{1,1}) \geq d_{s+1,1} +g-1$.
By (\ref{Hal40}), (\ref{Hal46}) and (\ref{Hal68}), we have 
\begin{equation}\nonumber
           \alpha_{1,1} \in \cL\big( \dot{G} \big) 
					\quad \with \quad \dot{G}=  W +(d_{1,2}+1)P_{1} - (d_{s+1,1} +g-1)P_{s+1}.
\end{equation}
From (\ref{Hal60a}), we get
\begin{equation}\nonumber
    \deg(\dot{G})   = 2g -2+ d_0e [m \epsilon] - ge_1 - (s-1) d_0e [m \epsilon] -g+1
		 \leq 2g-2-2g+1<0.
\end{equation}
Hence $\alpha_{1,1}=0$.  Therefore (\ref{Hal81})  is not true for $h=1$. 
We have  a contradiction.  
Thus  assertion (\ref{Hal64})  is not true, and  Lemma 19 follows. \qed\\
\\
{\bf End of the proof of Theorem 4.}

Let $ \tilde{d}_{i,2} = d_{i,2} +g =    d_0  [m \epsilon]e/e_i -1\;   (1 \leq i \leq s)$,
\begin{equation} \nonumber
 \Lambda_1^{'} = \Big\{ \Big(\big( \underset{P_{s+1},t_{s+1}}\Res(f_n k^{(i)}_{j_1,j_2} ) \big)_{0 \leq j_1 \leq  \tilde{d}_{i,2}, 1 \leq j_2 \leq e_i, 1 \leq i \leq s},
			\bar{a}_{d_{s+1,1}}(n),...,   \bar{a}_{d_{s+1,2}}(n) \Big)   \Big|  n \in [0, b^m) \Big\}
\end{equation}
and 
\begin{equation}\nonumber
  		\Lambda_2 = \Big\{ (\bar{a}_{d_{s+1,1}}(n),...,   \bar{a}_{d_{s+1,2}}(n) ) 
 		\; \Big|  \;  \underset{P_{s+1},t_{s+1}}\Res(f_n  k^{(i)}_{j_1,j_2} )=0 
				\qquad\qquad\qquad\qquad \qquad\qquad
\end{equation}
\begin{equation} \nonumber		
		\qquad\qquad\qquad\qquad \qquad\qquad \for \quad  0 \leq j_1 \leq  \tilde{d}_{i,2}, 1 \leq j_2 \leq e_i,
		1 \leq i \leq s, \;			
	n \in [0, b^m)  \Big\}. 
\end{equation}
By (\ref{NiOz250}) and Lemma 19, we have  $\dim_{\FF_b}(\Lambda_1^{'}) \geq \dim_{\FF_b}(\Lambda_1)$ and 
\begin{equation}\nonumber
			\dim_{\FF_b}(\Lambda_2)= \dim_{\FF_b}(\Lambda_1^{'}) -
			\dim_{\FF_b}\Big( \Big\{ \big( \underset{P_{s+1},t_{s+1}}\Res(f_n k^{(i)}_{j_1,j_2} ) \big)_{\substack{0 \leq j_1 \leq  \tilde{d}_{i,2}, 1 \leq j_2 \leq e_i\\ 1 \leq i \leq s}}   \Big|  n \in [0, b^m)           \Big\} \Big)
\end{equation}
\begin{equation} \label{Hal100}
		\geq \dim_{\FF_b}(\Lambda_1) -\sum_{i=1}^{s} (\tilde{d}_{i,2}  +1)e_i  
		\geq d_{s+1,2} -d_{s+1,1}+1  - \sum_{i=1}^{s} (q_i  +  g) e_i.
\end{equation}
 Using Lemma 18, we get $\Lambda_3 \supseteq \Lambda_2$ and 
$\dim_{\FF_b}(\Lambda_3) \geq \dim_{\FF_b}(\Lambda_2)$, where
\begin{equation}\nonumber	
 		\Lambda_3 = \Big\{ 			
	(\bar{a}_{d_{s+1,1}}(n),...,   \bar{a}_{d_{s+1,2}}(n) ) 		
 		\; \Big| \;  \Tr_{F_{P_i}/\FF_b}(\beta^{(i)}_{j_2}f_{n, j_1}^{(i)} ) =0
			 \qquad\qquad\qquad \qquad\qquad
		\qquad\qquad\qquad\qquad
\end{equation}		
\begin{equation} \nonumber		
\qquad\qquad	\qquad\qquad \for \; 0 \leq j_1 \leq  \tilde{d}_{i,2}, 1 \leq j_2 \leq e_i, 
	1 \leq i \leq s, \;				n \in [0, b^m)  \Big\}. 
\end{equation}

Taking into account that 
 $(\bx_n)_{0 \leq n <b^m} $ is a  $(t,m,s)$ net in base $b$, we  get from (\ref{Hal14}) and
 (\ref{Hal16})  that
\begin{equation}\nonumber
 \Big\{ 			
\Big( f_{n, j_1}^{(i)} )
\Big)_{0\leq j_1 \leq  \tilde{d}_{i,2}, 1 \leq i \leq s}
	 \; \Big| \; 		
					n \in [0, b^m)  \Big\} = \prod_{i=1}^s F_{P_i}^{\tilde{d}_{i,2} +1}. 	
\end{equation}		
Bearing in mind that $\{\beta^{(i)}_{1}, ..., \beta^{(i)}_{e_i}\}$ is a basis of $F_{P_i}/\FF_b$ (see Lemma 18), we obtain
\begin{equation}\nonumber
  		\Lambda_4 = \Big\{ 			
\Big( \Tr_{F_{P_i}/\FF_b}(\beta^{(i)}_{j_2}f_{n, j_1}^{(i)} )
\Big)_{0 \leq j_1 \leq  \tilde{d}_{i,2}, 1 \leq j_2 \leq e_i,1 \leq i \leq s}
	 \; \Big| \; 		
					n \in [0, b^m)  \Big\} =\FF_b^{sd_0 e[m\epsilon]}.
\end{equation}		
Let
\begin{equation}\nonumber
   		\Lambda_5 =  \Big\{ 			
 \Big( \Tr_{F_{P_i}/\FF_b}(\beta^{(i)}_{j_2}f_{n, j_1}^{(i)} ) \Big)_{0\leq j_1 \leq  \tilde{d}_{i,2}, 1 \leq j_2 \leq e_i,1 \leq i \leq s}, 
  	\big( \bar{a}_{j}(n) \big)_{d_{s+1,1} \leq j \leq  d_{s+1,2}}	 \; \Big| \; 		
					n \in [0, b^m) 	\Big\}.
\end{equation}			
By (\ref{Hal100}), (\ref{NiOz250}) and (\ref{Hal40}), we have  
\begin{equation} \nonumber
			\dim_{\FF_b}(\Lambda_5) = \dim_{\FF_b}(\Lambda_3) +\dim_{\FF_b}(\Lambda_4)  \geq d_{s+1,2} -d_{s+1,1}+1 +sd_0e \dot{m} 		-r
\end{equation}
with $r=(g+1)(e_0 +s) $, $e=e_1e_2...e_s$ and  $\dot{m} = [m \epsilon]$.

Let
$\dot{m}_{1} = d_0 e \dot{m}$, $\epsilon = \eta_1(2sd_0 e)^{-1}$,
  $\dddot{m}_i=0$,	$1 \leq i \leq  s ,$  and 
$\dddot{m}_{s+1} =d_{s+1,1}+g  $, 
 $ d_{s+1,1}  =t+(s-1)d_0 [m \epsilon]e $, 			
  	$d_{s+1,2} =t-1+sd_0 [m \epsilon]e =d_{s+1,1}+\dot{m}_{1}-1$ (see (\ref{Hal60a})),  $\tilde{d}_{i,2} = d_0  [m \epsilon]e/e_i-1$$=d_{i,2}+g=\dot{m}_{1}/e_i-1$
		$(i \in [ 1,s])$, 
\begin{equation}\nonumber
   \dot{\theta}^{(i)}_{n, j_1 e_s +j_2}:
 = \Tr_{F_{P_i}/\FF_b}(\beta^{(i)}_{j_2}f_{n, j_1}^{(i)} )  	
			\qquad \ad \qquad \dot{\theta}^{(s+1)}_{n, j+1}: = f^{(s+1)}_{n, j} 
 = \bar{a}_{ j-g} (n) \;\; {\rm (see \;(\ref{Hal67}))}
\end{equation}
for $0 \leq j_1 \leq  \tilde{d}_{i,2},\; 1 \leq j_2 \leq e_i, \; 1 \leq i \leq s,  2g \leq j $, 	
and let
\begin{equation}\nonumber
   		\Lambda_6 = \Big\{ \Big(
			\Big( 	\dot{\theta}^{(i)}_{\dddot{m}_i+d_0 e \dot{j}_i+ \ddot{j}_i}  \Big)_{
			0 \leq \dot{j}_i <  \dot{m},	1 \leq \ddot{j}_i \leq  d_0 e, 1 \leq i \leq s+1}
\; \Big| \; 	n \in [0, b^m) \Big\}.
\end{equation}
It is easy to verify that $\Lambda_6=\Lambda_5$ and $\dim_{\FF_b}(\Lambda_6)=(s+1)\dot{m}_1 -\dot{r}$
with $0 \leq \dot{r} \leq r=(g+1)(e_0 +s)$.

Let $ m \geq  |2g -2 +2(t+g-2) (\eta_1^{-1} -1)|+2t+ 	2 / \epsilon$. Applying    Lemma 2, with $\dot{s} =s+1$, we get that there exists  $B_i \subset \{0,...,\dot{m}-1\}$, $1 \leq i \leq s+1$ 
such that 
\begin{equation}\nonumber	
    \Lambda_7 =\FF_b^{(s+1)\dot{m}_1-d_0e B}, \qquad \where \quad  B =\#B_1+...+\#B_{s+1} \leq (g+1)(e_0 +s),
\end{equation}
and  
\begin{equation} \nonumber
  		  \Lambda_7  = \Big\{   \Big( 	\dot{\theta}^{(i)}_{\dddot{m}_i+d_0 e \dot{j}_i+ \ddot{j}_i}   \; \Big| \; 
	 \dot{j}_i  \in \bar{B}_{i}, \; \ddot{j}_i \in [1, d_0 e],  \;     i \in [1,s+1] \Big)
			 \; \Big| \; 		
					n \in [0, b^m) 	\Big\} 
\end{equation}
with $\bar{B}_i =\{0,...,\dot{m}-1\} \setminus B_{i}$.
Hence
\begin{equation}\nonumber
			\Big\{   \Big( 	f^{(i)}_{n,\dddot{m}_i+ \dot{j}_id_0e/e_i+ \ddot{j}_i -1}     \Big|  
	 \dot{j}_i  \in \bar{B}_{i},  \ddot{j}_i \in [1, \frac{d_0 e}{ e_i}],       i \in [1,s+1] \Big) 
				\Big|	n \in [0, b^m) 	\Big\}=\prod_{i=1}^s F_{P_i}^{\chi_i}\FF_b^{\chi_{s+1}}	
\end{equation}
with  $e_{s+1} =1$, $\chi_i=d_0e(\dot{m} -\#B_i)/e_i$, $1 \leq i \leq s+1$.\\  
Taking into account that $\sigma_{P_i}: F_{P_i} \to Z_{b^{e_i}}$ is a bijection (see (\ref{Ha08})), we obtain
\begin{equation}\nonumber
			\Big\{   \Big( \sigma_{P_i}(	f^{(i)}_{n,\dddot{m}_i+ \dot{j}_id_0e/e_i+ \ddot{j}_i -1})    \; \Big| \; 
	 \dot{j}_i  \in \bar{B}_{i},  \ddot{j}_i \in [1, \frac{d_0 e}{ e_i}],       i \in [1,s]\Big), 	
			\qquad \qquad \qquad 	\qquad \qquad 			
\end{equation}
\begin{equation}\nonumber
 \big(	 a_{\dddot{m}_{s+1}+ \dot{j}_{s+1}d_0e+ \ddot{j}_{s+1} -1-g}(n)   \big|  \dot{j}_{s+1} \in \bar{B}_{s+1},  \ddot{j}_{s+1} \in [1, d_0 e]	\big)  \Big|  		
					n \in [0, b^m) 	\Big\}= Z_b^{(s+1)\dot{m}_1-d_0e B}. 	
\end{equation}
Let $\tilde{B}_{i} =\bar{B}_{i}$, $1 \leq i \leq s$, and let $\tilde{B}_{s+1} =\{ \dot{m}-j-1| j \in  \bar{B}_{s+1} \}$.
From (\ref{Hal25}), we derive
\begin{equation}\nonumber
			\Big\{   \Big( 	x^{(i)}_{n,\ddot{m}_i+ \dot{j}_id_0e+ \ddot{j}_i -1}    \; \Big| \; 
	 \dot{j}_i  \in \tilde{B}_{i},  \ddot{j}_i \in [1, d_0 e ],       i \in [1,s+1]\Big)
	 \Big|  					n \in [0, b^m) 	\Big\}= Z_b^{(s+1)\dot{m}_1-d_0e B}			 , 	
\end{equation}
where $x^{(s+1)}_n=\sum_{j=1}^{m} x^{(s+1)}_{n,j} b^{-j} :=n/b^{m}$, and $x^{(s+1)}_{n,j} =a_{m-j-1}(n)$ $(1 \leq j \leq m)$, $\ddot{m}_i = \dddot{m}_i=0$ for $1 \leq i \leq s$ and 
$\ddot{m}_{s+1} = m-t-s\dot{m}_1= m-1-(\dddot{m}_{s+1}+ \dot{m}_1  -1-g)$.

By Lemma 17 and Theorem K, we obtain that $ (\bx_n)_{ n \geq 0}$ is a $d-$admissible
  $(t, s) $  sequence with $\bx_n= (x^{(1)}_n,...,x^{(s)}_n) $, $d=g+e_0$ and $t=g+e_0-s$.\\
Now  applying Corollary 1   with $\dot{s}=s+1$, $\tilde{r}=0$,  $\tilde{m}=m$  and 
 $\hat{e}=e=e_1  ... e_{s+1}$, we derive
\begin{equation}\nonumber
    \min_{0 \leq Q <b^m}  \min_{ \bw \in E_m^s}   b^m \emph{D}^{*}((\bx_n \oplus \bw, n\oplus Q/b^m)_{0 \leq  n < b^m}) \geq   
		2^{-2} b^{-d}K_{d,t,s+1}^{-s} \eta^{s}_1 m^{s},
\end{equation}\\ 
with $ m \geq 2^{2s+3} b^{d+t+s+1}(d+t)^{s+1}s^{2s} e(g+1)(e_0+s)\eta_1^{-s}$, 
 and $\eta_1 =(1+\deg((t_{s+1})_\infty))^{-1} $.
Using Lemma B, we get the  
 first assertion in Theorem~4.\\

Consider the second assertion in Theorem 4.\\
By (\ref{Hal12})-(\ref{Hal16}), we get that the net  $(\bx_n)_{0 \leq n < b^m}$ is constructed similarly to the construction  of the Niederreiter-\"{O}zbudak net (see (\ref{NiOz18})-(\ref{NiOz24}) and 
(\ref{NiOz08b0})).
 The difference is t
hat in the construction of Section 2.3 the map 
$\sigma_i:F_{P_i} \to \FF_b^{e_i}$ is linear, while in 
the construction of  Section 2.4 this map 
may be nonlinear.

It is easy to verify that this 
 does not affect the proof of 
 bound (\ref{NiXi07}) and Theorem 4 follows
. \qed



\subsection{Niederreiter-Xing sequence. Sketch of the proof of Theorem 5}
%
First we will prove that
\begin{equation}\label{NiXi10}
\dot{\cC}_m =\cM_m^{\bot}(P_1,...,P_s; G_m) \quad  \for  \quad m \geq g+1 .
\end{equation}
By  (\ref{Ap309}) and (\ref{NiXi05}), we get
\begin{equation}  \nonumber
\dot{\cC}_m = \Big\{  \Big(\sum_{r=0}^{m-1}   \dot{c}^{(i)}_{j,r}  \bar{a}_r(n) 
\Big)_{0 \leq j
\leq m-1, 1 \leq i \leq s}   \; \Big| \; 0 \leq n < b^m   \Big\}.
\end{equation}
Using (\ref{NiOz23b}) with $\tilde{G}=(g-1)P_{s+1}$, we derive $G_m^{\bot} = L_m$,
 where $L_m =\cL((m-g+1) P_{s+1} +W)$.
From   (\ref{NiXi11a}), we have
\begin{equation}  \nonumber
 \{f^{\bot} \; |\; f^{\bot} \in  L_m\} =
 \{\dot{f}_n:= \sum_{  r = 0}^{m-1} a_r(n) \dot{v}_r\; |\; n \in [0,b^m) \}.
\end{equation}
Applying (\ref{NiXi05}), we obtain 
\begin{equation} \nonumber
   \dot{f}_n \tau_i= \; \sum_{j = 0}^{\infty} \dot{f}_{n,j}^{(i)}
	 t_i^{j}  , \quad {\rm where  } 	\quad 
	\dot{f}_{n,j}^{(i)} =\sum_{r=0}^{m-1}     \dot{c}^{(i)}_{j,r}  \bar{a}_r(n) 
	\in \FF_b, \; i \in[1, s],\;  j \geq 0.
\end{equation}
Therefore
\begin{equation}\label{NiXi11b}
\dot{\cC}_m = \{  (\dot{f}_{n,j}^{(i)} )_{0 \leq j
\leq m-1, 1 \leq i \leq s}   \; | \; 0 \leq n < b^m   \}.
\end{equation}
We use notations (\ref{NiOz23c})-(\ref{NiOz24}) with the following modifications.
In (\ref{NiOz18}) we take the field $\FF_b$ instead of $F_{P_i}$ ,
and in (\ref{NiOz20}) we consider the map $ \vartheta_{i}^{\bot}$ as the identical map $(1 \leq i \leq s)$.
By (\ref{NiOz20a}), we have   $\dot{\theta}^{\bot}_{i,j}(f_n) =\dot{f}_{n,j-1}^{(i)}$ for $ 1 \leq j \leq m$,  and
$\dot{\theta}^{\bot}_i(\dot{f}_n) =(\dot{f}_{n,0}^{(i)},...,\dot{f}_{n,m-1}^{(i)} )$, $1 \leq i \leq s$. 
 According to (\ref{NiOz24}) and (\ref{NiXi11b}) we get
\begin{eqnarray*}
       \Xi_m =  \dot{\Xi}_m = \{   \dot{\theta}^{{\bot}}(f^{\bot})    |  f^{\bot} \in 
				\cL(G_m^{\bot}) \} =
				\{   \dot{\theta}^{{\bot}}(\dot{f}_n)   |   n \in [0, b^m) \} 
				\qquad\qquad\qquad	\qquad\qquad \\
				\qquad=	
				\{    (\dot{\theta}^{\bot}_1(\dot{f}_n),...,\dot{\theta}^{\bot}_s(\dot{f}_n))   |    n \in [0, b^m) \} =
			\{  (\dot{f}_{n,j}^{(i)} )_{0 \leq j
\leq m-1, 1 \leq i \leq s}   \; | \; 0 \leq n < b^m   \}=	\dot{\cC}_m				.
\end{eqnarray*}
Now applying (\ref{NiOz07}), (\ref{NiXi03}) and Lemma 12,  we obtain  (\ref{NiXi10}).
By \cite[ref. 8.9]{DiPi}, we have
\begin{equation}\nonumber
        \delta_m(\cM_m) = \delta_m(\cM_m(P_1,...,P_s; G_m)) \geq m-g+1 \quad \for \quad  \;m \geq g+1.
\end{equation}
Taking into account Proposition B, we get that 
$\bx_n(\dot{C})_{n \geq 0}$ is a digital $(\bT,s)$ sequence  with  
$T(m)=g$   for  $m \geq  g+1$.

Now the $d-$admissible property follow from Lemma 16. 
In order to complete the proof of Theorem 5, we use Theorem 3 and Theorem 4. \qed \\



\subsection{ General  $d-$admissible  $(t,s)$ sequences. Proof of Theorem 6.}
 First we will prove  Lemma 20. We need the following notations: \\

  Let $\tilde{C}^{(1)},...,\tilde{C}^{(\dot{s})}$ are $m \times m$  generating matrices of a 
	digital  $(t,m,{\dot{s}})$ net $ (\tilde{\bx}_n)_{n=0}^{b^m-1} $ in base $b$, $\tilde{x}_n^{(\dot{s})}\neq \tilde{x}_k^{(\dot{s})} $ for $n \neq k$, 
		$\tilde{C}^{(i)} = (\tilde{c}_{r,j}^{(i)})_{1 \leq r,j \leq m}$, 
%
$ \tilde{\fc}^{(i)}_j =(\tilde{c}_{1,j}^{(i)},...,\tilde{c}_{m,j}^{(i)}) \in \FF^m_b$, $i \in [1,\dot{s}]$,
$ \tilde{\fc}_j =(\tilde{\fc}^{(1)}_j,...,\tilde{\fc}^{(\dot{s})}_j) \in \FF^{m \dot{s}}_b$   
 $(1 \leq j \leq m )$.
%
Let 
$\phi : \; Z_b  \mapsto \FF_b$ be a  bijection with
$\phi(0) = \bar{0}$, and let 	  	
	$n =\sum^m_{j=1} a_j(n) b^{j-1}$, $ \bn = (\bar{a}_1(n),...,\bar{a}_m(n)) \in  \FF^m_b $, 
	$\bar{a}_j(n) =\phi(a_j(n))$, 
	   $\tilde{\by}_n =(\tilde{\by}^{(1)}_n,...,\tilde{\by}^{({\dot{s}})}_n) \in \FF^{m \dot{s}}_b $,  
		$\tilde{\by}^{(i)}_n = 	(\tilde{y}^{(i)}_{n,1},...,\tilde{y}^{(i)}_{n,m}) \in  \FF^m_b$, 
\begin{equation} \label{6Di00}
    \tilde{\bx}_n =(\tilde{x}^{(1)}_n,...,\tilde{x}^{({\dot{s}})}_n), \quad \tilde{x}^{(i)}_n =\sum_{j=1}^m \phi^{-1} (\tilde{y}^{(i)}_{n,j})/b^j \quad \for \quad 1 \leq i \leq \dot{s},
\end{equation} 
\begin{equation} \label{6Di01}
    		\tilde{\by}^{(i)}_n = 	\bn (\tilde{\fc}_{1}^{(i)},...,\tilde{\fc}_{m}^{(i)})^{\top} :=
				\sum_{j=1}^m \bar
				{a}_j(n) \tilde{\fc}_j^{(i)}  = \bn  \tilde{C}^{(i)\top} \quad
				\for \quad  1 \leq i \leq \dot{s}.
\end{equation}	
Hence
\begin{equation} \nonumber
	 \tilde{\by}_n =\sum_{j=1}^m \bar{a}_j(n) \tilde{\fc}_j    , \quad \for \quad  
              0 \leq n < b^m.
\end{equation}
We put 
\begin{equation} \nonumber
\tilde{\Phi}_m =\{ \tilde{\bx}_n|   n \in[0, b^m) \}, \;
	\tilde{\Psi}_m =\{  \tilde{\by}_n 	 |   n \in[0, b^m) \},\;  \tilde{Y}_m =\{  \tilde{\by}_n^{({\dot{s}})}
	  |   n \in[0, b^m) \}.
\end{equation}
We see that $\tilde{\Psi}_m $ is a vector space over $\FF_b$, with 
 $\dim (\tilde{\Psi}_m)  \leq m$. 
Taking into account that
$\tilde{x}_n^{(\dot{s})}\neq \tilde{x}_k^{(\dot{s})} $ for $n \neq k$, we obtain $\dim (\tilde{\Psi}_m ) =m$, 
$\tilde{\fc}_1,...,\tilde{\fc}_m$ is the basis of $\tilde{\Psi}_m$ and  $\tilde{Y}_m =\FF_b^m $. \\
Let $d \geq 1$,	$d_0=d+t$, $m \geq 4d_0(s+1)$,	$\dot{m} =[(m-t)/(2d_0({\dot{s}}-1))]$,
\begin{equation} \label{6Di04}
  d^{(\dot{s})}_1 =m -t+1 -({\dot{s}}-1)d_0\dot{m} \quad \;\; \ad  \;\;
	\quad d^{(\dot{s})}_2 =m -t -({\dot{s}}-2)d_0\dot{m} .
\end{equation}
Bearing in mind that $\tilde{\Phi}_m $ is a $(t,m,{\dot{s}})$ net, we get that
 for each $j \in [1, (\dot{s}-1)d_0\dot{m}]$ with $j=(j_1-1)(\dot{s}-1)+j_2$, $j_1 \in [1, d_0\dot{m}]$ and $ j_2 \in [1,\dot{s}-1]$
 there exists 
$n(j) \in[0,b^m)$  such that
\begin{equation}\label{6Di07a}
      \tilde{x}_{n(j),r_1}^{(\dot{s})}  = \delta_{(j_1-1)(\dot{s}-1)+j_2,r_1}
		\quad	\quad  {\rm and  } \quad \quad
			 \tilde{x}_{n(j),r_2}^{(i)}    = \delta_{i,j_2}
		\delta_{j_1,r_2}   
\end{equation}
for all $ r_1 \in [1,({\dot{s}}-1)d_0\dot{m}],$  $ r_2 \in [1,d_0\dot{m}]$,  $ i \in [1,\dot{s}-1]$.

Taking into account that $Y_m =\FF_b^m $, we derive that
  there exists $n(j) \in[0,b^m)$ with
\begin{equation}\label{6Di07}
      \tilde{y}_{ n(j),r}^{(\dot{s})}   = \delta_{j,r} \quad \for  \quad 
			(\dot{s}-1)d_0\dot{m} + 1 \leq  j \leq m, \;\; 
 	\; 1 \leq  r \leq m.
\end{equation}

We take a basis $\dot{\ff}_1,...,\dot{\ff}_m $ of $\tilde{\Psi}_m$ in the following way:\\

Let $ \dot{\ff}_j =(\dot{\ff}^{(1)}_j,...,\dot{\ff}^{(\dot{s})}_j) \in \FF^{m \dot{s}}_b$ 
with  
$ \dot{\ff}^{(i)}_j =(\dot{\ff}_{1,j}^{(i)},...,\dot{\ff}_{m,j}^{(i)}) \in \FF^m_b$,
$ i \in [1,\dot{s}]$, $ j \in [1,m]$.\\
For $j \in [1, m]$, we put
 $\dot{\ff}_j:=  \tilde{\by}_{n(j)}$.
We have from (\ref{6Di07a}) and (\ref{6Di07}) that
\begin{equation}\nonumber
      \dot{\ff}_{(j_1-1)(\dot{s}-1)+j_2,r_1}^{(\dot{s})} = \delta_{(j_1-1)(\dot{s}-1)+j_2,r_1} \qquad 
				 {\rm and  } \qquad 
			 \dot{\ff}_{(j_1-1)(\dot{s}-1)+j_2,r_2}^{(i)} = \delta_{i,j_2}
		\delta_{j_1,r_2} 	
\end{equation}
for $ r_1 \in [1,({\dot{s}}-1)d_0\dot{m}],$  $ r_2 \in [1,d_0\dot{m}]$,  
$ i \in [1,\dot{s}-1]$,         $j_1 \in [1, d_0\dot{m}]$, $ j_2 \in [1,\dot{s}-1]$
and
\begin{equation}\label{6Di09}
      \dot{\ff}_{j,r}^{(\dot{s})} = \delta_{j,r} \quad\quad \for  \quad \quad
			(\dot{s}-1)d_0\dot{m} + 1 \leq  j \leq m, 
 	\; 1 \leq  r \leq m.
\end{equation}
 
It is easy to see that the vectors 
$\dot{\ff}_{1},..., \dot{\ff}_{m} \in \tilde{\Psi}_m$  are linearly independent over $\FF_b$.
Thus $\dot{\ff}_{1},..., \dot{\ff}_{m}$ is a basis of $ \tilde{\Psi}_m$. \\
Let
\begin{equation} \label{6Di11}
    		\dot{\by}^{(i)}_n = 	(\dot{y}^{(i)}_{n,1},...,\dot{y}^{(i)}_{n,m}): = 
				\bn(\dot{\ff}^{(i)}_1,...,\dot{\ff}^{(i)}_m)=
			\sum_{j=1}^m \bar{a}_j(n)\dot{\ff}^{(i)}_j
 = \bn  \dot{\cF}^{(i) \;\top}, 
\end{equation}	
where
$\dot{\cF}^{(i)} = (\dot{\ff}_{r,j}^{(i)})_{1 \leq r,j \leq m}$ for $1 \leq i \leq \dot{s}$.
 Hence
\begin{equation} \nonumber
	\dot{\by}_n := 	(\dot{\by}^{(1)}_n,...,\dot{\by}^{(\dot{s})}_n) =\sum_{j=1}^m \bar{a}_j(n)\dot{\ff}_j     \quad \for \quad  
              0 \leq n < b^m.
\end{equation}
We put
\begin{equation} \nonumber
	\dot{\Psi}_m =\{  	\dot{\by}_n 	  \; | \; 0 \leq n < b^m \}.
\end{equation}	
It is easy to see that $\dot{\Psi}_m =\tilde{\Psi}_m$.

 For $\ddot{\ff}_j=(\ddot{\ff}_j^{(1)},...,\ddot{\ff}_j^{(\dot{s})})$ with
$\ddot{\ff}_j^{(i)} =(\ddot{\ff}_{1,j}^{(i)},...,\ddot{\ff}_{m,j}^{(i)})$, we define
\begin{equation} \nonumber
\ddot{\ff}_j = \dot{\ff}_j  \;\; \for  \;\;   j \in[(\dot{s}-1)d_0\dot{m} + 1, m]
\; \;\ad \;\; 
\ddot{\ff}_j^{(i)} = \dot{\ff}_j^{(i)}  \; \for  \; i \in  [1,\dot{s}-1], \;   j \in[1, m],
\end{equation} 
\begin{equation}\label{6Di14}
 \ddot{\ff}_{j,r}^{(\dot{s})} = \bar{0} \quad\quad \for \quad \quad j \in [1,(\dot{s}-1)d_0\dot{m}],\;
r \in [d^{(\dot{s})}_1,d^{(\dot{s})}_2],  \quad \ad \quad \ddot{\ff}_{j,r}^{(\dot{s})} = \dot{\ff}_{j,r}^{(\dot{s})} 
\end{equation}
for   $j \in [1,(\dot{s}-1)d_0\dot{m}]$
and  $r \in [1,m] \setminus  [d^{(\dot{s})}_1,d^{(\dot{s})}_2]$.
Let
\begin{equation} \label{6Di15}
    		\ddot{\by}^{(i)}_n = 	(\ddot{y}^{(i)}_{n,1},...,\ddot{y}^{(i)}_{n,m}): = 
					\bn(\ddot{\ff}^{(i)}_1,...,\ddot{\ff}^{(i)}_m)=
			\sum_{j=1}^m \bar{a}_j(n)\ddot{\ff}^{(i)}_j
 = \bn  \ddot{\cF}^{(i) \;\top} ,
\end{equation}	
where
$\ddot{\cF}^{(i)} = (\ddot{\ff}_{r,j}^{(i)})_{1 \leq r,j \leq m}$ for $1 \leq i \leq \dot{s}$.
Hence
\begin{equation} \label{6Di16}
	\ddot{\by}_n: = 	(\ddot{\by}^{(1)}_n,...,\ddot{\by}^{(\dot{s})}_n) =\sum_{j=1}^m \bar{a}_j(n)\ddot{\ff}_j     \quad \for \quad                0 \leq n < b^m.
\end{equation}
We put
\begin{equation} \label{6Di17}
	\ddot{\Psi}_m =\{  	\ddot{\by}_n 	  \; | \; 0 \leq n < b^m \}\quad\quad \ad \quad \quad
	 \ddot{Y}_m =\{  \ddot{\by}_n^{({\dot{s}})}
	  |   n \in[0, b^m) \}.
\end{equation}
Now let $\dot{\bx}_n =	(\dot{x}^{(1)}_n,...,\dot{x}^{({\dot{s}})}_n)$ and 
    $\ddot{\bx}_n =	(\ddot{x}^{(1)}_n,...,\ddot{x}^{({\dot{s}})}_n)$, where
\begin{equation} \nonumber
     \dot{x}^{(i)}_n =\sum_{j=1}^m \phi^{-1} (\dot{y}^{(i)}_{n,j})/b^j, \quad \ad \quad
		   \ddot{x}^{(i)}_n =\sum_{j=1}^m \phi^{-1} (\ddot{y}_{n,j}^{(i)})/b^j 
\end{equation} 
for $  1 \leq i \leq \dot{s}$. We have
\begin{equation}\label{6Di19}
\tilde{\Phi}_m =\{ \tilde{\bx}_n \; | \; 0 \leq n < b^m \}=\{ \dot{\bx}_n \; | \; 0 \leq n < b^m \}  \quad \ad \quad 
	 \ddot{Y}_m =\FF_b^m.
\end{equation}

Bearing in mind that $\dot{\ff}_1,...,\dot{\ff}_m$ and $\tilde{\fc}_1,...,\tilde{\fc}_m$
 are two basis of the vector space $\tilde{\Psi}_m$, we get that there  exists a nonsingular 
matrix $B =(b_{j,r})_{1 \leq j,r \leq m}$ with $b_{j,r} \in \FF_b$ such that
$(\dot{\ff}_1,...,\dot{\ff}_m)^{\top} =B(\tilde{\fc}_1,...,\tilde{\fc}_m)^{\top}$. 
Hence
\begin{equation} \nonumber
\dot{\ff}_k =\sum_{r=1}^m b_{k,r} \tilde{\fc}_r,  \quad \ad \quad
\dot{\ff}_{k,j}^{(i)} =\sum_{r=1}^m b_{k,r} \tilde{c}_{r,j}^{(i)},
\end{equation}
for $1 \leq k,j \leq m,  1 \leq i \leq \dot{s}$.
Therefore
\begin{equation}\label{6Di36}
(\dot{\ff}_1^{(i)},...,\dot{\ff}_m^{(i)})^{\top} =B(\tilde{\fc}_1^{(i)},...,\tilde{\fc}_m^{(i)} )^{\top}     \; \ad \;
      \tilde{C}^{(i)} =     \dot{\cF}^{(i)} 
  B^{-1 \;\top}      \; \for \;   i \in [1, \dot{s}].
\end{equation}
Let $n' \in [0,b^m)$, 
$\bn' =(\bar{a}_1(n'),...,\bar{a}_m(n'))$, and let $\bn' = \bn B^{-1} $.\\ 
Using (\ref{6Di01}) and  (\ref{6Di11}),   we get
\begin{equation} \nonumber
  \dot{\by}_{n'}^{(i)}  =\bn' \dot{\cF}^{(i)\top} =\bn' (\dot{\ff}_1^{(i)},...,\dot{\ff}_m^{(i)})^{\top}
 = \bn B^{-1} B (\tilde{\fc}_1^{(i)},...,\tilde{\fc}_m^{(i)} )^{\top} 
\end{equation}
\begin{equation}\nonumber
  =  \bn(\tilde{\fc}_1^{(i)},...,\tilde{\fc}_m^{(i)} )^{\top}   =   \bn \tilde{C}^{(i)\top}
 =  \tilde{\by}_{n}^{(i)},		 \quad \for \quad  1 \leq i \leq \dot{s} \quad \ad \quad
 0 \leq n < b^m.
\end{equation} 
Let $\breve{C}^{(i)} =(\breve{c}_{r,j}^{(i)})_{1 \leq r,j \leq m} : = \ddot{\cF}^{(i)} B^{-1 \;\top}, \;1 \leq i \leq \dot{s}  $,  $ \breve{\fc}^{(i)}_j =(\breve{c}_{1,j}^{(i)},...,\breve{c}_{m,j}^{(i)})$,  
 $1 \leq i \leq \dot{s}, \;1 \leq j \leq m $  and let $\breve{\by}_{n} :=\ddot{\by}_{n'}$, $\breve{\bx}_{n} :=\ddot{\bx}_{n'}$ for  $\bn' = \bn B^{-1} $.
We have
\begin{equation}  \label{6Di43}
  \breve{\by}_{n}^{(i)} =\ddot{\by}_{n'}^{(i)} =\bn' \ddot{\cF}^{(i)\top} =\bn B^{-1}  \ddot{\cF}^{(i)\top} 
 =   \bn \breve{C}^{(i) \top}
  \;\; \for \;\; 1 \leq i \leq \dot{s},\;  0 \leq n < b^m. 
\end{equation}
Hence,  $\breve{C}^{(1)} ,..., \breve{C}^{(\dot{s})}$ are generating matrices of 
the net $ (\breve{\bx}_n)_{0 \leq n < b^m} $.
According to (\ref{6Di14}) and (\ref{6Di36}), we obtain   $\quad \ddot{\cF}^{(i)}=\dot{\cF}^{(i)}$,
\begin{equation}\label{6Di39}
\breve{C}^{(i)}=\tilde{C}^{(i)} \quad \for \quad 1 \leq i \leq \dot{s}-1, 
  \quad \ad \quad \breve{C}^{(\dot{s})}-\tilde{C}^{(\dot{s})} =
	(\ddot{\cF}^{(\dot{s})}-\dot{\cF}^{(\dot{s})}) B^{-1 \; \top}.
\end{equation}

Let $(B^{-1})^{\top} =(\hat{b}_{r,j})_{1 \leq r,j \leq m}$, 
$\Delta c_{r,j} =  \breve{c}^{(\dot{s})}_{r,j} - \tilde{c}^{(\dot{s})}_{r,j}$ and 
$\Delta \ff_{r,j} =  \ddot{\ff}^{(\dot{s})}_{r,j} - \dot{\ff}^{(\dot{s})}_{r,j}$ for $1 \leq r,j \leq m$.
Applying (\ref{6Di11}), (\ref{6Di15}) and (\ref{6Di39}), we derive
\begin{equation}  \label{6Di83}
       \Delta c_{r,j} = \sum_{l=1}^{m}   \Delta \ff_{r,l} \hat{b}_{l,j}  
	  \quad  \for \quad 1 \leq r,j \leq m.
\end{equation}
From (\ref{6Di14}) and (\ref{6Di36}), we get
\begin{equation}\label{6Di84}
   \Delta c_{r,j} =\breve{c}_{r,j}^{(\dot{s})}-\tilde{c}_{r,j}^{(\dot{s})}=0 \;\; \for \;\;
		r \in [(\dot{s}-1)d_0\dot{m} +1,m], \;
		1 \leq j \leq m.
\end{equation}
By (\ref{6Di36}) and (\ref{6Di09}), we have 
\begin{equation}  \label{6Di87}
       \tilde{c}^{(\dot{s})}_{r,j} = \sum_{l=1}^{m}   \dot{\ff}^{(\dot{s})}_{r,l} \hat{b}_{l,j}  
	 = \hat{b}_{r,j}  \quad  \for \quad  r \in [(\dot{s}-1)d_0\dot{m}] +1,m] \;\;  \ad \;\;
		1 \leq j \leq m.
\end{equation}
Using (\ref{6Di04}), we obtain $d^{(\dot{s})}_1 > (\dot{s}-1)d_0\dot{m}$.
By (\ref{6Di14}), (\ref{6Di83}) and (\ref{6Di87}), we get 
\begin{equation}  \label{6Di87a}
       \Delta c_{r,j} = \sum_{l=d^{(\dot{s})}_1}^{d^{(\dot{s})}_2}   \Delta \ff_{r,l} \tilde{c}_{l,j}  
	  \quad \; \for \quad \; r \in [1, (\dot{s}-1)d_0\dot{m}] \; \;\;  \ad \;\;\;
		1 \leq j \leq m.
\end{equation}  \\  


{\bf Lemma 20.}  {\it  With notations as above. Let  $\dot{s} \geq 3$, $ (\tilde{\bx}_n)_{0 \leq n <b^m} $  be a digital 
     $(t,m,{\dot{s}})$ net in base $b$, $\tilde{x}^{\dot{s}}_n \neq 
	\tilde{x}^{\dot{s}}_k $ for $n \neq k$. Then
 $(\breve{\bx}_n)_{0 \leq n <b^m}$ is a digital   $(t,m,\dot{s})$ net
 in base $b$ with $\breve{x}^{\dot{s}}_n \neq \breve{x}^{\dot{s}}_k $ for $n \neq k$,
 \begin{equation}\label{6Di50a}
  \left\| \breve{\bx}_n^{(\dot{s})} \right\|_b 
	= \left\| \tilde{\bx}_n^{(\dot{s})} \right\|_b \quad \for \quad 0 <n<b^m
\end{equation}
and
\begin{equation}\label{6Di50}
    \Lambda  =\FF_b^{\dot{s}d_0\dot{m}},\quad  \for \quad  m \geq 2d_0\dot{s}, \; \dot{m} =[(m-t)/(2d_0({\dot{s}}-1))],\;
\end{equation}
where 
\begin{equation} \nonumber
    \Lambda = \{ (\breve{y}^{(1)}_{n, d^{(1)}_{1}},...,\breve{y}^{(1)}_{n, d^{(1)}_{2}}, ...,
		\breve{y}^{(\dot{s})}_{n, d^{(\dot{s})}_{1}},...,  \breve{y}^{(\dot{s})}_{n, d^{(\dot{s})}_{2}} ) 
		\;|\;   n \in [0, b^{m})  \}
\end{equation}
with $
 d^{(i)}_1 =1, \; d^{(i)}_2 =d_0\dot{m}$ for $1 \leq i <\dot{s}$,
 $d^{(\dot{s})}_1 =m -t+1 -(\dot{s}-1)d_0 \dot{m}$ and 
	$d^{(\dot{s})}_2 =m -t -({\dot{s}}-2)d_0\dot{m}$.} \\

{\bf Proof.} By (\ref{6Di43}),  we have $\breve{\by}_{n} =\ddot{\by}_{n'}$,  $\breve{\bx}_{n} =\ddot{\bx}_{n'}$ and $\tilde{\by}_{n} =\dot{\by}_{n'}$,  $\tilde{\bx}_{n} =\dot{\bx}_{n'}$
 for  $\bn' = \bn B^{-1} $. Hence, in order to prove the lemma, it is sufficient to take $\ddot{\bx}_{n}$ instead of and $\breve{\bx}_{n}$ and $\dot{\bx}_{n}$ instead of  $\tilde{\bx}_{n}$.
 Applying (\ref{6Di17}) and (\ref{6Di19}), we derive that $\ddot{x}^{\dot{s}}_n \neq \ddot{x}^{\dot{s}}_k $ for $n \neq k$.

Suppose that $a_j(n) = 0$ for $1 \leq j \leq (\dot{s} -1)d_0 \dot{m}$.
By (\ref{6Di14}) and (\ref{6Di16}), we get   $ \left\| \ddot{x}^{\dot{s}}_n \right\|_b = \left\| \dot{x}^{\dot{s}}_n \right\|_b$.

Let  $a_j(n) = 0$ for $1 \leq j <j_0 \leq (\dot{s} -1)d_0 \dot{m}$ and let
 $a_{j_0}(n) \neq 0$. 
From (\ref{6Di14}) and (\ref{6Di16}), we have
$ \left\| \ddot{x}^{(\dot{s})}_n \right\|_b = \left\| \dot{x}^{(\dot{s})}_n \right\|_b =b^{-j_0}$.
Hence  $ \left\| \ddot{\bx}_n^{(\dot{s})} \right\|_b = \left\| \dot{\bx}_n^{(\dot{s})} \right\|_b$ for all $n \in [1, b^{m}) $ and (\ref{6Di50a}) follows.\\

Let $\bd =(d_1,...,d_{\dot{s}})$,  $d_i \geq 0 \; (i=1,...,\dot{s})$, 
 $\ddot{\bv}_{\bd}=(\ddot{v}_{1}^{(1)} ,...,\ddot{v}_{d_1}^{(1)},...,\ddot{v}_{1}^{(\dot{s})},...,
\ddot{v}_{d_{\dot{s}}}^{(\dot{s})}) \in \FF_b^{\dot{d} }$, with 
  $\dot{d} =d_1+...+d_{\dot{s}}$, and let
\begin{equation}  \label{6Di52}
   \ddot{U}_{\ddot{\bv}_{\bd}} = \{ 0 \leq n < b^m  \; | \;  \ddot{y}_{n,j}^{(i)} =v_{j}^{(i)}, \;  1\leq j \leq d_i, \;
	   1 \leq i \leq {\dot{s}}\}.
\end{equation}
In order to prove that  $(\ddot{\bx}_n)_{0 \leq n <b^m}$ is a $(t,m,{\dot{s}})$ net, it is sufficient to verify that\\ $\#\ddot{U}_{\ddot{\bv}_{\bd}} =b^{m-\dot{d}}$ for all
  $\ddot{\bv}_{\bd} \in \FF_b^{\dot{d} }$ and all $\bd$ with $\dot{d} \leq m-t$.
By (\ref{6Di11}), (\ref{6Di14}) and (\ref{6Di15}), we get
\begin{equation} \label{6Di52a}
    		\dot{\by}^{(i)}_n = 	
							\sum_{j=1}^m \bar{a}_j(n) \dot{\ff}^{(i)}_j
      \quad \; \; \ad \quad \;\;
	\ddot{\by}^{(i)}_n = 	
							\sum_{j=1}^m \bar{a}_j(n) \ddot{\ff}^{(i)}_j,  \quad \with \quad
							\ddot{\ff}^{(i)}_j =\dot{\ff}^{(i)}_j
\end{equation}
for $1 \leq i \leq {\dot{s}}-1$, $1 \leq j \leq m$ and 
$ i = {\dot{s}}$, $(\dot{s} -1)d_0 \dot{m} +1 \leq j \leq m$, $0 \leq n <b^m$. \\
Hence
\begin{equation} \label{6Di53}
  \dot{\by}^{(i)}_n -	\ddot{\by}^{(i)}_n = 0 \;\; \for \;\; 1 \leq i \leq {\dot{s}}-1,
	 \;  \; \dot{\by}^{(\dot{s})}_n -	\ddot{\by}^{(\dot{s})}_n = 
	\sum_{r=1}^{(\dot{s} -1)d_0 \dot{m}} \bar{a}_r(n)( \dot{\ff}^{(\dot{s})}_r 
	-\ddot{\ff}^{(\dot{s})}_r  )
\end{equation}
 and $ \dot{\by}^{(\dot{s})}_{n,j} -	\ddot{\by}^{(\dot{s})}_{n,j} = 0 $ for
$j \in [1,(\dot{s} -1)d_0 \dot{m}  ]$, $0 \leq n <b^m$.
Let
\begin{equation}  \nonumber
  \dot{v}_{j}^{(i)} := \ddot{v}_{j}^{(i)} \; \for \; 
	  j \in[1, d_i],  i \in [1, {\dot{s}}-1]    \;  \ad \; \dot{v}_{j}^{(\dot{s})} := \ddot{v}_{j}^{(\dot{s})} \; \for \;  j \in [1, \min(d_{\dot{s}},(\dot{s} -1)d_0 \dot{m} )].
\end{equation}
For $d_{\dot{s}} > (\dot{s} -1)d_0 \dot{m}$  and $j \in [(\dot{s} -1)d_0 \dot{m} +1 ,d_{\dot{s}}]$, we define
\begin{equation}  \nonumber
    \dot{v}_{j}^{(\dot{s})} = \ddot{v}_{j}^{(\dot{s})} +
		\sum_{r=1}^{(\dot{s} -1)d_0 \dot{m}} 
		 \ddot{v}_{r}^{(\dot{s})}( \dot{\ff}^{(\dot{s})}_{r,j}   - \ddot{\ff}^{(\dot{s})}_{r,j}  ).
\end{equation}
By (\ref{6Di09}) and (\ref{6Di52a}), we get
\begin{equation}  \nonumber
 \dot{y}_{n,j}^{(\dot{s})} = \dot{v}_{j}^{(\dot{s})}
 \Longleftrightarrow 
 \bar{a}_j(n) = \dot{v}_{j}^{(\dot{s})}= \ddot{v}_{j}^{(\dot{s})}, \quad \for \quad   j \in [1, \min(d_{\dot{s}},(\dot{s} -1)d_0 \dot{m} )],\; n\in [0,b^m).
\end{equation}
Using  (\ref{6Di53}), we obtain for $ n \in[0,b^m)$ that
\begin{equation}  \label{6Di55}
  \ddot{\by}^{(i)}_{n,j}  = \ddot{v}_{j}^{(i)} 
		\Longleftrightarrow 
	 \dot{\by}^{(i)}_{n,j}  = \dot{v}_{j}^{(i)} 
	\quad \for \quad 
	 1 \leq j \leq d_i,\; 1 \leq i \leq {\dot{s}}.
\end{equation}
Let
\begin{equation}  \nonumber
   \dot{U}_{\dot{\bv}_{\bd}} = \{ 0 \leq n < b^m  \; | \;  \dot{y}_{n,j}^{(i)} =\dot{v}_{j}^{(i)}, \;  1  \leq j \leq d_i, \;
	   1 \leq i \leq {\dot{s}}\}
\end{equation}
with  $\dot{\bv}_{\bd}=(\dot{v}_{1}^{(1)} ,...,\dot{v}_{d_1}^{(1)},...,\dot{v}_{1}^{(\dot{s})},..., \dot{v}_{d_{\dot{s}}}^{(\dot{s})})$.\\

Taking into account that  $(\dot{\bx}_n)_{0 \leq n <b^m}$ is a $(t,m,{\dot{s}})$ net in base $b$,
we get from (\ref{6Di52}) and (\ref{6Di55}) that $\#\ddot{U}_{\ddot{\bv}_{\bd}} =
\#\dot{U}_{\dot{\bv}_{\bd}} =b^{m-\dot{d}}$.\\

Now consider the statement (\ref{6Di50}). 
Let  $\ddot{\bv}=(\ddot{v}_{d^{(1)}_1}^{(1)} ,...,\ddot{v}_{d^{(2)}_2}^{(1)},...,
\ddot{v}_{d^{(\dot{s})}_1}^{(\dot{s})},...,
\ddot{v}_{d^{(\dot{s})}_2}^{(\dot{s})}) \in \FF_b^{\dot{d} }$, with
  $\dot{d} =d^{(1)}_2 +... + d^{(\dot{s}-1)}_2+ d^{(\dot{s})}_2 -  d^{(\dot{s})}_1 +1$.
It is easy to see that to obtain (\ref{6Di50}), it is sufficient to verify that
$\ddot{U}_{\ddot{\bv}}^{'} \neq \emptyset $ for all $\ddot{\bv} \in \FF_b^{\dot{d} }$.
where
\begin{equation}  \nonumber
   \ddot{U}_{\ddot{\bv}}^{'} = \{ 0 \leq n < b^m  \; | \;  \ddot{y}_{j}^{(i)} =\ddot{v}_{j}^{(i)}, \;  d^{(i)}_1  \leq j \leq   d^{(i)}_2, \;
	   1 \leq i \leq {\dot{s}}\}.
\end{equation}
According to   (\ref{6Di15}) and (\ref{6Di16}), $\ddot{U}_{\ddot{\bv}}^{'} \neq \emptyset $ if there exists $n \in [0,b^m)$
such that
\begin{equation}\label{6Di58}
     \sum_{r =1}^m   \bar{a}_{r}(n)\ddot{\ff}_{j,r}^{(i)} =\ddot{v}_{j}^{(i)}
		\quad \fall \quad  d^{(i)}_1  \leq j \leq   d^{(i)}_2 \; \; \ad \;\; 1 \leq i \leq \dot{s}.
\end{equation}
By (\ref{6Di09}) and (\ref{6Di14}),  we have that  (\ref{6Di58}) is true 
only if
$\bar{a}_{j} (n)= \ddot{v}_{j}^{(\dot{s})}$ 

for  $ d^{(\dot{s})}_1  \leq j \leq   d^{(\dot{s})}_2 $. 
Let $n_0=  \sum_{j=d^{(\dot{s})}_1}^{d^{(\dot{s})}_2} \phi^{-1}(\ddot{v}_{j}^{(\dot{s})})  b^{j-1}$  and let
\begin{equation}\nonumber
  n= n_0+\sum_{i=1}^{\dot{s}-1}  \sum_{j =d^{(i)}_1}^{d^{(i)}_2} \phi(\ddot{v}_{j}^{(i)}- \ddot{y}_{n_0,j}^{(i)} ) b^{(i-1)d_0\dot{m} +j-1}.
\end{equation}
Therefore $\bar{a}_{j}(n) = \ddot{v}_{j}^{(\dot{s})}$ for 
$j \in [d^{(\dot{s})}_1, d^{(\dot{s})}_2] $ and  $\bar{a}_{(i-1)d_0\dot{m} +j}(n) = \ddot{v}_{j}^{(i)}$ for 
$j \in [d^{(i)}_1, d^{(i)}_2] $, $i \in [1,\dot{s}-1 ]$.
Using (\ref{6Di09}) and (\ref{6Di14}), we get that   (\ref{6Di58}) is true and 
 $\ddot{U}_{\ddot{\bv}}^{'} \neq \emptyset $ for all $\ddot{\bv} \in \FF_b^{\dot{d} }$.
Hence (\ref{6Di50}) is proved, and Lemma 20 follows. \qed \\

{\bf End of the proof of Theorem 6. }	
Let  $C^{(1)},...,C^{(s)} \in \FF_b^{\infty \times \infty}$  be the 
generating matrices of a digital  
 $(t,s)$ 
sequence $(\bx_n)_{n \geq 0}$.
For any $m \in \NN$ 
we denote the $m \times m$ left-upper sub-matrix of
$C^{(i)}$ by $[C^{(i)}]_m$.

 Let $m_k=s^2d_0(2^{2k+2}-1)$, $k=0,1,...$ , 
\begin{equation}  \label{6Di100ab}
   x_n^{(i,k)} =\sum_{j=1}^{m_k} \phi^{-1}(y^{(i,k)}_{n,j})/ b^{j}, \quad
	\by^{(i,k)}_n = \bn [C^{(i) \; \top}]_{m_k}
\end{equation}
and $\by^{(i,k)}_n =(y^{(i,k)}_{n,1},...,y^{(i,k)}_{n,m_k})$ for $n \in [0,b^{m_k})$, $i \in [1,s]$.

For $x =\sum_{j \geq 1}  x_{j}p_i^{-j}$, 
where $x_{i} \in  Z_b =\{0,...,b-1\}$,   we define the truncation 
\begin{equation}  \nonumber
        [x]_m =\sum_{1 \leq j \leq m}  x_{j}b^{-j} \quad \with \quad m \geq 1.
\end{equation}
If $x = (x^{(1)}, . . . , x^{(s)})  \in [0, 1)^s$, then the truncation $[\bx]_m$ is defined coordinatewise, that is, $[\bx]_m = 
( [x^{(1)}]_m, . . . , [x^{(s)}]_m)$.

By  (\ref{Ap301}) - (\ref{Ap303}), we have
\begin{equation}  \label{6Di100abc}
     [\bx_n]_{m_k} = \bx^{(k)}_n :=(x^{(1,k)}_{n},...,x^{(s,k)}_{n})
		\quad \for \quad n \in [0,b^{m_k}).
\end{equation}

Let  $\hat{C}^{(s+1,0)} =(\hat{c}^{(s+1,0)}_{i,j})_{1 \leq i,j \leq m_{0}}$ with $\hat{c}^{(s+1,0)}_{i,j} =\delta_{i,m_0-j+1}$, $i,j=1,...,m_0$. 
%
We will use (\ref{6Di00}) - (\ref{6Di39}) to construct a sequence of matrices 
 $\hat{C}^{(s+1,k)} \in \FF_b^{ m_k \times m_k}$ 
$(k = 1,2,...)$, satisfying  the following induction assumption:\\

For given sequence of matrices 
$\hat{C}^{(s+1,0)},...,\hat{C}^{(s+1,k-1)}$
there exists 
 a matrix \\ $\hat{C}^{(s+1,k)} =(\hat{c}^{(s+1,k)}_{i,j})_{1 \leq i,j \leq m_{k}}$
such that
\begin{equation}  \label{6Di100a}
  \hat{c}^{(s+1,k)}_{m_k-i+1,j}= \hat{c}^{(s+1,k-1)}_{m_{k-1}-i+1,j} \;\; \for \;\;
	 i,j \in [1, m_{k-1}] \quad \ad 	\quad \hat{c}^{(s+1,k)}_{m_k-i+1,j}= 0
\end{equation}
for $i \in[m_{k-1} +1,m_{k}]$, $j \in [1, m_{k-1}]$, 
$(x_n^{(1,k)},..., x_n^{(s,k)},\hat{x}^{(s+1,k)}_n )_{0 \leq n < b^{m_{k}}} $  is a
  $(t,m_{k},s+1~)$ net in base $b$ with 
\begin{equation}  \label{6Di100b}
\hat{x}^{(s+1,k)}_n \neq \hat{x}^{(s+1,k)}_l \;\;\for \; n \neq l \; \ad \;
 \left\| \hat{x}^{(s+1,k)}_n \right\|_b = \left\| n \right\|_b b^{-m_{k}} \; \for \;
  0 \leq n <b^{m_k},
\end{equation}	
where	
\begin{equation}  \label{6Di100}
  \hat{x}^{(s+1,k)}_n =\sum_{j=1}^{m_k} \phi^{-1}(y^{(s+1,k)}_{n,j})/ b^{j}, \quad
	\by^{(s+1,k)}_n = \bn \hat{C}^{(s+1,m_k) \; \top}
\end{equation}
and $\by^{(s+1,k)}_n =(y^{(s+1,k)}_{n,1},...,y^{(s+1,k)}_{n,m_k})$ for $n \in [0,b^{m_k})$.\\

Let $k=1$.   We take 
$\hat{c}^{(s+1,1)}_{i,j} =\delta_{i, m_1-j+1}$ for  $i,j=1,...,m_1$.\\ 
Now  assume  we known $\hat{C}^{(s+1,k)}$ and we want to construct   $\hat{C}^{(s+1,k+1)}$.
We first construct 
$\tilde{C}^{(s+1,k+1)} =(\tilde{c}^{(s+1,k+1)}_{i,j})_{1 \leq i,j \leq m_{k+1}}$
 as following
\begin{equation}  \label{6Di105}
  \tilde{c}^{(s+1,k+1)}_{m_{k+1}-i+1,j}= \hat{c}^{(s+1,k)}_{m_{k}-i+1,j} \quad \for  \quad 
	i,j \in [1,m_{k}], 
	\quad \quad	 \tilde{c}^{(s+1,k+1)}_{i,j}= \delta_{i, m_{k+1}-j+1}
\end{equation}
\begin{equation}   \nonumber 
\for \quad 	i\in [1,m_{k+1}-m_{k}],\; 	j \in [1,m_{k+1}] 
	 \quad  \ad \quad \tilde{c}^{(s+1,k+1)}_{i,j}= \bar{0} 
\end{equation}
for $(i,j) \in [1,m_{k+1}-m_{k}]\times[1,m_{k}]$ and $ (i,j) \in
 [m_{k+1}-m_{k}+1,m_{k+1}]\times [m_{k}+1,m_{k+1}]$.\\

{\bf Lemma 21.}  {\it With notations as above,     $(x_n^{(1,k+1)},..., x_n^{(s,k+1)}, 
\tilde{x}_n^{(s+1,k+1)})_{0 \leq n < b^{m_{k+1}}} $  is a 
  $(t,m_{k+1},s+1)$ net in base $b$  with $\tilde{x}^{(s+1,k+1)}_n \neq \tilde{x}^{(s+1,k+1)}_l $ for $n \neq l$, and 
\begin{equation}  \label{6Di104a}
 \left\| \tilde{x}^{(s+1,k+1)}_n \right\|_b = \left\| n \right\|_b b^{-m_{k+1}} \quad \for \quad
  0 < n <b^{m_{k+1}}.
\end{equation} 	}  \\

{\bf Proof.} 
Let $\bd =(d_1,...,d_{s+1})$,
 $\bv_{\bd}=(v_{1}^{(1)},...,v_{d_1}^{(1)},...,v_{1}^{(s+1)},...,v_{d_{s+1}}^{(s+1)}) \in \FF_b^{\dot{d} }$ 
  with\\ $\dot{d} =d_1+...+d_{s+1}$,
\begin{eqnarray}  \nonumber
 &  \tilde{U}_{\bv_{\bd}} = \{ 0 \leq n < b^{m_{k+1}}  \; | \;  y_{n,j}^{(i,k)} =v_{j}^{(i)}, \quad 1\leq j \leq d_i, \;
	   1 \leq i \leq s \\
	&	\ad \quad \tilde{y}_{n,j}^{(s+1,k+1)} =v_{j}^{(s+1)}, \quad  1\leq j \leq d_{s+1}
		\}. \label{6Di105a}
\end{eqnarray}
In order to prove that  $(x_n^{(1,k+1)},..., x_n^{(s,k+1)}, 
\tilde{x}_n^{(s+1,k+1)})_{0 \leq n < b^{m_{k+1}}} $ is a $(t,m_{k+1},s+1)$ net, it is sufficient to verify that $\#\tilde{U}_{\bv_{\bd}} =b^{m_{k+1}-\dot{d}}$ for all  $\bv_{\bd} \in \FF_b^{\dot{d} }$ and all $\bd$ with $\dot{d} \leq m_{k+1}-t$.\\

Suppose that $d_{s+1} \leq m_{k+1}  -m_{k}$.\\
Let $n \in [0,b^{m_{k+1}})$, $n_0 \equiv n \;(\mod \; b^{m_{k+1}-d_{s+1}})$, $n_0 \in [0,b^{m_{k+1}-d_{s+1}})$ and let $n_1 = n -n_0$. It is easy to see that 
\begin{equation} \nonumber
    \tilde{y}_{n,j}^{(s+1,k+1)} = \tilde{y}_{n_0,j}^{(s+1,k+1)}
		+ \tilde{y}_{n_1,j}^{(s+1,k+1)}. 
\end{equation}
Let  $j \in [1,  m_{k+1} - m_{k}]$. By  (\ref{6Di105}), we get
\begin{equation} \label{6Di106b}
    \tilde{y}_{n,j}^{(s+1,k+1)} 
		= \sum_{r=1}^{m_{k+1}} \bar{a}_{r}(n) \tilde{c}^{(s+1,k+1)}_{j,r}
		= \sum_{r=1}^{m_{k+1}-m_k} \bar{a}_{r}(n) \delta_{j,m_{k+1} +1-r} 
		=		\bar{a}_{m_{k+1} +1-j}(n) .
\end{equation}
Let  $\ddot{n} =\sum_{j=1}^{d_{s+1}} \phi(v^{(s+1)}_j)b^{m_{k+1}-j}$.
By (\ref{6Di105a}), we get $ n \in  \tilde{U}_{\bv_{\bd}}  \Leftrightarrow n_1=\ddot{n}$ and \\
$ n_0 \in   \tilde{U}_{\bv_{\bd}}^{'}$, where
\begin{equation}  \nonumber
  \tilde{U}_{\bv_{\bd}}^{'}= \{ 0 \leq \dot{n} < b^{m_{k+1}-d_{s+1}}\; | \;   y_{\dot{n},j}^{(i,k+1)} =v_{j}^{(i)} -  y_{\ddot{n},j}^{(i,k+1)}, \;  1  \leq j \in [1, d_i], i \in [1,s] \}.	 
\end{equation}
Bearing in mind (\ref{6Di100}), (\ref{6Di105}), (\ref{6Di105a}) and that 
$(\bx(n))_{0 \leq n <b^{m_{k+1}-d_{s+1}}}$ is a 
$(t,m_{k+1} -d_{s+1},s)$ net in base $b$ , we obtain $\#\tilde{U}_{\bv_{\bd}} =
 \# \tilde{U}_{\bv_{\bd}}^{'} =b^{m_{k+1}-\dot{d}}$.\\

Now let $d_{s+1} > m_{k+1}  -m_{k}$. 
Let $n \in [0,b^{m_{k+1}})$, $n_0 \equiv n \;(\mod \;  b^{m_k})$, $n_0 \in [0,b^{m_k})$ and let $n_1 = n -n_0$. We have 
\begin{equation} \nonumber
    \tilde{y}_{n,j}^{(s+1,k+1)} = \tilde{y}_{n_0,j}^{(s+1,k+1)}
		+ \tilde{y}_{n_1,j}^{(s+1,k+1)}. 
\end{equation}
Let  $\ddot{n} =\sum_{j=1}^{m_{k+1}-m_k} \phi(v^{(s+1)}_j)b^{m_{k+1}-j}$.
By (\ref{6Di105a}) and (\ref{6Di106b}), we get
\begin{eqnarray}  \nonumber
 n \in  \tilde{U}_{\bv_{\bd}}  \Leftrightarrow n_1=\ddot{n}\;\;  \ad \;\; n_0 \in \{ 0 \leq \dot{n} < b^{m_k}  \; | \;   y_{\dot{n},j}^{(i,k+1)} =v_{j}^{(i)} -  y_{\ddot{n},j}^{(i,k+1)}, \;  1  \leq j \leq d_i, 	    \\ 
		 1 \leq i \leq s \;\; \ad \quad y_{\dot{n},j}^{(s+1,k+1)} =v_{j}^{(s+1)}-y_{\ddot{n},j}^{(s+1,k+1)}, \quad  m_{k+1}  -m_{k}+1 \leq j \leq d_{s+1}\}. \nonumber
\end{eqnarray}
Let  $j \in [ m_{k+1} - m_{k}+1, m_{k+1}]$ and  let $j_0 = m_{k+1} +1-j \in [1,m_{k}]$.\\
 By  (\ref{6Di105}), we derive
\begin{equation}  \nonumber 
   \tilde{y}_{\dot{n},j}^{(s+1,k+1)} =\tilde{y}_{\dot{n},m_{k+1} +1-j_0}^{(s+1,k+1)}
		 = \sum_{r=1}^{m_{k+1}} \bar{a}_{r}(\dot{n}) \tilde{c}^{(s+1,k+1)}_{m_{k+1} +1-j_0,r}	
		 = \sum_{r=1}^{m_{k}} \bar{a}_{r}(\dot{n}) \tilde{c}^{(s+1,k+1)}_{m_{k+1} +1-j_0,r} 
\end{equation}
\begin{equation} \label{6Di106}
		= \sum_{r=1}^{m_{k}} \bar{a}_{r}(\dot{n}) \tilde{c}^{(s+1,k)}_{m_{k}+1-j_0,r}
		=\tilde{y}_{\dot{n},m_{k} +1-j_0}^{(s+1,k)} \quad \fall \quad  \dot{n} \in [0, b^{m_k}) . 
\end{equation}
We have that $y_{\dot{n},j}^{(i,k+1)}=y_{\dot{n},j}^{(i,k)}$ $(i=1,....,s)$ and  
$y_{\dot{n},j}^{(s+1,k+1)}=y_{\dot{n},m_k+1 -j_0}^{(s+1,k)}$ for $ \dot{n} \in [0, b^{m_k})$.
Hence
\begin{equation}  \nonumber
  n \in  \tilde{U}_{\bv_{\bd}}  \Leftrightarrow 
	n_1=\ddot{n}\;  \ad \; n_0 \in \tilde{U}_{\bv_{\bd}}^{'} = \Big\{ 0 \leq \dot{n} < b^{m_k} \;   | \;   y_{\dot{n},j}^{(i,k)} =v_{j}^{(i)} -  y_{\ddot{n},j}^{(i,k+1)},   j \in [1, d_i], 	 
\end{equation}
\begin{equation}  \nonumber
	 i \in [1, s], \;\;	\ad \quad y_{\dot{n},j- m_{k+1} + m_k}^{(s+1,k)} =v_{j-m_{k+1}+m_k}^{(s+1)}- y_{\ddot{n},j}^{(s+1,k+1)}, \;\;  j \in ( m_{k+1} - m_k, d_{s+1}] \Big\}. \nonumber
\end{equation}

Taking into account that $(x_n^{(1,k)},..., x_n^{(s,k)},\tilde{x}^{(s+1,k)}_n))_{0 \leq n <b^{m_{k}}}$ is a 
$(t,m_k,s+1)$ net in base $b$, we obtain $\# \tilde{U}_{\bv_{\bd}}=
 \# \tilde{U}_{\bv_{\bd}}^{'} =b^{m_k -(\dot{d} -m_{k+1}+m_k)}=b^{m_{k+1}-\dot{d}}$.
Therefore 
$(x_n^{(1,k+1)},..., x_n^{(s,k+1)}, 
\tilde{x}_n^{(s+1,k+1)})_{0 \leq n < b^{m_{k+1}}} $ is a $(t,m_{k+1},s+1)$ net in base $b$.

From (\ref{6Di105}), (\ref{6Di106b}), (\ref{6Di106}) and the induction assumption, we get
 that \\
 $\tilde{x}^{(s+1,k+1)}_n \neq \tilde{x}^{(s+1,k+1)}_l $ for $n \neq l$. \\

Consider the  assertion (\ref{6Di104a}). Let $n \in [0,b^{m_{k+1}})$ and let   
\begin{equation} \label{6Di106d}
  \left\| \tilde{x}^{(s+1,k+1)}_n \right\|_b = b^{-j_1}. 
\end{equation}
Hence $\tilde{y}^{(s+1,k+1)}_{n,j} =0$ for $1 \leq j \leq j_1 -1$ and 
$\tilde{y}^{(s+1,k+1)}_{n,j_1} \neq 0$ (see (\ref{In04})).

Let  $j_1 \in [1,  m_{k+1} - m_{k}]$. By  (\ref{6Di106b}), we get 
$\bar{a}_{m_{k+1} +1-j}(n) =0$ for $1 \leq j \leq j_1 -1$ and 
$\bar{a}_{m_{k+1} +1-j_1}(n) \neq 0$. Therefore 
$\left\| n \right\|_b  = \left\|  \sum_{i=1}^{m_{k+1}} a_i(n) b^{i-1}   \right\|_b =b^{m_{k+1}-j_1} $.

Now let $j_1 \in [  m_{k+1} - m_{k} +1,  m_{k+1} ]$. From  (\ref{6Di106b}), we obtain  $\bar{a}_{m_{k+1} +1-j}(n) =0$ for $1 \leq j \leq  m_{k+1} - m_{k}$. Hence  $n \in [0,b^{m_{k}})$.
Using  (\ref{6Di105}) and (\ref{6Di106b}), we have 
$\tilde{y}_{n,j}^{(s+1,k)} =\tilde{y}_{n,j -m_{k+1} + m_{k}}^{(s+1,k)} 
$ for $  m_{k+1} - m_{k}+1 \leq j \leq j_1$. 
Therefore $\tilde{y}^{(s+1,k)}_{n,j} =0$ for $1 \leq j \leq j_1-m_{k+1} + m_{k} -1 $ and 
$\tilde{y}^{(s+1,k)}_{n,j_1-m_{k+1} + m_{k}} \neq 0$.
Using the induction assumption  (\ref{6Di100b}), we get 
$ b^{-j_1+m_{k+1} - m_{k}}= \left\| \tilde{x}^{(s+1,k)}_{n} \right\|_b =
\left\| n\right\|_b b^{-m_k}$.

By (\ref{6Di106d}), we obtain $  \left\| \tilde{x}^{(s+1,k+1)}_{n} \right\|_b =
\left\| n\right\|_b b^{-m_{k+1}}$.
Thus  assertion (\ref{6Di104a}) is proved and Lemma 21 follows. \qed \\

Now we apply (\ref{6Di00}) - (\ref{6Di39}) with $\dot{s} =s+1$, $m=m_{k+1}$,
$\tilde{C}^{(i)} :=[C^{(i)}]_{m_{k+1}}$ $(i=1,...,s)$ and $\tilde{C}^{(s+1)}:=\tilde{C}^{(s+1,k+1)}$ 
 to construct  matrices $\breve{C}^{(i)}$ $(i=1,...,s+1)$. 

From (\ref{6Di39}),  we have
\begin{equation}  \label{6Di107}
    \breve{C}^{(i)} = \tilde{C}^{(i)} =[C^{(i)}]_{m_{k+1}} \quad \for \quad i=1,...,s.
\end{equation} 
Let $ \hat{C}^{(s+1,k+1)} := \breve{C}^{(s+1)}$.
According to (\ref{6Di84}) and (\ref{6Di105}),  we get
\begin{equation}  \label{6Di107a}
   \hat{c}^{(s+1,k+1)}_{r,j} -\tilde{c}^{(s+1,k+1)}_{r,j}=0
	\quad \for \quad r \in [s d_0 \dot{m}_{k+1} +1, m_{k+1}] \; \ad \;1 \leq j \leq  m_{k+1}.
\end{equation}
By (\ref{6Di04}) and (\ref{6Di87a}), we obtain for $ r \in [1,s d_0\dot{m}_{k+1}]$ and $1 \leq j \leq m_{k+1}$
\begin{equation}  \label{6Di107b}
      \hat{c}^{(s+1,k+1)}_{r,j} -\tilde{c}^{(s+1,k+1)}_{r,j} = 
			\sum_{l=d^{(s+1,k+1)}_1}^{d^{(s+1,k+1)}_2}   \Delta \ff_{r,l}^{(s+1,k+1)} 
			\tilde{c}^{(s+1,k+1)}_{l,j},
\end{equation}
where $d^{(s+1,k+1)}_1 =m_{k+1} -t+1 -sd_0 \dot{m}_{k+1}$,
	$d^{(s+1,k+1)}_2 =m_{k+1} -t -(s-1)d_0\dot{m}_{k+1}$, $m_{k+1}=s^2d_0(2^{2k+4}-1)$, $d_0=d+t$ and 
	$\dot{m}_{k+1} =[(m_{k+1}-t)/(2sd_0)]$. \\
We have $d^{(s+1,k+1)}_1 >(s-1)d_0\dot{m}_{k+1}  $,	$\dot{m}_{k+1}   =2^{2k+3}-1$ for $k=0,1,...$ and
\begin{equation}  \nonumber
 m_{k+1}-d^{(s+1,k+1)}_2  \geq (s-1)d_0\dot{m}_{k+1} \geq 2^{-1}s^2d_0(2^{2k+3}-1) > m_k.
\end{equation}
By (\ref{6Di105}), we obtain $\tilde{c}^{(s+1,k+1)}_{r,j} =0$ for $ r \leq d^{(s+1,k+1)}_2 <m_{k+1} -m_k$ and
 $1 \leq j \leq m_k$.\\
From (\ref{6Di107b}), we derive
\begin{equation}  \label{6Di107c}
      \hat{c}^{(s+1,k+1)}_{r,j} -\tilde{c}^{(s+1,k+1)}_{r,j} = 0
				 \;\;\; \for \;\;\; r \in [1, sd_0\dot{m}_{k+1}]  \;\;  \ad \;\;
		1 \leq j \leq m_{k}.
\end{equation}
Bearing in mind that
\begin{equation}  \nonumber
 m_{k+1} -sd_0 \dot{m}_{k+1} =s^2d_0(2^{2k+4}-1) - s^2d_0(2^{2k+3}-1)=s^2d_02^{2k+3}
>m_k,
\end{equation}
we get from (\ref{6Di107a}) and (\ref{6Di105})
\begin{equation}  \label{6Di124a}
  \hat{c}^{(s+1,k+1)}_{m_{k+1} -i+1,j} =\tilde{c}^{(s+1,k+1)}_{m_{k+1}-i+1,j}=\hat{c}^{(s+1,k)}_{m_{k} -i+1,j}
			\quad \for \quad 1 \leq i,j \leq  m_{k}.
\end{equation}
Applying (\ref{6Di105}), (\ref{6Di107a}) and (\ref{6Di107c}), we have
\begin{equation}  \nonumber
\hat{c}^{(s+1,k+1)}_{i,j}=\tilde{c}^{(s+1,k+1)}_{i,j}=0,  \; \for \; 1 \leq i \leq  m_{k+1} - m_{k}, \; 1 \leq j \leq m_{k}
\end{equation}
Now using (\ref{6Di124a}), we obtain (\ref{6Di100a}).\\
We see that (\ref{6Di100b}) follows from (\ref{6Di104a}) and (\ref{6Di50a}).
Consider the net $(\hat{\bx}_n^{(k+1)})_{ n=0}^{ b^{m_{k+1}}-1} $ with 
 $\hat{\bx}_n^{(k+1)} =(x_n^{(1,k+1)},..., x_n^{(s,k+1)}, \hat{x}_n^{(s+1,k+1)}) :=\breve{\bx}_n=(\breve{x}_n^{(1)},...,  \breve{x}_n^{(s+1)}) $. 
Let
\begin{equation} \nonumber
    \Lambda_{k+1} = \Big\{ \Big( \big(y^{(i,{k+1})}_{n, 1},...,y^{(i,{k+1})}_{n, d^{(i,k+1)}}\big)_{1 \leq i \leq s}, \;
		\hat{y}^{(s+1,{k+1})}_{n, d^{(s+1,k+1)}_1},...,  
		\hat{y}^{(s+1,{k+1})}_{n, d^{(s+1,k+1)}_2}  \Big)   
	\;	\Big| \;   n \in [0, b^{m_{k+1}})  \Big\}
\end{equation}
with   $  d^{(i,k+1)} =d_0\dot{m}_{k+1}$ for $1 \leq i  \leq s$.
Using (\ref{6Di04}), (\ref{6Di107}) and Lemma 20, we obtain
\begin{equation}\label{6Di110}
    \Lambda_{k+1}  =\FF_b^{(s+1)d_0\dot{m}_{k+1}},\quad  \for \quad \dot{m}_{k+1} =[(m_{k+1}-t)/(2sd_0)]  =s(2^{k+1}-1),
\end{equation}
and $(\hat{\bx}_n^{(k+1)})_{0 \leq n < b^{m_{k+1}}} $ is a $(t,m_{k+1},s+1)$ net in base $b$. 
Thus we have that   $\hat{C}^{(s+1,k+1)}$  satisfy the induction assumption. \\

Let
$C^{(s+1,k+1)} =(c^{(s+1,k+1)}_{i,j})_{1 \leq i,j \leq m_{k+1}}$
where $ c^{(s+1,k+1)}_{i,j}:=\hat{c}^{(s+1,k+1)}_{m_{k+1}-i+1,j}$ for\\ $1 \leq i,j \leq m_{k+1}$.
By (\ref{6Di100a}), we get
\begin{equation}  \label{6Di120}
  [C^{(s+1,k+1)}]_{m_{k}}= C^{(s+1,k)} \; \ad \; c^{(s+1,k+1)}_{i,j}= 0, \; 
	i \in (m_{k}, m_{k+1} ] , \; j \in [1, m_{k}].
\end{equation}
Now let  $C^{(s+1)} = (c^{(s+1)}_{i,j})_{ i,j \geq 1} = \lim_{k \to \infty}C^{(s+1,k)}$ i.e.
$[C^{(s+1)}]_{m_{k}}: =C^{(s+1,k)}$, $k=1,2,...$ . 
We define 
\begin{equation}  \label{6Di123a} 
 h_k(n):=h_{k,1}(n)+...+ h_{k,m_k}(n)b^{m_k-1}:= \hat{x}_n^{(s+1,k)} b^{m_k}\quad \for \quad
0 \leq n < b^{m_{k}}.
\end{equation}
From  (\ref{6Di100}), we have
\begin{eqnarray}  \nonumber
   &  \phi(h_{k,i}(n)) = \phi(\hat{x}_{n,m_k-i+1}^{(s+1,k)}) =  \hat{y}_{n,m_k-i+1}^{(s+1,k)}
		 = \sum_{j=1}^{m_k} \bar{a}_j(n)\hat{c}^{(s+1,k)}_{m_k-i+1,j} \\
	&	=\sum_{j=1}^{m_k} \bar{a}_j(n) c^{(s+1,k)}_{m_k-i+1,j}
\quad \for \quad  0 \leq n < b^{m_{k}}.   \label{6Di122} 
\end{eqnarray}
Applying (\ref{6Di120}), we obtain for $n \in [0,b^{m_{k}})$ that
\begin{equation}  \label{6Di123} 
   h_{k,i}(n)=0 \; \for \; i >m_k \; \ad \;
    h_{k}(n)=h_{k-1}(n) \in [0,b^{m_{k-1}}) \;\for \;n \in [0,b^{m_{k-1}}).
\end{equation}
For $n \in [1, b^{m_k})$, we  get from (\ref{6Di122}) and (\ref{6Di100b}) that
\begin{equation}   \label{6Di124}
   \left\| h_k(n) \right\|_b = \left\| n \right\|_b.	 
\end{equation}
Let $l \neq n \in [0,b^{m_{k}})$. Using (\ref{6Di100b}), we have 
 $(\hat{y}_{l,1}^{(s+1,k)},...,\hat{y}_{l,m_k}^{(s+1,k)}) \neq \\(\hat{y}_{n,1}^{(s+1,k)},...,\hat{y}_{n,m_k}^{(s+1,k)})  $. Hence  $(h_{k,1}(l),...,h_{k,m_k}(l)) \neq (h_{k,1}(n),...,h_{k,m_k}(n))$ and 
$h_{k}(l) \neq h_{k}(n)$. 

Therefore $h_k$  is a bijection from $[0,b^{m_k})$ to $[0,b^{m_k})$. 
We define $h_{k}^{-1}(n)$ such that $  h_{k}(h_{k}^{-1}(n)) =n $ for all $n \in [0,b^{m_{k}}).$\\
Let $n \in [0,b^{m_{k}})$ and $l =h_{k}^{-1}(n)$, then $l \in [0,b^{m_{k}})$ and
$h_{k+1}(l) =h_{k}(l) =n$. Thus 
\begin{equation}  \label{6Di125} 
   h_{k+1}^{-1}(n) = h_{k}^{-1}(n)=l \quad \for  \quad n \in [0,b^{m_{k}}).
\end{equation}

Let  $h(n) = \lim_{k \to \infty}h_k(n)$, and $h^{-1}(n) = \lim_{k \to \infty}h_k^{-1}(n)$. \\
%
%
Let $n \in [0,b^{m_{k}})$ and let $l =h_{k}^{-1}(n)$.
By (\ref{6Di123})  and (\ref{6Di125}),  we get 
\begin{equation} \nonumber
   h(n) =h_k(n) =l, \quad  h^{-1}(l) =h_{k}^{-1}(l)=n, \quad \ad  \quad h^{-1}(h(n))=n.
\end{equation}


Consider the $d-$admissible property of the sequence  $(\bx_{h^{-1}(n)})_{n \geq 0} $.  It is sufficient to take  $k=0$ in (\ref{In04}).

Let $n \in [0,b^{m_{k}})$.  
By (\ref{6Di124}),  we have $\left\| h(n) \right\|_b = \left\| h_k(n) \right\|_b=\left\| n \right\|_b$.
Taking into account Definition 5 and that $(\bx_n)_{n \geq 0} $ is a $d-$admissible sequence, we obtain
\begin{equation} \label{6Di126a} 
   \left\| n \right\|_b \left\|\bx_{h^{-1}(n)} \right\|_b 
	=\left\| h(l) \right\|_b \left\|\bx_{l} \right\|_b
=\left\| l \right\|_b \left\|\bx_{l} \right\|_b \geq b^{-d}, \quad \with \quad  l =h^{-1}(n). 
\end{equation}
Hence $(\bx_{h^{-1}(n)})_{n \geq 0} $ is a $d-$admissible sequence.\\


By  the induction assumption,  $([\bx_n]_{m_k}, h_k(n)/b^{m_k})_{0 \leq n <b^{m_k}}$ is a $(t,m_k,s+1)$ net 
in base $b$ for $k \geq 1$.
Hence  $(\bx_n, h(n)/b^{m_k})_{0 \leq n <b^{m_k}}$ and 
 $( \bx_{h^{-1}(n)}, n/b^{m_k})_{0 \leq n <b^{m_k}}$
 are also   $(t,m_k,s+1)$ nets in base $b$ for $k \geq 1$.
By Lemma 1, $(\bx_{h^{-1}(n)})_{n \geq 0} $ is a $(t,s)$ sequence in base $b$.

Let $N \in [b^{m_k}, b^{m_{k+1}})$. Applying Lemma B, we get
\begin{eqnarray}  \nonumber
   \sigma:= 1+  \min_{0 \leq Q <b^{m_k},  \bw \in E_{m_k}^{s}} \max_{1 \leq M \leq N} 
		M \emph{D}^{*}((\bx_{h^{-1}(n\ominus Q)} \oplus \bw)_{0 \leq  n < M})  \nonumber\\
		\geq 1+  \min_{0 \leq Q <b^{m_k},  \bw \in E_{m_k}^{s}} \max_{1 \leq M \leq b^{m_k}} 
		M \emph{D}^{*}((\bx_{h^{-1}(n\ominus Q)} \oplus \bw)_{0 \leq  n < M})  \nonumber\\
  \geq  \min_{0 \leq Q <b^{m_k},  \bw \in E_{m_k}^{s}} 
		b^{m_k} \emph{D}^{*}((\bx_{h^{-1}(n\ominus Q)} \oplus \bw,n/b^{m_k})_{0 \leq  n < b^{m_k}})         \nonumber \\
		   \geq   \min_{0 \leq Q <b^{m_k},  \bw \in E_{m_k}^{s}} 
     		b^{m_k} \emph{D}^{*}((\bx_{l} \oplus \bw, h(l)\oplus Q/b^{m_k})_{0 \leq  l < b^{m_k}})	
         \nonumber
\end{eqnarray}
where $l=h^{-1}(n\ominus Q)$ and $n =h(l)\oplus Q$.
Bearing in mind that $h(n)= h_k(n)$ for $0 \leq  n < b^{m_k} $, and that 
$\hat{x}_n^{(s+1,k)} =h_k(n)/b^{m_k}$  for $0 \leq  n < b^{m_k} $, we get
\begin{equation}  \label{6Di130}
   \sigma \geq  \min_{0 \leq Q <b^{m_k},  \bw \in E_{m_k}^{s}} 	
	b^{m_k} \emph{D}^{*}((\bx_{n}\oplus \bw,\hat{x}_n^{(s+1,k)}\oplus (Q/b^{m_k}))_{0 \leq  n < b^{m_k}})	.
\end{equation}

By (\ref{6Di126a}) and (\ref{In04}), we obtain that $(\bx_{n} , h(n)/b^{m_k})_{0 \leq  n < b^{m_k}}$ is a $d$-admissible net.

Applying (\ref{6Di100abc}) and the induction assumption, we get that $(\bx_{n} , h(n)/b^{m_k})_{0 \leq  n < b^{m_k}}$ is a $(t,m_k,s+1)$ net in base $b$.
Let
\begin{equation} \nonumber
    \Lambda_{k}^{'} = \Big\{ \Big( \big(y^{(i)}_{n,1},...,y^{(i)}_{n, d^{(i,k)}}\big)_{1 \leq i \leq s}, \;
		\hat{y}^{(s+1,k)}_{n, d^{(s+1,k)}_1},...,  \hat{y}^{(s+1,k)}_{n, d^{(s+1,k)}_2}  \Big)   
	\;	\Big| \;   n \in [0, b^{m_k})  \Big\}.
\end{equation}
Using (\ref{6Di100ab}), (\ref{6Di100abc}) and (\ref{6Di123a})    , we obtain $y^{(i)}_{n, j} =y^{(i,k)}_{n, j}$ for $1 \leq j \leq m_k$, $1 \leq i \leq s$, and $ h(n)/b^{m_k} =\hat{x}^{(s+1,k)}_{n} $.
By (\ref{6Di110}), we have
\begin{equation}\nonumber
   \Lambda_{k}^{'} = \Lambda_{k}  =\FF_b^{(s+1)d_0\dot{m}},\quad  \for \quad \dot{m} =[(m_k-t)/(2sd_0)]  =d^{(s+1,k)}_2  - d^{(s+1,k)}_1 +1.
\end{equation}

Now we apply Corollary 2 with $\dot{s} =s+1$, $\epsilon =(2sd_0)^{-1}$ , $\eta=\hat{e} =1$,
 $\tilde{r}=t$, $m=m_k$, $\tilde{m}=m-t$, $\ddot{m}_{s+1} =d^{(s+1,k)}_1 -1$,   $B_i = \emptyset$ for 
$i \in [1,s+1]$, and $B=0$.
Taking into account (\ref{6Di130}), we get the assertion in Theorem 6. \qed



\begin{thebibliography} {20}

\bibitem [Be]{Be}
  Beck, J., Probabilistic Diophantine approximation. I. Kronecker sequences, 
           Ann. of Math. (2) 140 (1994), no. 1, 109-160.

\bibitem [BC]{BC}
 Beck, J., Chen, W.~W.~L., 
 Irregularities of Distribution,  
 Cambridge Univ. Press,  Cambridge, 1987.


\bibitem [Bi]{Bi}
Bilyk, D.,  On Roth's orthogonal function method in discrepancy theory, Unif. Distrib. Theory 6 (2011), no. 1, 143-184.



\bibitem
[DiPi]{DiPi}   Dick, J. and Pillichshammer, F.,  Digital Nets and Sequences, Discrepancy Theory and Quasi-Monte Carlo Integration, Cambridge
University Press, Cambridge, 2010.


\bibitem
[DiNi]{DiNi}   Dick, J. and  Niederreiter, H., Duality for digital sequences,
Journal of Complexity ,  25 (2009),  406-414.
\bibitem
[FaCh]{FaCh} Faure, H. and Chaix, H., Lower bound for discrepancy in two dimensions, Acta Arith. 76 (1996), no. 2, 149-164. 

\bibitem
[KrLaPi]{KrLaPi}  Kritzer, P., 
 Larcher, G. and  Pillichshammer, F.,
Discrepancy estimates for index-transformed uniformly
distributed sequences, arXiv:1407.8287


\bibitem
[LaPi]{LaPi}
 Larcher, G. and  Pillichshammer, F.,
A metrical lower bound on the star discrepancy of
digital sequences, 
Monat Math., 174 (2014), 105-123.



\bibitem
[Le1] {Le1} Levin, M.B.,  Adelic constructions of low discrepancy sequences, Online J. Anal. Comb. No. 5 (2010), 27 pp. 

\bibitem
[Le2] {Le2} Levin, M.B.,  On the lower bound in the  lattice point remainder problem  for a parallelepiped, arxiv: 1307.2080.


\bibitem
[Le3] {Le3} Levin, M.B.,  On the lower bound of the discrepancy  of Halton's sequences: I, 
  arXiv:1412.8705

	
\bibitem
[Le4] {Le4} Levin, M.B.,  On the lower bound of the discrepancy  of $(t,s)$ sequences:~I, 
 arXiv:1505.06610

\bibitem
[Le5] {Le5} Levin, M.B.,  On the lower bound of the discrepancy  of $(t,s)$ sequences: III,
 Admissible lattices,  in preparation.

\bibitem
[LiNi]{LiNi} 
 Lidl, R., and  Niederreiter, H., Introduction to finite fields and their applications.
Cambridge University Press, Cambridge, first edition, 1994.

\bibitem
[MaNi]{MaNi} 
Mayor, D.J.S.  and Niederreiter, H., 
A new construction of $(t,s)$-sequences and some improved bounds on their quality parameter,  
Acta Arith. 128 (2007), no. 2, 177-191.

\bibitem
[Ni]{Ni}  Niederreiter, H., Random Number Generation and Quasi-Monte Carlo Methods, in: CBMS-NSF Regional Conference Series
in Applied Mathematics, vol. 63, SIAM, 1992.

\bibitem
[NiXi]{NiXi}
 Niederreiter, H. and  Xing. C.P., Low-discrepancy sequences and global function
fields with many rational places, Finite Fields Appl. 2 (1996), 241-273.


\bibitem
[NiPi]{NiPi}
 Niederreiter, H. and  Pirsic, G., Duality for digital nets and its applications,
Acta Arith. 97 (2001), 173-182.


\bibitem
[NiYe]{NiYe}   Niederreiter, H. and  Yeo, A.S.,  Halton-type sequences from global function fields,
Sci. China Math. 56 (2013), 1467-1476.



\bibitem
[Sa]{Sa}  Salvador, G.D.V., Topics in the Theory of Algebraic Function Fields. Mathematics: Theory $\&$ Applications. Birkhauser Boston, Inc., Boston, MA, 2006.

\bibitem
[Skr]{Skr}  Skriganov, M.M., Coding theory and uniform distributions, Algebra i Analiz,
13 (2001), 191-239,  translation in St. Petersburg Math. J. 13 (2002), no. 2,
301-337.

\bibitem
[St]{St} Stichtenoth, H. Algebraic Function Fields and Codes, 2nd ed. Berlin: Springer, 2009.

\bibitem
[Te1] {Te1}  Tezuka, S., Polynomial arithmetic analogue of Halton sequences. ACM Trans Modeling Computer Simulation,
3 (1993),  99-107

\bibitem
[Te2] {Te2}  Tezuka, S.,  Uniform Random Numbers: Theory and Practice. Kluwer International
Series in Engineering and Computer Science. Kluwer, Boston, 1995.

\bibitem
 [Te3] {Te3} Tezuka, S., On the discrepancy of generalized Niederreiter sequences,
Journal of Complexity 29 (2013), 240-247.
\end{thebibliography}
\end{document}